\pdfoutput=1

\newcommand{\pointsize}{11pt}

\documentclass[oneside, \pointsize]{amsbook}
\pdfoutput=1
\usepackage[
   includehead,
   includefoot,
     left = 1.5in, 
      top = 0.5in, 
    right = 1in,
   bottom = 1in
]{geometry}
\usepackage{fancyhdr}
\usepackage{setspace}
\usepackage{calc}
\usepackage{rotating}

\fancyheadoffset[R]{0.5in} 

\fancypagestyle{prelim}{%    
   \renewcommand{\headrulewidth}{0pt} 
   \fancyhf{}           
   \pagenumbering{roman}    
   \cfoot{-\thepage-}       
}

\fancypagestyle{maintext}{%
   \renewcommand{\headrulewidth}{0.4pt}
   \pagenumbering{arabic}
   \fancyhf{}
   \fancyhead[L]{\rightmark}
   \rhead{\thepage}
}

\numberwithin{figure}{chapter} 
\numberwithin{table}{chapter}
\numberwithin{equation}{chapter}
\numberwithin{section}{chapter}

\usepackage{color}
\usepackage[hyphens]{url}
\usepackage[square,numbers,comma,sort&compress]{natbib}
\usepackage{enumerate}
\usepackage{colortbl}
\usepackage{booktabs}
\usepackage{iba-algo}
\usepackage{ifpdf}
\usepackage{graphics}
\usepackage{epsfig}
\usepackage{subfigure}
\usepackage{booktabs}
\usepackage{mathrsfs}  
\usepackage{ifthen}
\usepackage{amsthm}
\usepackage[breaklinks=true,colorlinks,citecolor=blue]{hyperref}
\def\ve#1{\mathchoice{\mbox{\boldmath$\displaystyle\bf#1$}}
{\mbox{\boldmath$\textstyle\bf#1$}}
{\mbox{\boldmath$\scriptstyle\bf#1$}}
{\mbox{\boldmath$\scriptscriptstyle\bf#1$}}}
\newcommand\vealpha{{\boldsymbol{\alpha}}}

\newcommand\velambda{{\boldsymbol{\lambda}}}

\newcommand\U{U}
\newcommand\M{\mathcal M}   %%%  Changed \mathbf to \mathbb   - David
\newcommand\Z{\mathbb Z}   %%%  Changed \mathbf to \mathbb   - David
\newcommand\N{\mathbb N}   %%%  Changed \mathbf to \mathbb   - David
\newcommand\R{\mathbb R}   %%%  Changed \mathbf to \mathbb   - David
\newcommand\F{\mathcal F}   %%%  Changed \mathbf to \mathbb   - David
\newcommand\Q{\mathbb Q}   %%%  Changed \mathbf to \mathbb   - David
\newcommand\Po{\mathcal P}   %%%  Added by David. This is the notation for polytope P.
\newcommand\Co{\mathcal C}   %%%  Added by David. This is the notation for a cone C.
\newcommand\Tr{\mathcal T}   %%%  Added by David. This is the notation for a triangulation.
\newcommand\In{\mathcal I}  %%%  Added by David. This is the notation for independent sets.
\newcommand\B{\mathcal B}   %%%  Added by David. 
\newcommand\Backslash{\, \backslash \,}
\newcommand\Side[1]{\begin{sideways}{\small #1}\end{sideways}}
\DeclareMathOperator{\Inc}{Inc}
\DeclareMathOperator{\argmin}{argmin}
\DeclareMathOperator{\Adj}{Adj}
\DeclareMathOperator{\LS}{LS}
\DeclareMathOperator{\PT}{PT}
\DeclareMathOperator{\TS}{TS}
\DeclareMathOperator{\PB}{PB}
\DeclareMathOperator{\DFBFS}{DFBFS}
\DeclareMathOperator{\BTRPT}{BTRPT}

\DeclareMathOperator{\vertices}{vert}

\DeclareMathOperator{\td}{td}
\DeclareMathOperator{\cone}{cone}
\DeclareMathOperator{\interior}{relint}

\DeclareMathOperator{\conv}{conv}   %%% Added by David. Denotes convex hull
\DeclareMathOperator{\vol}{vol}     %%% Added by David. Denotes volume
\DeclareMathOperator{\supp}{supp}        %%% Added by David. Denotes support.
\DeclareMathOperator{\affine}{aff}       %%% Added by David. Denotes affine space 
\DeclareMathOperator{\rank}{rank}
\DeclareMathOperator{\cl}{cl}       %%% Added by David. Denotes cl operator from Topkis
\DeclareMathOperator{\lin}{lin}     %%% Added by David. Denotes the linear space containing said object.

\let\epsilon=\varepsilon
\makeatletter
\newcommand{\DeclareBracket}[3]{
  \newcommand{#1}[2][]{%
  \ifthenelse%
  {\equal{##1}{}}%
  {\left#2##2\right#3}%
  {\csname ##1l\endcsname#2##2\csname ##1r\endcsname#3}}}    
\makeatother
\DeclareBracket\abs||           % absolute value
\DeclareBracket\norm\|\|        % norm
\DeclareBracket\floor\lfloor\rfloor
\DeclareBracket\ceil\lceil\rceil
\DeclareBracket\set\{\}
\DeclareBracket\paren()
\DeclareBracket\inner\langle\rangle
\DeclareBracket\fractional\{\}
\newcommand\C{\mathbb C}
\newcommand\inputfig[1]{\ifpdf
    \input{#1.pdf_t}
    \else
    \input{#1.pstex_t}
    \fi}

\newtheorem{theorem}{Theorem}%
\makeatletter
\newtheorem{lemma}{Lemma}
\renewcommand*{\c@lemma}{\c@theorem}
\renewcommand*{\p@lemma}{\p@theorem}

\newtheorem{conjecture}{Conjecture}
\renewcommand*{\c@conjecture}{\c@theorem}
\renewcommand*{\p@conjecture}{\p@theorem}

\newtheorem{proposition}{Proposition}
\renewcommand*{\c@proposition}{\c@theorem}
\renewcommand*{\p@proposition}{\p@theorem}

\newtheorem{corollary}{Corollary}
\renewcommand*{\c@corollary}{\c@theorem}
\renewcommand*{\p@corollary}{\p@theorem}

\renewcommand*{\c@observation}{\c@theorem}
\renewcommand*{\p@observation}{\p@theorem}

\theoremstyle{definition}

\renewcommand*{\c@problem}{\c@theorem}
\renewcommand*{\p@problem}{\p@theorem}

\newtheorem{definition}{Definition}
\renewcommand*{\c@definition}{\c@theorem}
\renewcommand*{\p@definition}{\p@theorem}

\newtheorem{remark}{Remark}
\renewcommand*{\c@remark}{\c@theorem}
\renewcommand*{\p@remark}{\p@theorem}

\newtheorem{example}{Example}
\renewcommand*{\c@example}{\c@theorem}
\renewcommand*{\p@example}{\p@theorem}

\newtheorem{algorithm}{Algorithm}
\renewcommand*{\c@algorithm}{\c@theorem}
\renewcommand*{\p@algorithm}{\p@theorem}

\newtheorem{heuristic}{Heuristic}
\renewcommand*{\c@heuristic}{\c@theorem}
\renewcommand*{\p@heuristic}{\p@theorem}

\makeatother

\pdfoutput=1
\usepackage{graphicx}

\begin{document}
   \frontmatter

   \pagestyle{prelim}
   
   % Redefine plain page style so that the first pages of chapters
   % have desired page style.
   %
   \fancypagestyle{plain}{%
      \fancyhf{}
      \cfoot{-\thepage-}
   }%
   \pdfoutput=1
\begin{center}
   \null\vfill
   \textbf{%
      Matroid Polytopes: Algorithms, Theory, and Applications
   }%
   \\
   \bigskip
   By \\
   \bigskip
   DAVID CARLISLE HAWS \\
   \bigskip
   B.S. (University of California, Davis) 2004 \\
   \bigskip
   DISSERTATION \\
   \bigskip
   Submitted in partial satisfaction of the requirements for the
   degree of \\
   \bigskip
   DOCTOR OF PHILOSOPHY \\
   \bigskip
   in \\
   \bigskip
   Mathematics \\
   \bigskip
   in the \\
   \bigskip
   OFFICE OF GRADUATE STUDIES \\
   \bigskip        
   of the \\
   \bigskip
   UNIVERSITY OF CALIFORNIA \\
   \bigskip
   DAVIS \\
   \bigskip
   Approved: \\
   \bigskip
   \bigskip
   \makebox[3in]{\hrulefill} \\
   Jes\'us A. De Loera (\emph{chair}) \\
   \bigskip
   \bigskip
   \makebox[3in]{\hrulefill} \\
   Matthias K\"oppe \\
   \bigskip
   \bigskip
   \makebox[3in]{\hrulefill} \\
   Eric Babson \\
   \bigskip
   \bigskip
   Committee in Charge \\
   \bigskip
   2009 \\
   \vfill
\end{center}

   \newpage

    % %------------------------------------------------------------------------------
    % 
    % %------------------------------------------------------------------ NEW PAGE --
    %
    % % -- COPYRIGHT PAGE
        
    % The following commands create a copyright notice page. This 
    % page can be deleted if you would prefer to not include it.

    ~\\[7.75in] % to place the copyright near the bottom of the page
    \centerline{
                \copyright\ David Carlisle Haws,
                            2009. All rights reserved.
               }
    \thispagestyle{empty}
    \addtocounter{page}{-1}

   \newpage
   
   % Begin Double Spacing
   %
   \doublespacing
   
   \tableofcontents
   \newpage
   
   \pdfoutput=1
% This ONLY goes in the single abstract page submission. For the abstract page that is included in the dissertation, there should be no name, major nor department.
%{\singlespacing
%   \begin{flushright}
%      David Carlisle Haws \\
%      June 2009 \\
%      Mathematics \\
%   \end{flushright}
%}

\bigskip

\begin{center}
Matroid Polytopes: Algorithms, Theory and Applications.
\end{center}

\section*{Abstract}
This dissertation presents new results on three different themes all related to matroid polytopes. First we investigate properties of Ehrhart polynomials of matroid polytopes, independence matroid polytopes, and polymatroids. We prove that for fixed rank their Ehrhart polynomials are computable in polynomial time. The proof relies on the geometry of these polytopes as well as a new refined analysis of the evaluation of Todd polynomials. 

Second, we discuss theoretical results regarding the algebraic combinatorics of matroid polytopes. We discuss two conjectures about the $h^*$-vector and coefficients of Ehrhart polynomials of matroid polytopes and provide theoretical and computational evidence for their validity. We also explore a variant of White's conjecture which states that every matroid polytope has a regular unimodular triangulation. We provide extensive computational evidence supporting this new conjecture and propose a combinatorial condition on simplices sufficient for unimodularity. Lastly we discuss properties of two dimensional faces of matroid polytopes.

Finally, motivated by recent work on algorithmic theory for non-linear and multicriteria matroid optimization, we have developed algorithms and heuristics aimed at practical solutions of large instances of these difficult problems.  Our methods primarily use the local adjacency structure inherent in matroid polytopes to pivot to feasible solutions which may or may not be optimal.  We also present a modified breadth-first-search heuristic that uses adjacency to enumerate a subset of feasible solutions. We present other heuristics, and provide computational evidence supporting these new techniques. We implemented all of our algorithms in the software package MOCHA (Matroids Optimization Combinatorial Heuristics and Algorithms). %\footnote{Matroids Optimization Combinatorial Heuristics and Algorithms}.

%    \thispagestyle{empty}
%    \addtocounter{page}{-1}

   \newpage
   
   \section*{Acknowledgments}
   \pdfoutput=1
First and foremost to my advisor Jes\'us De Loera who started me on this journey as an undergraduate. He guided me around the pitfalls, encouraged me, mentored me, provided me with problems and kept challenging me at every turn. His patience was remarkable and he was my number one cheerleader. Most of all, he taught me that mathematics without passion is not math at all.

To my wife Tami Joy whose support and love made this all possible. Her utmost confidence in me is amazing. It lifted my spirits and kept me going. She can always make me laugh, and is always there to lean on.

My thanks to Matthias K\"oppe for sharing his mathematical and computational skills. His careful and thoughtful approach to mathematics set a great example and he always made himself available for any math or computer questions. 

To my family and friends for all their support, understanding and love. My thanks to the UC Davis mathematics department staff who were very helpful and instrumental in my success.

Thanks to Francisco Santos for all our insightful conversations on triangulations of matroid polytopes and more.

My appreciation to Eric Babson for serving on my qualification exam and dissertation committee.

Support for this dissertation was provided by VIGRE (DMS-01-35345, DMS-0636297), University of California at Davis (DMS-0608785) and a gift from IBM.

   \mainmatter
   
   \pagestyle{maintext}
   
   % Redefine plain page style so that the first pages of 
   % chapters have desired page style.
   %
   \fancypagestyle{plain}{%
      \renewcommand{\headrulewidth}{0pt}
      \fancyhf{}
      \rhead{\thepage}
   }%
   
   \chapter[Matroids and Their Polytopes]{What are Matroids and who are their polytopes?}
   \label{ch:IntroductionLabel}
   \pdfoutput=1
\section{Matroids}
Matroids are combinatorial objects that naturally encapsulate the idea of independence. They are found in matrices, graphs, transversals, point configurations, hyperplane arrangements, greedy optimization, and pseudosphere arrangements to name a few.  One of the reasons matroids have become fundamental objects in pure and applied combinatorics are their many equivalent axiomatizations. For an excellent guide and thorough treatment of matroids we recommend \cite{Oxley1992Matroid-Theory,Welsh1976Matroid-Theory,Schrijver2003Combinatorial-O}, which we refer to for most all of our definitions and theory. We now offer a brief introduction to matroids and their associated polytopes. We begin with the most basic and prominent matroid definition (axiomatization). First we adopt the notation $X + y := X \cup \{y\}$ and $X - y := X \, \backslash \, \{y\}$. A pair $(S,\In)$ is called a \emph{matroid} if $S$ is a finite set and $\In$ is a nonempty collection of subsets of $S$ (denoted $2^S$) satisfying:
\begin{equation} \label{matroid:ind}
    \begin{split}
        (i)  \quad & \text{if } I \in \In \text{ and } J \subseteq I, \text{ then } J \in \In\\
        (ii) \quad & \text{if } I,J \in \In \text{ and } |I| < |J|, \text{ then } I + z \in \In \text{ for some } z \in J \backslash I.
    \end{split}
\end{equation}
For a matroid $\M = (S,\In)$, $I \subseteq S$ is \emph{independent} if $I \in \In$ and \emph{dependent} otherwise. $B \subseteq S$ is called a \emph{base} or \emph{basis} if it is an inclusion maximal independent set. That is, $B \in \In$ and there does not exist $A \in \In$ such that $B \subset A \subseteq S$. We denote the collection of bases of $\M$ as $\B_\M$, or sometimes $\B$ when the matroid is clear. An easy consequence of \autoref{matroid:ind} is that all bases have the same cardinality: If not, let $B_1, B_2$ be two bases of $(S,\In)$ such that $|B_1| < |B_2|$. Then by $(ii)$ in \autoref{matroid:ind} there exists $z \in B_2$ such that $B_1 + z \in \In$. But this contradicts the inclusionwise maximality of $B_1$. \\

Given any $U \subseteq S$ the \emph{rank} of $U$ is given by the \emph{rank function} of $\M$,  $\varphi\colon 2^{[n]} \rightarrow \Z$ where
\begin{equation*}
\varphi (A) := \max \{\, |X| \mid X \subseteq A, \, X \in \In \,\}.
\end{equation*}
% $2^S := \{\, A \subseteq S \,\}$.
The common cardinality of the bases of $\M$ is called the \emph{rank} of the matroid. Throughout we will refer to the rank of $\M$ as $\rank(\M) := \rank(S)$.
\begin{example}[Uniform matroid]
Let $[n] := \{1,\ldots,n\}$. The \emph{uniform} matroid of rank $r$ is $\U^{r,n} := \{\, A \subseteq [n] \mid |A| \leq r \, \}$. It is easy to see that $\U^{r,n}$ satisfies both $(i)$ and $(ii)$ of \autoref{matroid:ind}. 
\begin{equation*}
\U^{2,4} = \{\, \{1,2\}, \{1,3\}, \{1,4\}, \{2,3\}, \{2,4\}, \{3,4\} \, \}
\end{equation*}
\end{example}
We state an additional useful axiomatization on the independence sets.
\begin{lemma}[Theorem 39.1 \cite{Schrijver2003Combinatorial-O}] \label{lem:equivind}
Let $S$ be a finite set and let $\In$ be a nonempty collection of subsets satisfying \autoref{matroid:ind}$(i)$. Then \autoref{matroid:ind}$(ii)$ is equivalent to:
\begin{equation} \label{eq:altii}
\text{if } I,J \in \In \text{ and } |I \Backslash J| = 1,\, | J \Backslash I | = 2, \text{ then } I + z \in \In \text{ for some } z \in J \Backslash I.
\end{equation}
\end{lemma}
\begin{proof}
Obviously \autoref{matroid:ind}$(ii)$ implies \autoref{eq:altii}. Conversely, \autoref{matroid:ind}$(ii)$ follows from \autoref{eq:altii} by induction on $|I \Backslash J|$, the case $|I \Backslash J|=0$ being trivial. If $|I \Backslash J| \geq 1$, choose $i \in I \Backslash J$. We apply the induction hypothesis twice: first to $I-i$ and $J$ to find $j \in J \Backslash I$ with $I - i + j \in \In$, and then to $I - i + j$ and $J$ to find $j' \in J \Backslash (I + j)$ with $I - i + j + j' \in \In$. Then by \autoref{eq:altii} applied to $I$ and $I - i + j + j'$, we have $I + j \in \In$ or $I + j' \in \In$.
\end{proof}
We can now prove one of the most useful axiomatizations for the purposes of this dissertation. Considering the properties of \autoref{matroid:ind} it is natural to think that the bases of $\M$ are all that are needed to fully describe a matroid and indeed the following theorem gives an important characterization of matroids in terms of bases.
\begin{theorem}[Theorem 39.6 \cite{Schrijver2003Combinatorial-O}]
Let $S$ be a set and $\B$ be a nonempty collection of subsets of $S$. Then the following are equivalent:
\begin{equation} \label{eq:basisexch}
\begin{split}
(i)  \quad & \B \text{ is the collection of bases of a matroid;}\\
(ii) \quad & \text{if } B, B' \in \B \text{ and } x \in B' \Backslash B, \text{ then } B' - x + y \in \B \text{ for some } y \in B \Backslash B'; \\
(iii)\quad & \text{if } B, B' \in \B \text{ and } x \in B' \Backslash B, \text{ then } B' - y + x \in \B \text{ for some } y \in B \Backslash B'.
\end{split}
\end{equation}
\end{theorem}
\begin{proof}
$(i) \Rightarrow (ii)$: Let $\B$ be the collection of bases of a matroid $(S,\In)$. Then all sets in $\B$ have the same size. Now, let $B,B' \in \B$ and $x \in B' \Backslash B$. Since $B' - x \in \In$, there exists a $y \in B \Backslash B'$ with $B'' := B' - x + y \in \In$. Since $|B''| = |B'|$, we know $B'' \in \B$.\\

$(iii) \Rightarrow (i)$: $(iii)$ directly implies that no set in $\B$ is contained in another. Let $\In$ be the collection of sets $I$ with $I \subseteq B$ for some $B \in \B$. We check \autoref{eq:altii}. Let $I,J \in \In$ with $|I \Backslash J| = 1$ and $|J \Backslash I| = 2$. Let $I \Backslash J = \{x\}$.\\

Consider sets $B, B' \in \B$ with $I \subseteq B$, $J \subseteq B'$. If $x \in B'$, we are done. So assume $x \notin B'$. Then by $(iii)$, $B' - y + x \in \B$ for some $y \in B' \Backslash B$. As $|J \Backslash I| = 2$, there is a $z \in J \Backslash I$ with $z \neq y$. Then $I + z \subseteq B' - y + x$, and so $I + z \in \In$.\\

$(ii) \Rightarrow (iii)$: By the foregoing we know that $(iii)$ implies $(ii)$. Now axioms $(ii)$ and $(iii)$ interchange if we replace $\B$ by the collection of complements of sets in $\B$. Hence the implication $(ii)\Rightarrow (iii)$ also holds.
\end{proof}
Parts $(ii)$ and $(iii)$ of \autoref{eq:basisexch} might seem identical, but $(ii)$ says that we can pick an element to \emph{remove} from $B'$ and there exist some element of $B$ to add to get a basis while $(iii)$ says that we can pick an element to \emph{add} to $B'$ and there exists some element of $B$ to remove to get a basis. Having the flexibility to choose the element to add or remove is a useful property for many proofs, algorithms and heuristics on the bases of matroids.

\begin{example}[Graphical Matroid]
Let $G= (V,E)$ be graph with nodes $V$ and edges $E$. The collection of spanning forests of $G$ are the bases of the \emph{graphical} matroid on ground set $V$, denoted $\M_G$. We can verify the basis exchange axiom \autoref{eq:basisexch}$(iii)$ by considering two spanning forests $F_1, F_2$. If $e \in F_2 \Backslash F_1$, then $F_1 + e$ contains a cycle $C$ in $G$. Also there exists $f \in C \Backslash F_2$ and thus $F_1 - f + e$ is a spanning forest.
\end{example}
\begin{figure}[!htb]
\hfil(a)\subfigure{\scalebox{0.5}{\ifpdf
    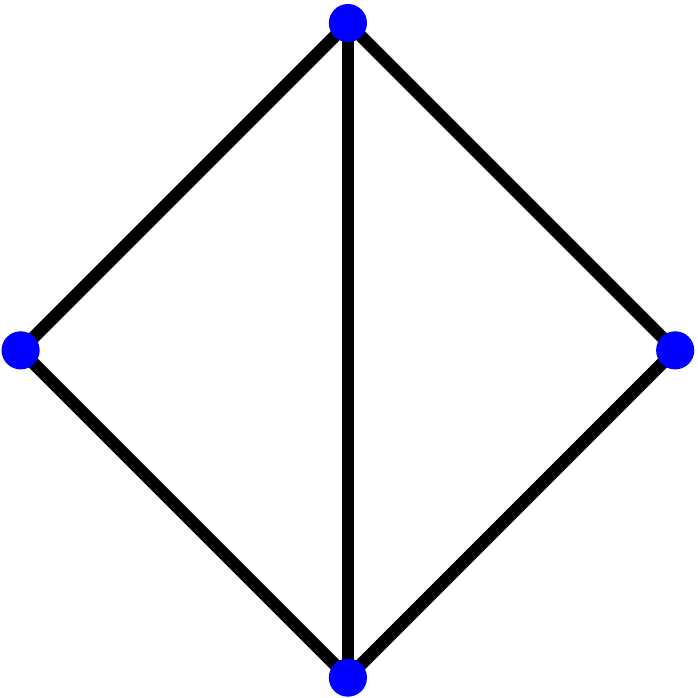
    \else
    \input{favgraph.pstex_t}
    \fi}}
\hfil(a)\subfigure{\scalebox{0.5}{\ifpdf
    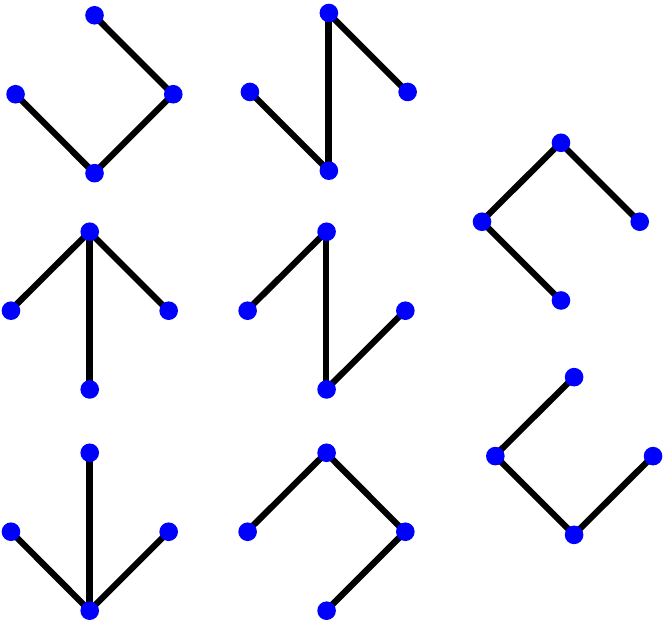
    \else
    \input{favgraphspanningtrees.pstex_t}
    \fi}}
\caption{(a) A graph $G$ and (b) its spanning forests, i.e. bases $\B_{\M_G}$ of $\M_G$.}
\end{figure}

Matroids play a prominent role in combinatorial optimization and are intrinsically tied to the greedy algorithm which we state now: If $w: S \rightarrow \R$ is a weight function we use the notation $w(I) := \sum_{i\in I} w(i)$ where $I \subseteq S$. Letting $\In \subset 2^S$ be closed under taking subsets and given a weighting $w$, the \emph{greedy algorithm} starts by setting $I := \emptyset$, and repeatedly choosing $y\in S \Backslash I$ such that $I + y \in \In$ where $w(y)$ is as large as possible. We terminate if no such $y$ exists.

For arbitrary $\In \subseteq 2^S$ the greedy algorithm does not always terminate at the maximum solution, although for matroids this will always be the case. In fact, the converse is true: If the greedy algorithm always works on $\In$ then it is a matroid.

\begin{theorem}[Theorem 40.1 \cite{Schrijver2003Combinatorial-O}] \label{thm:greedy}
Let $\In$ be a nonempty collection of subsets of a set $S$, closed under taking subsets. Then the pair $(S,\In)$ is a matroid if and only if for each weight function $w:S \rightarrow \R_+$, the greedy algorithm leads to a set $I \in \In$ of maximum weight $w(I)$.
\end{theorem}
\begin{proof}
\emph{Necessity.} Let $(S,\In)$ be a matroid and let $w: S \rightarrow \R_+$ be any weight function on $S$. Call an independent set $I$ good if it is contained in a maximum-weight basis. It suffices to show that if $I$ is good, and $y$ is an element in $S \Backslash I$ with $I + y \in \In$ and with $w(y)$ as large as possible, then $I + y$ is good.

As $I$ is good, there exists a maximum-weight basis $B \supseteq I$. If $y \in B$, then $I + y$ is good again. If $y \notin B$, then there exists a base $B'$ containing $I+y$ and contained in $B+y$. So $B' = B - z + y$ for some $z \in B \Backslash I$. As $w(y)$ is chosen maximum and as $I+z \in \In$ since $I + z \subseteq B$, we know $w(y) \geq w(z)$. Hence $w(B') \geq w(B)$, and therefore $B'$ is a maximum-weight base. So $I+y$ is good.

\emph{Sufficiency.} Suppose that the greedy algorithm leads to an independent set of maximum weight for each weight function $w: S \rightarrow \R_+$. We show that $(S,\In)$ is a matroid.

Condition \autoref{matroid:ind}(i) is satisfied by assumption. To see condition \autoref{matroid:ind}(ii), let $I,J \in \In$ with $|I| < |J|$. Suppose that $I+z \notin \In$ for each $z \in J \Backslash I$. Let $k := |I|$. Consider the following weight function $w$ on $S$:

\begin{equation*}
w(S) := \left\{ \begin{array}{cl} 
k+2 & \text{if } s \in I, \\
k+1 & \text{if } s \in J \Backslash I,\\
0   & \text{if } s \in S \Backslash (I \cup J).
\end{array} \right.
\end{equation*}

Now in the first iterations of the greedy algorithm we find the $k$ elements in $I$. By assumption, at any further iteration, we cannot chose any element in $J \Backslash I$. Hence any further element chosen, has weight $0$. So the greedy algorithm yields an independent set of weight $k(k+2)$.

However, $J$ has weight at least $|J|(k+1) \geq (k+1)(k+1) > k(k+2)$. Hence the greedy algorithm does not give a maximum-weight independent set, contradicting our assumption.
\end{proof}

The rank function $\varphi$ of a matroid $\M$ also determines $\M$ and we can easily see that a set $A \subseteq S$ is independent if and only if $\varphi(A) = |A|$.
\begin{theorem}[Theorem 39.8 \cite{Schrijver2003Combinatorial-O}]
Let $S$ be a set and $\varphi: 2^S \rightarrow \Z_{+}$. Then $\varphi$ is the rank function of a matroid if and only if for all $T,U \subseteq S$:
\begin{equation} \label{matroid:rank}
\begin{split}
(i) \quad  & \varphi(T) \leq \varphi(U) \leq |U| \text{ if } T \subseteq U,\\
(ii) \quad & \varphi(T \cap U) + \varphi(T \cup U) \leq \varphi(T) + \varphi(U).
\end{split}
\end{equation}
\end{theorem}
\begin{proof}
\emph{Necessity}. Let $\varphi$ the rank function of a matroid $(S,\In)$. Choose $T,U \subseteq S$. Clearly \autoref{matroid:rank}(i) holds. To see $(ii)$, let $I$ be an inclusionwise maximal set in $\In$ with $I \subseteq T \cap U$ and let $J$ be an inclusionwise maximal set in $\In$ with $I \subseteq J \subseteq T \cup U$. Since $(S,\In)$ is a matroid, we know that $\varphi(T\cap U) = |I|$ and $\varphi(T \cup U) = |J|$. Then
\begin{equation*}
\varphi(T) + \varphi(U) \geq |J \cap T| + |J \cap U| = |J \cap (T \cap U)| + |J\cap(T \cup U)| \geq |I| + |J| = \varphi(T \cap U) + \varphi(T \cup U);
\end{equation*}
that is, we have \autoref{matroid:rank}(ii).

\emph{Sufficiency} Let $\In$ be the collection of subsets $I$ of $S$ with $\varphi(I) = |I|$. We show that $(S,\In)$ is a matroid, with rank function $\varphi$. 

Trivially, $\emptyset \in \In$. Moreover, if $I \in \In$ and $J \subset I$, then 
\begin{equation*}
\varphi(J) \leq \varphi(I) - \varphi(I \Backslash J) \geq |I| - |I \Backslash J|= |J|.
\end{equation*}
So $J \in \In$.

In order to check \autoref{eq:altii}, let $I, J \in \In$ with $|I \Backslash J| = 1$ and $|J \Backslash I| = 2$. Let $J \Backslash I = \{z_1,z_2\}$. If $I + z_1, I+z_2 \notin \In$, we have $\varphi(I+z_1)=\varphi(I+z_2) = |I|$. Then by \autoref{matroid:rank}(ii),
\begin{equation*}
\varphi(J) \leq \varphi(I+z_1 + z_2) \leq \varphi(I+z_1) + \varphi(I+z_2) - \varphi(I) = |I| < |J|,
\end{equation*}
contradicting the fact that $J \in \In$.

So $(S,\In)$ is a matroid. Its rank function is $\varphi$, since $\varphi(U) = \max\{\, |I| \mid I \subseteq U,\, I \in \In \, \}$ for each $U \in S$. Here $\geq$ follows from \autoref{matroid:rank}(ii), since if $I \subseteq U$ and $I \in \In$, then $\varphi(U) \geq \varphi(I) = |I|$. Equality can be shown by induction on $|U|$, the case $U = \emptyset$ being trivial. If $U \neq \emptyset$, choose $y \in U$. By induction, there is an $I \subseteq U - y$ with $I \in \In$ and $|I| = \varphi(U-y)$. If $\varphi(U) = \varphi(U-y)$ we are done, so assume $\varphi(U) > \varphi(U-y)$. Then $I + y \in \In$, since $\varphi(I+y) \geq \varphi(I) + \varphi(U) - \varphi(U-y) \geq |I| + 1$. Moreover, $\varphi(U) \leq \varphi(U-y) + \varphi(\{y\}) \leq |I| + 1$. This proves equality for $U$.
\end{proof}
%Functions satisfying \autoref{matroid:rank}(ii) are called \emph{submodular} and we will see later that functions of this type, with other conditions, will define polymatroids. 
One of the most natural methods of presenting a matroid in a computational model is by \emph{independence oracle} which given any $A \subseteq S$ asserts whether $A \in \In$ or not. The independence oracle allows us to compute the rank of any $A \subseteq S$. Begin with $I := \emptyset$ and for each $i \in A$ (any order) if $I + i \in \In$ then set $I := I + i$. Repeat and thus, $\rank(A) := |I|$. For an interesting overview of the computational relation of various matroid axiomatization representations see \cite{Mayhew}.

In $\autoref{thm:greedy}$ it was assumed that the weights are non-negative. But even with a weight function $w:S \rightarrow \R$ it can still be shown that for a matroid $\M$ the greedy algorithm finds a maximum-weight basis. Not surprising, one can replace the condition ``as large as possible'' with ``as small as possible'' and the greedy algorithm will find a minimum-weight basis. Assuming an independence oracle and \autoref{thm:greedy} the corollary follows: 
\begin{corollary}[Corollary 40.1a \cite{Schrijver2003Combinatorial-O}]
A maximum-weight independent set in a matroid can be found in strongly polynomial time.
\end{corollary}

\begin{example}[Linear, Realizable, or Vector Matroids]
Let $A := \{a_1 \, | \, a_2 \, | \, \cdots \, | \, a_n\} \in \R^{m \times n}$ be an $n \times m$ matrix. Let $S := \{1,\ldots,n\}$ and we say $I \subseteq S$ is independent if the columns of $A$ indexed by $I$ are linearly independent. %In other words if $I := \{i_1,\ldots,i_n\}$, then $I$ is independent if $det(a_{i_1},\ldots,a_{i_n}) \neq 0$.
\begin{equation*}
A = \bordermatrix{ & \phantom{-}\mathbf{1} & \phantom{-}\mathbf{2} & \phantom{-}\mathbf{3} & \phantom{-}\mathbf{4} & \phantom{-}\mathbf{5} \cr
& \phantom{-}1  & \phantom{-}0 & \phantom{-}1 & -1 & \phantom{-}2 \cr
& \phantom{-}1  & \phantom{-}1 & \phantom{-}0 & \phantom{-}1 & \phantom{-}2} 
\qquad
\{\mathbf{1},\mathbf{2}\}, \{\mathbf{2},\mathbf{3}\} \text{ are indepdendent while } \{\mathbf{1},\mathbf{5}\} \text{ is not.}
\end{equation*}
Such a matriod $\M_A$ is called \emph{linear} or \emph{realizable}. If the entries of $A$ are in a field $\mathbb F$ then we say $\M_A$ is realizable over $\mathbb F$. Graphical matroids are linear and realizable over any field. This can be done by arbitrarily orienting the edges of $G$ and take its incidence matrix. Not all matroids are linear, e.g. see $F_8$ in the appendix of \cite{Oxley1992Matroid-Theory}. Some matroids can only be realized in certain fields. For example, the Fano matroid can only be realized in a field of characteristic two \cite{Oxley1992Matroid-Theory}. %Smallest. See Welsch pg 509 F_8
\end{example}

\section{Matroid Polytopes \& Polymatroids}
\begin{figure} 
\scalebox{0.3}{\ifpdf
    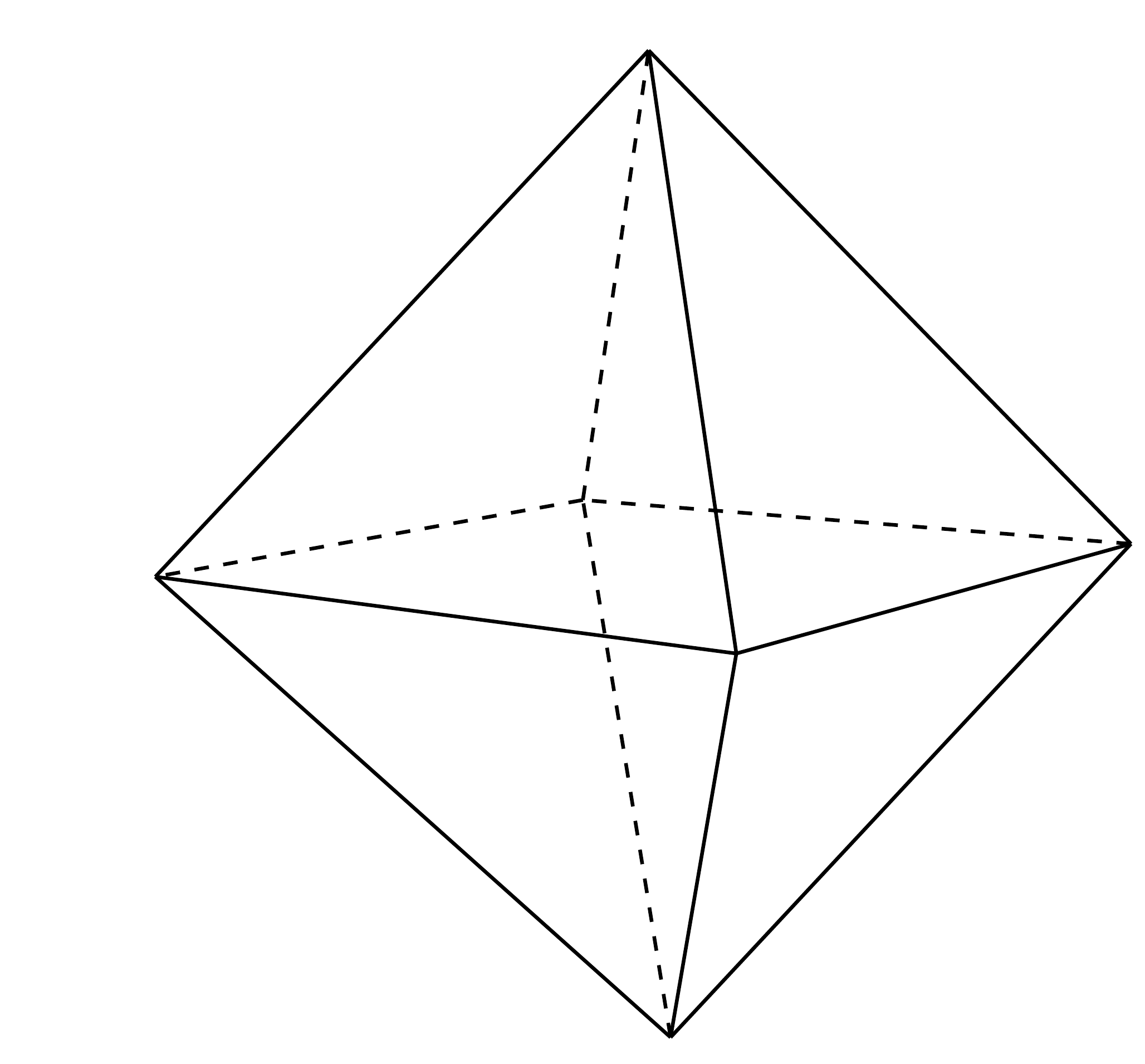
    \else
    \input{crosspoly.pstex_t}
    \fi}
\caption{The cross polytope $\Po_{U^{2,4}}$.} \label{fig:crosspoly}
\end{figure}

Let $A \subset [n] := \{1,\ldots,n\}$ and define the \emph{incidence vector} of $A$ as 
\begin{equation*}
\ve e_A := \sum_{i \in A} \ve e_i
\end{equation*}
where $\ve e_i \in \R^n$ are the standard unit vectors. For example if $A = \{1,3,4\} \subseteq [5]$ then $\ve e_A = (1,0,1,1,0)^\top \in \R^5$. We also adopt the notation $\Inc(X) := \{\, \ve e_A \mid A \in X  \,\}$ where $X \subseteq 2^S$.
Now we introduce the main object of this thesis.
\begin{definition}
Let $\M$ be a matroid with bases $\B_\M$. The \emph{matroid polytope} of $\M$ is
\begin{equation*}
    \Po_\M := \conv( \ve e_B \mid B \in \B_\M)
\end{equation*}
where $\conv(\cdot)$ denotes the convex hull.
\end{definition}
\begin{example}
Recall $U^{2,4} = \{\, \{1,2\}, \{1,3\}, \{1,4\}, \{2,3\}, \{2,4\}, \{3,4\} \,\}$. Thus 
\begin{equation*}
\Po_{U^{2,4}} = \conv(\,(1,1,0,0)^\top, (1,0,1,0)^\top, (1,0,0,1)^\top, (0,1,1,0)^\top, (0,1,0,1)^\top, (0,0,1,1)^\top \, ).
\end{equation*}
The polytope $\Po_{U^{2,4}} \in \R^4$ but we can easily see that all vertices lay in a hyperplane since they all sum to two. $\Po_{U^{2,4}}$ is in fact isomorphic to the cross polytope shown in \autoref{fig:crosspoly}.
\end{example}
The matroid polytope is different than the following well-known polytope.
\begin{definition}
Let $\M = (S,\In)$ be a matroid. The \emph{independence matroid polytope} of $\M$ is
\begin{equation*}
    \Po_\M^\In := \conv( \ve e_I \mid I \in \In).
\end{equation*}
\end{definition}
We can see that $\Po(M) \subseteq \Po^\In(M)$ and
$\Po_\M$ is a face of $\Po^\In_\M$ in the hyperplane $\sum_{i=1}^n x_i = \rank(M)$.

Polymatroids are closely related to matroid polytopes and independence matroid polytopes. A function $\psi \colon 2^{[n]} \longrightarrow \R$ is \emph{submodular} if $\psi(X \cap Y) + \psi(X \cup Y) \leq \psi(X) + \psi(Y)$ for all $X, Y \subseteq [n]$ and \emph{non-decreasing} if $\psi(X) \leq \psi(Y)$ for all $X \subseteq Y \subseteq [n]$. We say $\psi$ is a \emph{polymatroid rank function} if it is submodular, non-decreasing, and $\psi(\emptyset) = 0$. For example, the rank function of a matroid is a polymatroid rank function. The \emph{polymatroid} determined by a polymatroid rank function $\psi$ is the convex polyhedron in $\R^n$ given by
\begin{equation*}
    \Po_\psi := \Big\{ \, \ve x \in \R^n \ \Big| \ \sum_{i \in A} x_i \leq \psi(A) ,\, \forall A \subseteq [n], \, \ve x \geq \ve0 \, \Big\}  .
\end{equation*}

Independence matroid polytopes are a special class of polymatroids. Indeed, if $\varphi$ is a rank function on some matroid $\M$, then $\Po^\In_\M = \Po_\varphi $ \cite{Edmonds2003Submodular-func}. Moreover, the matroid polytope $\Po_\M$ is the face of $\Po_\varphi $ lying in the hyperplane $ \sum_{i=1}^n x_i = \varphi([n]) $.  Matroid polytopes and polymatroids appear in combinatorial optimization \cite{Schrijver2003Combinatorial-O}, algebraic combinatorics \cite{Feichtner2004Matroid-polytop}, and algebraic geometry \cite{Gelfand1987Combinatorial-g}. 
\begin{figure}[!htb]
\scalebox{0.3}{\ifpdf
    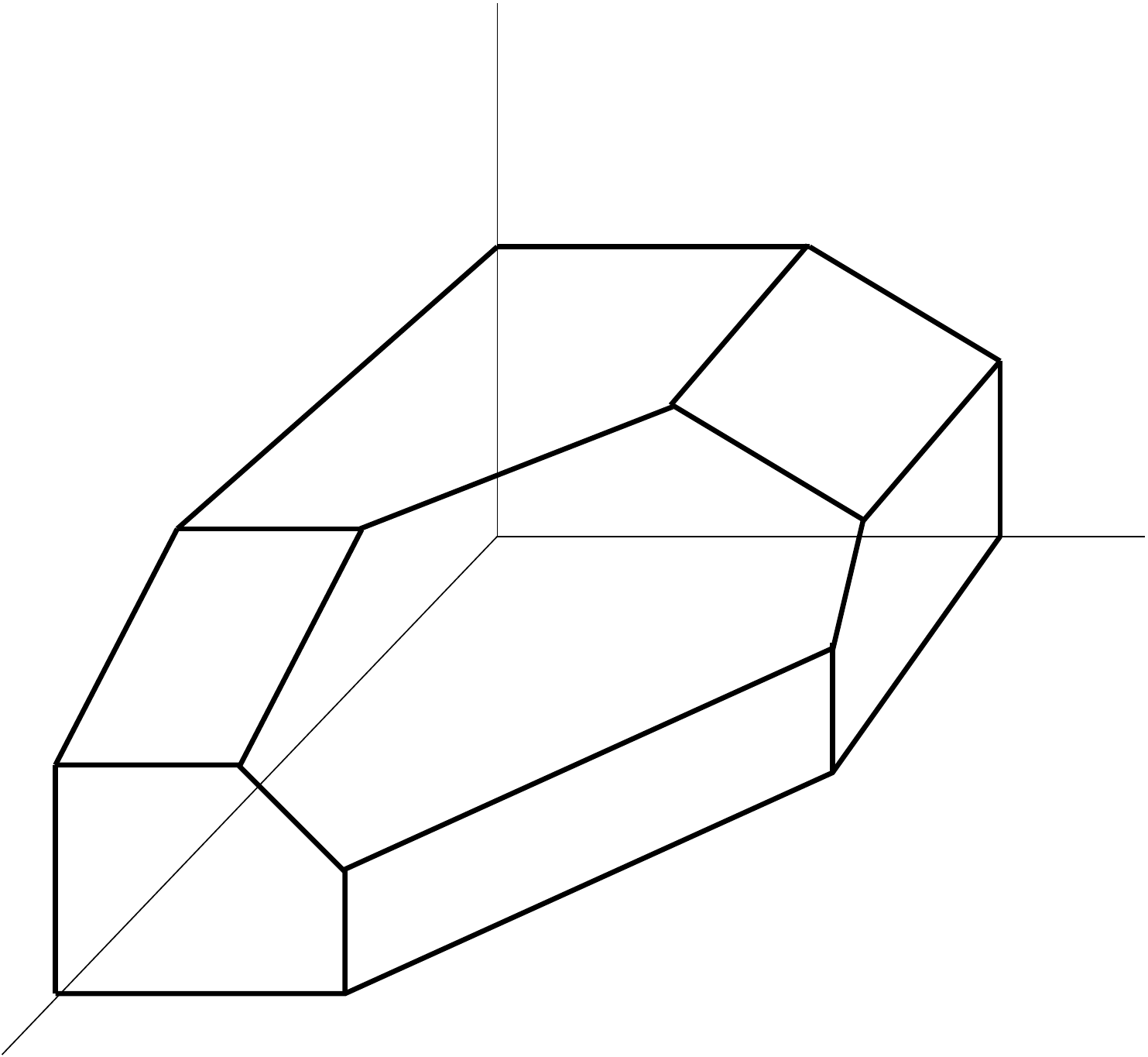
    \else
    \input{polymatroidex1.pstex_t}
    \fi}
\caption{An example of a polymatroid} \label{fig:polymatroid}
\end{figure}
Edmonds showed that a linear function can be optimized over a polymatroid $\Po_\psi$ via an extension of the greedy algorithm \cite{Edmonds2003Submodular-func}. This also provides a technique which easily lists off the vertices of a polymatroid $\Po_\psi$. One important implication we will see is that if $\psi$ is an integral polymatroid rank function then $\Po_\psi$ will be integral, that is, all the vertices of $\Po_\psi$ are integral. 

First we recall some basic definitions and key theorems from linear programing \cite{schrijver}. If $S$ is a finite subset then $\R^S$ is the vector space of vectors with coordinates indexed by $S$. That is, if $\ve x \in \R^S$ and $s \in S$ then the $s$th coordinate of $\ve x$ is given by $ x(s) \in \R$. If $A \subseteq S$ then $\ve x(A) := \sum_{a \in A} x(a)$.

Let $S$, $T$ be finite sets, $\ve c \in \R^S$, $\ve b \in \R^T$ and $A = \{\, a_{ij} \in \R \mid i \in T,\, j \in S \, \}$ be a matrix. A \emph{linear program} is 
\begin{equation*}
    \max \{\, \ve c \cdot \ve x \mid A\ve x \leq \ve b, \, \ve x \geq 0,\, \ve x \in \R^S \, \}.
\end{equation*}
The \emph{dual linear program} is
\begin{equation*}
    \min \{\, \ve b \cdot \ve y \mid A^\top \ve y \leq \ve c, \, \ve y \geq 0,\, \ve y \in \R^T \, \}.
\end{equation*}
A vector $\ve x \in \R_+^S$ satisfying $A\ve x \leq \ve b$ is a \emph{feasible solution} and a vector $\ve y \in \R_+^T$ satisfying $ A^\top \ve y \leq \ve c$ is a \emph{feasible dual solution}.

\begin{theorem}[\cite{Schrijver2003Combinatorial-O}]
For any linear program maximization problem exactly one of the following is true;
\begin{itemize}
    \item[(i)]   There exists no feasible solution,
    \item[(ii)]  For any $\alpha \in \R$ there is a feasible solution $\ve x$ such that $ \ve c \cdot \ve x > \alpha$,
    \item[(iii)] There is an optimal (feasible) solution.
\end{itemize}
\end{theorem}
A linear program and its dual are related and we state the famous \emph{strong duality theorem}.
\begin{theorem}[Corollary 7.1g \cite{schrijver}]
Let $S,T$ be finite sets, $A \in \R^S \times \R^T$, $\ve c \in \R^S$, $\ve b \in \R^T$. Then
\begin{equation*}
    \max \{\, \ve c \cdot \ve x \mid A\ve x \leq \ve b, \, \ve x \geq 0,\, \ve x \in \R^S \, \} = \min \{\, \ve b \cdot \ve y \mid A^\top \ve y \leq \ve c, \, \ve y \geq 0,\, \ve y \in \R^T \, \}
\end{equation*}
if both sets are nonempty.
\end{theorem}

In order to justify the following approach to finding the vertices of a polymatroid we will use a fact from linear programming and polyhedral theory; every vertex of a polyhedra is the optimal solution to some linear program. Let $\varphi$ be a polymatroid rank function on the finite set $S$ and $\ve w : S \rightarrow \R$. We start by ordering the elements of $S$ as $s_1,\ldots,s_n$ such that $w(s_1) \geq w(s_2) \geq \cdots \geq w(s_n)$ and define $U_i := \{s_1,\ldots,s_i\}$ for $i=0,\ldots,n$. Next, define $\widehat{\ve x} \in \R^S$ as 
\begin{equation} \label{eq:polyx}
\widehat x(s_i) := \varphi(U_i) - \varphi(U_{i-1}) \text{  for } i=1,\ldots,n. 
\end{equation}
The following theorem will show that $\widehat{\ve x}$ maximizes $\ve w \cdot \ve x$ over $\Po_\varphi$. Consider the linear program and its dual:
\begin{equation} \label{eq:polylindual}
\begin{split}
& \max \{\, \ve w \cdot \ve x \mid \ve x \in \Po_\varphi \,\} \\
        & = \min \left\{ \, \sum_{T \subseteq S} \ve y(T) \varphi(T) \mid \ve y \in \R^{2^S}_+, \sum_{T \subseteq S} \ve y(T) \cdot e_T \geq \ve w\,\right\}
\end{split}
\end{equation}
Also, define:
\begin{equation} \label{eq:polyyopt}
\begin{array}{rlr}
    \widehat{\ve y} (U_i) = & w(s_i) - w(s_{i+1}) &\quad  (i=1,\ldots,n-1)\\
    \widehat{\ve y} (S)   = & w(s_n),\\
    \widehat{\ve y} (T)   = & 0                  &\quad       (T \neq U, \text{ for each } i).
\end{array}
\end{equation}

\begin{theorem}[Theorem 44.3 \cite{Schrijver2003Combinatorial-O}] \label{thm:polyopt}
Let $\varphi$ be a polymatroid rank function on $S$ and let $\ve w: S \rightarrow \R^S$. Then $\widehat{\ve x}$ and $\widehat{\ve y}$ are optimum solutions to \autoref{eq:polylindual}.
\end{theorem}
\begin{proof}
We first show that $\widehat{\ve x} \in \Po_\varphi$; that is, $\widehat{\ve x} \leq \varphi(T)$ for each $T \subseteq S$. This is shown by induction on $|T|$, the case $T = \emptyset$ being trivial. Let $T \neq \emptyset$ and let $k$ be the largest index with $s_k \in T$. Then by induction,
\begin{equation*}
    \widehat{\ve x}(T - s_k) \leq \varphi(T - s_k).
\end{equation*}
Hence
\begin{equation*}
    \widehat{\ve x}(T) \leq \varphi(T - s_k) + \widehat{\ve x}(x_k) = \varphi(T - s_k) + \varphi(U_k) - \varphi(U_{k-1}) \leq \varphi(T)
\end{equation*}
(the last inequality follows from the submodularity of $\varphi$). So $\widehat{\ve x} \in \Po_\varphi$.

Also, $\widehat{\ve y}$ is feasible for \autoref{eq:polylindual}. Trivially, $\widehat{\ve y} \geq \ve 0$. Moreover, for any $i$ we have by \autoref{eq:polyyopt}:
\begin{equation*}
 \sum_{T \ni s_i} \widehat{\ve y}(T) = \sum_{j \geq i} \widehat{\ve y}(U_j) = \ve w(s_i).
\end{equation*}
So $\widehat{\ve y}$ is a feasible solution of $\autoref{eq:polylindual}$.

Optimality of $\widehat{\ve x}$ and $\widehat{\ve y}$ follows from:
\begin{equation*}
\begin{split}
& \ve w \cdot \widehat{\ve x} = \sum_{s \in S} w(s) \widehat x(s) = \sum_{i=1}^n w(s_i)(\varphi(U_i) - \varphi(U_{i-1})) \\
& = \sum_{i=1}^{n-1} \varphi(U_i)(w(s_i) - w(s_{i+1})) + \varphi(S)w(s_n) = \sum_{T \subseteq} \widehat{\ve y}(T)\varphi(T).
\end{split}
\end{equation*}
The third equality follows from a straightforward reordering of the terms, using that $\varphi(\emptyset) = 0$.
\end{proof}

Thus, \autoref{eq:polyx} gives a vertex of $\Po_\varphi$. If we vary the weights we will get all the vertices of $\Po_\varphi$ and it suffices to choose weights that permute the elements of $S$. \autoref{eq:polyx} and \autoref{thm:polyopt} imply the following important corollary.

\begin{corollary}[Corollary 44.3d \cite{Schrijver2003Combinatorial-O}]
Let $\varphi$ be an integral polymatroid rank function. Then the vertices of $\Po_\varphi$ are integral.
\end{corollary}

Given a matroid $\M$ of rank $r$ it is easy to see that $\Po_\M$ is contained in the $(n-1)$-simplex
\begin{equation*}
\Delta = \{\, (x_1,\ldots,x_n) \in \R^n \mid x_i \geq 0,\, x_1 + \cdots + x_n = r. \, \}
\end{equation*}
Moreover, the following theorem of Gel'fand, Goresky, MacPherson, and Serganova gives an equivalent geometric axiomatizatio of matroids:
\begin{theorem}[\cite{Gelfand1987Combinatorial-g}]
Let $\B$ be a family of $r$-element subsets of $[n]$. The following are equivalent:
    \begin{itemize}
        \item[(i)] $\B$ are the bases for a matroid,
        \item[(ii)] Every edge of $\conv (\, \ve e_B \mid B \in \B \, )$ is parallel to an edge of the simplex $\Delta $.
    \end{itemize}
\end{theorem}
The following is a more useful characterization of adjacency of vertices of $\Po_\M$, as well as a characterization of adjacency in the independence matroid polytope $\Po_\M^\In$.

\begin{lemma}[See Theorem 4.1 in \cite{Gelfand1987Combinatorial-g}, Theorem 5.1 and Corollary 5.5 in \cite{Topkis1984Adjacency-on-Po}] \label{lem:adj}
    Let $\M$ be a matroid. 
    \begin{itemize}
    \item[(i)]   Two vertices $\ve e_{B_1}$ and $\ve e_{B_2}$ are adjacent in $\Po_\M$ if and only if $\ve e_{B_1} - \ve e_{B_2} = \ve e_i - \ve e_j$ for some $i,j$.
    \item[(ii)]   If two vertices $\ve e_{I_1}$ and $\ve e_{I_2}$ are adjacent
      in $\Po^\In_\M$ then $\ve e_{I_1} - \ve e_{I_2} \in \{ \, \ve e_i - \ve
      e_j, \, \ve e_i, \, -\ve e_j \, \}$ for some $i,j$. Moreover if $\ve v$
      is a vertex of $\Po^\In_\M$ then all adjacent vertices of $\ve v$ can
      be computed in polynomial time in $n$, even if the matroid $\M$ is only
      presented by an evaluation oracle of its rank function~$\varphi$.
    \end{itemize}
\end{lemma}

We prove one direction of $(i)$, that every edge of $\Po_\M$ is parallel to $\ve e_i - \ve e_j$ for some $i,j \in \N$. The following can be found in the proof of Proposition $2.2$ in \cite{onn-2003-17}.
\begin{proof}
Consider any pair $A,B \in \B_\M$ of bases such that $[\ve e_A, \ve e_B]$ is an edge (that is, a $1$-face) of $\Po_\M$, and let $\ve a \in \R^n$ be a linear functional maximized over $\Po_\M$ uniquely on that edge. If $A \Backslash B = \{i\}$ is a singleton, then $B \Backslash A = \{j\}$ is a singleton as well, in which case $\ve e_A - \ve e_B = \ve e_i - \ve e_j$ and we are done. Suppose then, indirectly, that it is not, and pick and element $i$ in the symmetric difference $(A \Backslash B) \cup (B \Backslash A)$ of $A$ and $B$ of minimum value $a_i$. Without loss of generality assume $i \in A \Backslash B$. Then there is a $j \in B \Backslash A$ such that $C := A - i + j$ is a basis of $\M$. Since $|(A \Backslash B) \cup (B \Backslash A)| > 2$, $C$ is neither $A$ nor $B$. By the choice of $i$, this basis satisfies $a \cdot \ve e_C = \ve a \cdot \ve e_A - a_i + a_j \geq \ve a \cdot \ve e_A$, and hence $\ve e_C$ is also a maximizer of $\ve a$ over $\Po_\M$ and so lies in the $1$-face $[\ve e_A, \ve e_B]$. But no $\{0,1\}$-vector is a convex combination of other, yielding a contradiction.
\end{proof}

A similar property holds for integral polymatroids which we also state but we first recall some needed
definitions from \cite{Topkis1984Adjacency-on-Po}.  Let $\ve v,\ve w \in
\R^n$ and define $\Delta(\ve v,\ve w) := \{ \, i \in [n] \mid v_i \neq w_i
\, \}$ and $\cl(\ve v) := \{ \, G \mid G \subseteq [n], \ \sum_{i \in G}
v_i = \psi(G) \, \}$. Let $F = \{ f_1,\ldots,f_{|F|} \}$ be an ordered
subset of $[n]$ and $F_i := \{f_1,\ldots,f_i\}$. If $\psi$ is a polymatroid
rank function then we construct $\ve v \in \R^n$ where $v_i = \psi(F_i) -
\psi(F_{i-1})$ where $ v_j = 0$ when $j \notin F$ and one says $F$
\emph{generates} $\ve v$. %A classical result of Edmonds \cite{Edmonds2003Submodular-func} says that the set of vectors generated by all ordered subsets of $[n]$ is exactly the set of vertices of $\Po(\psi)$. 
%Now we can restate an important lemma. 

\begin{lemma}[See Theorem 4.1 and Section 2 in \cite{Topkis1984Adjacency-on-Po}.] \label{lem:topadj}
    Let $\psi$ be a polymatroid rank function. If $\ve v$ and $\ve w$ are vertices
    of the polymatroid $\Po(\psi)$ then either
    \begin{itemize}
        \item[(i)] $|\Delta(\ve v, \ve w)| = 1$ or
        \item[(ii)] $\cl(\ve v) = \cl(\ve w)$ and $\Delta(\ve v,\ve w) =
        \{c,d\}$ for some $c,d \in [n]$ where there exists some ordered set $F
        = \{f_1,\ldots,f_{|F|}\}$ which generates $\ve v$ with $f_{k+1} = d$
        and $f_k = c$ for some integer $k$, $1 \leq k \leq |F| -1$; moreover
        the ordered set\\ $\tilde F :=
        \{f_1,\ldots,f_{k-1},f_{k+1},f_k,f_{k+2},\ldots,f_{|F|}\}$ generates
        $\ve w$.
    \end{itemize}
\end{lemma}

In order to describe the matroid polytope $\Po_\M$ by a system of linear inequalities we require a few more definitions. A \emph{circuit} $C \subseteq S$ of $\M = (S,\In)$ is an inclusionwise minimal dependent set. That is, a dependent set $C$ is a circuit of $\M$ if there does not exist a dependent set $A$ such that $ A \subset C \subseteq S$. A \emph{flat} $F$ of $\M$ is a subset of $S$ such that $\varphi(F + s) = \varphi(F) + 1$ for all $s \in S$. The \emph{span} of $A \subseteq S$ is the smallest flat $F$ such that $A \subseteq F$. A matroid is called \emph{connected} if for every $x,y \in S$ there exists a circuit $C$ such that $x,y \in C$. In fact the circuits of a matroid $\M$ define an equivalence class on the ground set $S$ where $x,y$ are equivalent if there exists a circuit $C$ such that $x,y \in C$. The \emph{connected components} of $\M$ are defined as the number of such equivalence classes.

\begin{lemma}[\cite{Feichtner2004Matroid-polytop}]
The matroid polytope $\Po_\M$ equals the following subset of the simplex $\Delta$:
\begin{equation*}
 \Po_\M = \left\{\, (x_1,\ldots,x_n) \in \Delta \mid \sum_{i \in F} x_i \leq \rank(F) \text{ for all flats } F \subseteq S \, \right\}
\end{equation*}
\end{lemma}
\begin{proof}
Consider any facet of the polytope $\Po_\M$ and let $\sum_{i=1}^n a_ix_i \leq b$ be an inequality defining this facet. The normal vector $(a_1,a_2,\ldots,a_n)$ is perpendicular to the edges of that facet. But each edge of tht facet is parallel to some difference of unit vectors $e_i - e_j$ by \autoref{lem:adj}. Hence the only constraints on the coordinates of the normal vector are of the form $a_i=a_j$. Using the equation $\sum_{i=1}^n x_i = r$ and scaling the right hand side $b$, we can therefore assume that $(a_1,a_2,\ldots,a_n)$ is a vector in $\{0,1\}^n$. Hence the polytope $\Po_\M$ is characterized by the inequalities of the form $\sum_{i \in G} x_i \leq b_G$ for some $G \subseteq [n]$. The right hand side $b_G$ equals
\begin{equation*}
    b_G = \max\{\, |\sigma| \cap G \mid \sigma \text{ basis of } \M \, \} = \rank(G).
\end{equation*}
The second equality holds because every independent subset of $G$ can be completed to a basis $\sigma$. Let $F$ be the flat spanned by $G$. Then $G \subseteq F$ and $\rank(G) = \rank(F)$, and hence the inequality $\sum_{i\in F} x_i \leq \rank(F)$ implies the inequality $\sum_{i\in G} x_i \leq \rank(G)$.
\end{proof}

\begin{lemma}[\cite{Feichtner2004Matroid-polytop}]
The dimension of the matroid polytope $\Po_\M$ equals $n - \#$connected components of $\M$.
\end{lemma}
\begin{proof}
Two elements $i$ and $j$ are equivalent if and only if there exist bases $\sigma$ and $\tau$ with $i \in \sigma$ and $\tau = \sigma + i - j$. The linear space parallel to the affine span of $\Po_\M$ is spanned by the vectors $\ve e_i - \ve e_j$ arising in this manner. The dimension of this space equals $n - c(M)$.
\end{proof}

\subsection{Matroid Optimization}
Matroids find application throughout combinatorics. The \emph{matroid intersection} problem is defined as finding the largest size common independent set of two matroids. We have seen that the greedy algorithm fully characterizes matroids and Edmonds proved that matroid intersection, can be solved efficiently. This implies that bipartite matchings, common transversals, and tree packing and covering can be computed efficiently too \cite{Schrijver2003Combinatorial-O}. Interestingly, it has been shown that finding matroid intersection of three or more matroids is NP-complete.

%Matroids also find application beyond combinatorics. 
%
%Given such weightings, matroid optimization can tackle a wide range of problems. For example consider we are given a set of points in $\R^n$ and wish to evenly partitions the points such that the variance of their euclidean distances are minimized for each partition. The \emph{minimal variance balanced clustering} problem can be solved via multi criteria matroid optimization \cite{onn-2003-17}. As another interesting example of the broad application of matroid optimization, the statistical problem of experimental design can be solved by multi criteria matroid optimization \cite{berstein-2008}.

%Throughout this thesis we make the modest assumption that our matroids $\M$ are connected.  
Matroids are particularly important in optimization. The greedy algorithm requires weightings on the ground set $S$ of the matroid and it has many applications to optimization. But in real life, multiple weightings can be applied to the ground set along with a balancing function, or concrete notion of optimality over multiple weightings. Consider $d$ weightings $\ve w_1,\ldots,\ve w_d \in \R^n$ on the ground set $[n]$. That is, every $\ve w_i$ assigns a real-value to each element of $[n]$. We let $W \in \R^{d \times n}$ be the matrix with rows $\ve w_1,\ldots, \ve w_d$. We also define $W\Po_\M := \{\, W \ve e_B \mid B \in \B_\M \,\} \subseteq \R^d$, $\# W\Po_\M $ to mean the cardinality of $W \Po_\M$, $\overline{W\Po_\M} := \conv (W\Po_\M)$ and $\vertices(\Po) = \{\text{vertices of } \Po \}$ for a polytope $\Po$. Intuitively one thinks of the rows of $W$ as the set of criteria that (possibly conflicting) parties may bring to a discussion. For example, think of two groups (loggers vs. environmentalists) trying to decide how to choose an optimal spanning tree (a network of logging roads) of a graph. Each group will have different weights.

Hence in \autoref{ch:4thChapterLabel} of this dissertation, we will consider four generalizations of the classical single-criterion (linear-objective) matroid optimization problem (i.e.  $d=1$ which is well-known to be solvable via the greedy method). %The generalizations involve the \emph{multicriteria weighting matrix} $W \in \R^{d \times n}$, a non-linear function $f$, called a \emph{balancing function}, which intuitively is used for breaking ties and reaching decisions comparing all $d$ linear objectives simultaneously.  \\

{\bf Non-Linear Matroid Optimization:} Given a matroid $\M$ on $[n]$ with bases $\B_\M$, $W \in \R^{d \times n}$, and a function $f:\, \R^d \rightarrow \R$, find a base $B \in \B_\M$ such that $f (W\ve e_B) = \min( \, f ( W \ve e_{B'}) \mid B' \in \B \, )$.\\

Two important special cases of our investigations are the following problem.\\

{\bf Convex Matroid Optimization:} Given a matroid $\M$ on $[n]$ with bases $\B$, $W \in \R^{d \times n}$ where $d$ is fixed, and a convex function $f:\, \R^d \rightarrow \R$, find a base $B \in \B$ such that $f (W\ve e_B) = \max f( \, W \ve e_{B'} \mid B' \in \B \, )$. Similarly we consider the minimization problem $f (W\ve e_B) = \min\ f( \, W \ve e_{B'} \mid B' \in \B \, )$\\

{\bf Min-Max Multi-criteria Matroid Optimization:} Given a matroid $\M$ on $[n]$ with bases $\B_\M$, $W \in \R^{d \times n}$, find a base $B \in \B_\M$ such that \begin{equation*} \max_{i=1 \dots d} ((W\ve e_B)_i) = \min( \, \max_{i=1\dots d}(( W \ve e_{B'})_i) \mid B' \in \B_\M \, ).  \end{equation*}

Also, we investigate the following problem.

{\bf Pareto Multi-criteria Matroid Optimization:} Given a matroid $\M$ on $n$-elements with bases $\B_\M$, $W \in \R^{d \times n}$, find a base $B \in \B_\M$ such that\\ $B= \argmin_{Pareto}( \,  ( W \ve e_{B'}) \mid B' \in \B_\M \, )$.\\

In the previous statement, $\min_{Pareto}$ is understood in the sense of Pareto optimality for problems with multiple objective functions, namely, we adopt the convention that for vectors $\ve a,\ve b \in \R^d$, we have $\ve a \leq \ve b$ if and only if $a_i \leq b_i$ for all entries of the vectors. Further, we say that $\ve a<\ve b$ if $\ve a \leq \ve b$ and $\ve a\not=\ve b$. We note that an optimum of the minmax problem will be a Pareto optimum. This is easy to see since if $\ve a \leq \ve b$ then $\max_i(a_i) \leq \max_i(b_i)$. The Pareto multi-criteria matroid optimization problem has been studied by several authors before. For example, Ehrgott \cite{ehrgottmatroid} investigated two optimization problems for matroids with multiple objective functions, and he pioneered a study of Pareto bases via the base-exchange property of matroids.  See \cite{ehrgottsurvey} for a detailed introduction to multicriteria combinatorial optimization.

The matroid optimization problems we consider here have wide applicability. For example, in \cite{berstein-2008-22} the authors consider the ``minimum aberration model fitting'' problem in statistics, which can be reduced to a non-linear matroid optimization problem.  Multi-criteria problems concerning minimum spanning trees of graphs are common in applications (see \cite{ehrgottsurvey} and \cite{Knowles} and references there in).

The Min-Max optimization problem includes the NP-complete partition problem\cite{Garey1979Computers-and-i}, certain multi-processor scheduling problems \cite{Graham1979Optimization-an}, and specific worst-case stochastic optimization problems \cite{Warburton1985Worse-case-anal}.

\begin{example}[Multi-criteria Spanning-Tree Optimization]
Consider a graph $G$ with several linear criteria on the edges,
\begin{itemize}
\item[1.] the first row of $W$ encodes the \emph{fixed installation cost} of each edge of $G$;
\item[2.] the second row of $W$ encodes the \emph{monthly operating cost} of each edge of $G$;
\item[3.] assuming that the edge $j$ fails independently with probability $1-p_j$, then by having the $\log p_j$ as the third row of $W$ (scaled and rounded suitably), $\sum_{j \in T} \log p_j$ captures the \emph{reliability} of the spanning tree $T$ of $G$.
\end{itemize}
It can be difficult for a decision maker to balance these three competing objectives. There are many issues to consider such as the time horizon, repairability, fault tolerance, etc. These issues can be built into a concrete function $f$, for example a weighted norm, or can be thought of as determining a black-box $f$.

\end{example}

Unfortunately, although useful, the problems we are considering are also very difficult in general. For example, multi-criteria matroid optimization is generally NP-hard even in the case when $d=2$ for uniform matroids. It is also NP-hard for the case of spanning trees.  Nevertheless, recently many theoretical and complexity properties about these problems have been explored, yielding good complexity bounds and algorithms under nice assumptions about $W, f$ and $d$. For instance, it has been shown that although convex matroid optimization is NP-complete in general, it is polynomial-time solvable under certain restrictions on $W$ and with fixed $d$.  We refer the reader to the recent series of papers on nonlinear matroid optimization \cite{berstein-2008,berstein-2008-5,onn-2003-17,berstein-2008-22}. % which serve as background for the algorithms and strategies implemented here. The present paper reports on some of the current computational possibilities by comparing various heuristics and algorithms.  
This dissertation explored the experimental performance of their algorithms and proposed some new heuristics and algorithms which performed well in practice. We will see this in \autoref{ch:4thChapterLabel}.

   \chapter[% 
     Volume and Ehrhart Polynomial Computation
   ]{% 
     Volume and Ehrhart Polynomial Computation 
   }%
   \label{ch:2ndChapterLabel}
   \pdfoutput=1

To state our main results recall that given an integer $k > 0$ and a
polytope $\Po \subseteq \R^n$ we define $k \Po := \{\, k \ve \alpha  \mid  \ve \alpha \in \Po
\, \}$ and the function $i(\Po,k) := \#(k\Po \cap \Z^n) $, where we
define $i(\Po,0) :=1$. It is well known that for integral polytopes,
as in the case of matroid polytopes, $i(\Po,k)$ is a polynomial,
called the \emph{Ehrhart polynomial} of $\Po$.  Moreover the leading
coefficient of the Ehrhart polynomial is the \emph{normalized volume}
of $\Po$, where a unit is the volume of the fundamental domain of the
affine lattice spanned by $\Po$ \cite{Stanley1996Combinatorics-a}. Our
first theorem states:

\begin{theorem} \label{volume}
  Let $r$ be a fixed integer.
  Then there exist algorithms whose input data consists of a number $n$ and an evaluation oracle for 
  \begin{enumerate}[\rm(a)]
  \item a rank function~$\varphi$ of a matroid~$M$ on $n$ elements
    satisfying $\varphi(A) \leq r$ for all~$A$, or
  \item an integral polymatroid rank function~$\psi$ satisfying $\psi(A) \leq r$ for
    all~$A$,
  \end{enumerate}
  that compute in time polynomial in~$n$ the Ehrhart polynomial (in
  particular, the volume) of the matroid polytope $\Po(M)$, the independence
  matroid polytope $\Po^\In(M)$, and the polymatroid~$\Po(\psi)$,
  respectively.
\end{theorem}

The computation of volumes is one of the most fundamental geometric
operations and it has been investigated by several authors from the
algorithmic point of view. While there are a few cases for which the volume
can be computed efficiently (e.g., for convex polytopes in fixed
dimension), it has been proved that computing the volume of polytopes
of varying dimension is $\# P$-hard
\cite{dyerfrieze88,brightwellwinkler91,khachiyan93,lawrence91}. Moreover
it was proved that even approximating the volume is hard
\cite{elekes86}. Clearly, computing Ehrhart polynomials is a harder
problem still.  To our knowledge two previously known families of
varying-dimension polytopes for which there is efficient computation of the
volume are simplices or simple polytopes for which the
number of vertices is polynomially bounded (this follows from Lawrence's
volume formula \cite{lawrence91}). Already for simplices
it is at least NP-hard to compute the whole list of coefficients of the 
Ehrhart polynomial, while recently
\cite{barvinok-2006-ehrhart-quasipolynomial}  presented a
polynomial time algorithm to compute any fixed number of the highest
coefficients of the Ehrhart polynomial of a simplex of varying
dimension.  Theorem \ref{volume} provides another interesting family
of varying dimension with volume and Ehrhart polynomial that can be
computed efficiently.  The proof of Theorem \ref{volume}, presented in
Section \ref{sec:compehrhart}, relies on the geometry of tangent cones
at vertices of our polytopes as well as a new refined analysis of the
evaluation of Todd polynomials in the context of the
computational theory of rational generating functions developed by
%%%  \cite{bar,barvinok:99},\cite{barvinok-woods-2003},\cite{barvinok-2006-ehrhart-quasipolynomial},\cite{latte1},\cite{latte2},\cite{Woods:thesis},\cite{verdoolaege-woods-2005}.
\cite{bar,barvinok:99,barvinok-woods-2003,barvinok-2006-ehrhart-quasipolynomial,latte1,latte2,Woods:thesis,verdoolaege-woods-2005}. A nice introduction to these topics can be found at \cite{Beck2007Computing}

\section{Computing the Ehrhart Polynomials} \label{sec:compehrhart}

\section{Preliminaries on Rational Generating Functions.}

Generating functions are crucial to proving our main results. For a good
reference for the basic concepts used here see \cite{Beck2007Computing}.  
Let
$\Po \subseteq \R^n$ be a rational polyhedron. The \emph{multivariate
  generating function} of $\Po$ is defined as the formal Laurent series in $\Z[[z_1,\ldots,z_n,z_1^{-1},\ldots,z_n^{-1}]]$
\begin{equation*}
\tilde g_\Po ( \mathbf{z} ) = \sum_{\ve \alpha \in \Po \cap \Z^n} \mathbf{z}^\alpha ,
\end{equation*}
where we use the multi-exponent notation $\ve z^\vealpha =
\prod_{i=1}^n z_i^{\alpha_i}$.  If $\Po$ is bounded, $\tilde g_\Po$ is a
Laurent polynomial, which we consider as a rational function~$g_\Po$.
If $\Po$ is not bounded but is pointed (i.e., $\Po$ does not contain a
straight line), there is a non-empty open subset $U\subseteq\C^n$ such
that the series converges absolutely and uniformly on every compact
subset of~$U$ to a rational function~$g_\Po$ (see \cite{barvinok:99} and
references therein).  If $\Po$ contains a straight line, we set $g_\Po =
0$.  The rational function $g_\Po\in\Q(z_1,\dots,z_n)$ defined in this
way is called the \emph{multivariate rational generating function}
of~$\Po\cap\Z^n$.  Barvinok \cite{bar} proved that in polynomial time, when the
dimension of a polyhedron is fixed, $g_\Po$ can be represented
as a short sum of rational functions 
\begin{equation*}
%%%  g_\Po(\ve z) =  \sum_{i \in I} \epsilon_i \frac{\sum_{l \in B_i} \ve z^{\ve a_{il}}}{\prod_{j=1}^{n} (1 - \ve z^{\ve b_{ij}} )},
 g_\Po(\ve z) =  \sum_{i \in I} \epsilon_i \frac{\ve z^{\ve a_{i}}}{\prod_{j=1}^{n} (1 - \ve z^{\ve b_{ij}} )},
\end{equation*}
where $\epsilon_i \in \{-1,1\}$.

Our first contribution is to show that in
the case of matroid polytopes of fixed rank, this still holds even
when their dimension grows.
Let $\ve v$ be a vertex of
$\Po$. Define the \emph{tangent cone} or \emph{supporting cone} of
$\ve v$ to be
\begin{equation*}
    \Co_\Po(\ve v) := \left\{\, \ve v + \ve w \; | \; \ve v + \epsilon \ve w \in \Po \text{ for some } \epsilon > 0 \, \right\}.
\end{equation*}
We rely on the following result, which connects the rational generating
function of a rational polyhedron to those of the tangent cones of its
vertices. 
%%%   The following lemma is well-known as Brion's Theorem; it is due to \cite{Brion88} and, independently, Lawrence (1991) \cite{lawrence91-2}.  
%%%   The following lemma is due to \cite{Brion88} and independently \cite{lawrence91-2}.
This result was independently discovered by Brion \cite{Brion88} and Lawrence \cite{lawrence91-2}.
A proof can also be found in
\cite{barvinok:99} and \cite{beck-haase-sottile:theorema}.

\begin{lemma}[Brion--Lawrence's Theorem] \label{lem:bl}
  Let $\Po$ be a rational polyhedron
  and $V(\Po)$ be the set of vertices of~$\Po$. Then,
  \begin{equation*}
    g_\Po(\mathbf{z}) = \sum_{\ve v \in V(\Po)} g_{\Co_\Po(\ve v)}(\mathbf{z}),
  \end{equation*}
  where $\Co_\Po(\ve v)$ is the tangent cone of~$\ve v$.
\end{lemma}
Thus, we can write the multivariate generating function of $\Po$ by writing
all multivariate generating functions of all the tangent cones of the vertices
of~$\Po$.
Moreover, the map assigning to a rational polyhedron~$\Po$ its multivariate
rational generating function~$g_\Po(\ve z)$ is a \emph{valuation}, i.e., a
finitely additive measure, so it satisfies the equation
\begin{equation*}
  g_{\Po_1\cup \Po_2}(\ve z) = g_{\Po_1}(\ve z) + g_{\Po_2}(\ve z) - g_{\Po_1\cap \Po_2}(\ve z),
\end{equation*}
for arbitrary rational polytopes $\Po_1$~and~$\Po_2$,
the so-called \emph{inclusion--exclusion principle}.  
This allows to break a polyhedron~$\Po$ into pieces $\Po_1$~and~$\Po_2$ and to compute
the multivariate rational generating functions for the pieces (and their
intersection) separately in order to get the generating function~$g_{\Po}$. 
More generally, let us denote by~$[\Po]$ the \emph{indicator function} of~$\Po$,
i.e., the function
\begin{equation*}
  [\Po]\colon \R^n\to\R, \quad
  [\Po](\ve x) = 
  \begin{cases}
    1 & \text{if $\ve x\in \Po$} \\
    0 & \text{otherwise}.
  \end{cases}
\end{equation*}
Let
\begin{math}
%%  \label{eq:indicator-identity}
  \sum_{i\in I} \epsilon_i [\Po_i] = 0
\end{math}
be an arbitrary linear identity of indicator functions of rational polyhedra
(with rational coefficients~$\epsilon_i$);
the valuation property now implies that it carries over to a linear identity 
\begin{math}
%%  \label{eq:ratgenfun-identity}
  \sum_{i\in I} \epsilon_i\, g_{\Po_i}(\ve z) = 0
\end{math}
of rational generating functions.  

Now let $\Co$ be one of the tangent cones of~$\Po$, and let $\Tr$ be a
triangulation of~$\Co$, given by its simplicial cones of maximal dimension.  Let
$\hat \Tr$ denote the set of all (lower-dimensional) intersecting proper faces
of the cones~$\Co_i\in\Tr$.  Then we can assign an integer
coefficient~$\epsilon_i$ to every cone $\Co_i \in \hat\Tr$, such that the
following identity holds:
%%% \mkoeppesays{Is this $\epsilon_i$ called the M\"obius number?}
\begin{equation*}
  \label{eq:indicator-identity-triang}
  [\Co] = \sum_{\Co_i\in \Tr} [\Co_i] + \sum_{\Co_i\in \hat \Tr} \epsilon_i\, [\Co_i].
\end{equation*}
This identity immediately carries over to an identity of multivariate
rational generating functions,
\begin{equation}
  \label{eq:ratgenfun-identity-triang}
  g_\Co(\ve z) = \sum_{\Co_i\in \Tr} g_{\Co_i}(\ve z) 
  + \sum_{\Co_i\in \hat \Tr} \epsilon_i\, g_{\Co_i}(\ve z).
\end{equation}
Hence, the problem of computing rational generating functions of a polyhedron is reduced to
the case of simplicial cones.
\bigbreak

In the following (\autoref{subsec:ontangent}, \autoref{subsec:polymatroids}),
we study the tangent cones of the matroid polytopes and polymatroids and
introduce algorithms that construct triangulations for them.  Then, in
\autoref{subsec:short}, we construct a short multivariate generating function
using an \emph{efficient} variant of
identity~\eqref{eq:ratgenfun-identity-triang} and Brion--Lawrence's Theorem.
Finally, in \autoref{subsec:specialize}, we compute the Ehrhart polynomial.

        %%%%%%%%%%%%%%%%%%%%%%%%%%%%%%%%%%%%%%%%%%%%%%%%%%%%%%%%%%%%%%%%%%%
        %%%%%%       ON THE TANGENT CONES OF MATROID POLYTOPES      %%%%%%%
        %%%%%%%%%%%%%%%%%%%%%%%%%%%%%%%%%%%%%%%%%%%%%%%%%%%%%%%%%%%%%%%%%%%

\section{On the Tangent Cones of Matroid Polytopes} \label{subsec:ontangent}

\begin{figure}[!htb]
    \centerline{
    %\scalebox{0.15}{
    \ifpdf
    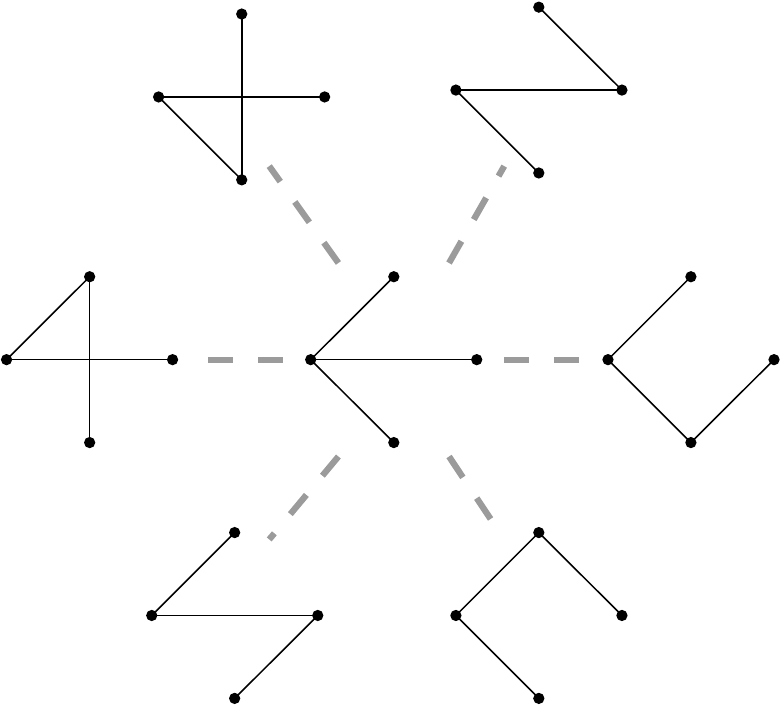
    \else
    \input{K4adj.pstex_t}
    \fi
    }
    %}
    \caption{${\{2, 3, 5\}}$, ${\{2, 3, 4\}}$, ${\{1, 3, 6\}}$, ${\{1, 3, 4\}}$, ${\{1, 2, 6\}}$ and ${\{1, 2, 5\}}$ are spanning trees of $K_4$ that differ from ${\{1, 2, 3\}}$ by adding one edge and removing one edge.}\label{fig:K4adj}
\end{figure}

Our goal is to compute the multivariate generating function of matroid
polytopes and independence matroid polytopes with fixed rank (later,
in \autoref{subsec:polymatroids}, we
will deal with the case of polymatroids), and to do this we will use a
crucial property of adjacent vertices.  To illustrate our techniques we will use a running example throughout this section.
%%%   \begin{figure}[!htb]    \label{fig:K4}
%%%       \centerline{
%%%           \scalebox{0.4}{
%%%       \inputfig{K4}
%%%           }
%%%       }
%%%       \caption{The complete graph on $4$ vertices, $\ve K_4$}
%%%   \end{figure}
\begin{example}[Matroid on $K_4$] \label{ex:k4}
    Let $K_4$ be the complete graph on four vertices. Label the ${4 \choose 2} = 6$ edges with $\{1,\ldots,6\}$ as in Figure \ref{fig:K4adj}. Every graph induces a matroid on its edges where the bases are all spanning trees (spanning forests for disconnected graphs) \cite{Welsh1976Matroid-Theory}. Let $M(K_4)$ be the matroid on the elements $\{1,\ldots,6\}$ with bases as all spanning trees of $K_4$. The rank of $M(K_4)$ is the size of any spanning tree of $K_4$, thus the rank of $M(K_4)$ is $3$. The $16$ bases of $M(K_4)$ are: 
$\{3, 5, 6\}$, $\{3, 4, 6\}$, $\{3, 4, 5\}$, $\{2, 5, 6\}$, $\{2, 4, 6\}$, $\{2, 4, 5\}$, $\{2, 3, 5\}$, $\{2, 3, 4\}$, $\{1, 5, 6\}$, $\{1, 4, 6\}$, $\{1, 4, 5\}$, $\{1, 3, 6\}$, $\{1, 3, 4\}$, $\{1, 2, 6\}$, $\{1, 2, 5\}$, $\{1, 2, 3\}$.
\end{example}

Let $M$ be a matroid on $n$ elements with fixed rank $r$. Then the number of
vertices of $\Po(M)$ is a polynomial in $n$ of degree $r$. We can see this since the number of
vertices is equal to the number of bases of $M$, and the number of bases is
bounded by $n \choose r$, a polynomial in $n$ of degree $r$. Clearly the number of
vertices of $\Po^\In (M)$ is also polynomial in~$n$.  It is also clear that, 
when the rank~$r$ is fixed, 
all vertices of either polytope can be enumerated in polynomial time in~$n$,
even when the matroid is only presented by an evaluation oracle for its rank
function~$\varphi$.

Throughout this section we shall discuss polyhedral cones $\Co$ with
extremal rays $\{\ve r_1,\ldots,\ve r_l\}$ such that 
%%%  \begin{equation*}
%%%      \ve r_k \; \in \; R_A := \{ \, \ve e_i - \ve e_j, \, \ve e_i, \, -\ve e_j \; | \;  i \in [n] \setminus A, \; j \in A \, \} \quad \text{for } k=1,\ldots,l
%%%  \end{equation*}
\begin{equation*}
    \ve r_k \; \in \; R_A := \{ \, \ve e_i - \ve e_j, \, \ve e_i, \, -\ve e_j \; | \;  i \in [n], \; j \in A \, \} \quad \text{for } k=1,\ldots,l
\end{equation*}
for some $A \subseteq [n]$. We will refer to $R_A$ as the
\emph{elementary set of $A$}. Note that by \autoref{lem:adj}
the rays of a tangent cone at a vertex $\ve e_A$ (corresponding to a set
$A\subseteq[n]$) of a matroid polytope or an independence matroid polytope form an
elementary set of $A$. Due to convexity and the
assumption that $\ve r_k$ are extremal, for each $i \in [n]$ and $j \in A$ at
most two of the three vectors $\ve e_i - \ve e_j, \ve e_i, - \ve e_j$ are
extremal rays $\ve r_k$ of $\Co$.  This implies by construction, that
    considering all pairs $\ve e_i - \ve e_j$ and $\ve e_i$ or $-\ve
    e_j$, the number of generators $\ve r_k$ of $\Co$ is bounded by
\begin{equation} \label{eq:genbnd}
%%% (n- |A|)|A|  + 2\max \{|A|,|[n] \setminus A| \}.
n|A| + n + |A|.
\end{equation}

Recall a cone is \emph{simple} if it is generated by linearly independent vectors and
it is \emph{unimodular} if its fundamental parallelepiped contains only one lattice point \cite{1997Handbook-of-dis}.  A triangulation of $\Co$ is
\emph{unimodular} if it is a polyhedral subdivision such that each sub-cone is
unimodular.

\begin{example}[Matroid on $K_4$]
    The vertices $\ve e_{\{2, 3, 5\}}$, $\ve e_{\{2, 3, 4\}}$, $\ve e_{\{1, 3, 6\}}$, $\ve e_{\{1, 3, 4\}}$, $\ve e_{\{1, 2, 6\}}$ and $\ve e_{\{1, 2, 5\}}$ are all adjacent to the vertex $\ve e_{\{1, 2, 3\}}$, see \autoref{fig:K4adj}. Moreover, the tangent cone $\Co_{\Po(M(\ve K_4))}(\ve e_{\{1, 2, 3\}})$ is generated by the differences of these vertices with $\ve e_{\{1, 2, 3\}}$:
\begin{multline*}
    \Co_{\Po(M(\ve K_4))}(\ve e_{\{1, 2, 3\}}) = \ve e_{\{1, 2, 3\}} + \cone \left\{\ve e_{\{2, 3, 5\}}-\ve e_{\{1, 2, 3\}}, \ve e_{\{2, 3, 4\}}-\ve e_{\{1, 2, 3\}}, \right. \\
           \left. \ve e_{\{1, 3, 4\}}-\ve e_{\{1, 2, 3\}},\ve e_{\{1, 2, 6\}}-\ve e_{\{1, 2, 3\}},\ve e_{\{1, 2, 5\}}-\ve e_{\{1, 2, 3\}} \right\}
%%%    \Co_{\Po(M(\ve K_4))}(\ve e_{\{1, 2, 3\}}) = \cone \left( \ve e_{\{2, 3, 5\}}, \ve e_{\{2, 3, 4\}}, \ve e_{\{1, 3, 4\}},\ve e_{\{1, 2, 6\}},\ve e_{\{1, 2, 5\}} \right) - \ve e_{\{1, 2, 3\}}
\end{multline*}
\end{example}

\begin{lemma} \label{lem:eluni}
    Let $\Co \subseteq \R^n$ be a cone generated by $p$ extremal rays $\{\ve
    r_1,\ldots,\ve r_p\} \subseteq R_A$ where $R_A$ is an elementary set of
    some $A \subseteq [n]$. Every triangulation of $\Co$ is unimodular.
\end{lemma}

\begin{proof}
%%%   The $\Co$ extremal ray generators $\{\ve r_1,\ldots,\ve r_p\}$ are of the form $ \ve r_k \; \in \; R_A $ for some $A \subseteq [n]$.
Without loss of generality, we can assume $\{\ve r_1,\ldots,\ve r_l\}$
are generators of the form $\ve e_i - \ve e_j$ and $\{\ve r_{l+1},\ldots,\ve r_p\}$ are
generators of the form $\ve e_i$ or $-\ve e_j$ for the cone $\Co$. 

It is easy to see that the matrix $\tilde T_\Co := [\ve r_1,\ldots,\ve r_l]$ is totally
unimodular. Let $G_\Co$ be a directed graph with vertex set $[n]$ and
an edge from vertex $i$ to $j$ if $\ve r_k = \ve e_i - \ve e_j$ is an extremal ray of
$\Co$. We can see that $G_\Co$ is a subgraph of the complete directed graph $K_n$ with two arcs between each pair of vertices; one for each direction.
Since $\tilde T_\Co := [\ve r_1,\ldots,\ve r_l]$ is the incidence matrix of the graph $G_\Co$,
      it is totally unimodular
      \cite[Ch.~19,~Ex.~2]{Schrijver1986Theory-of-linea}, i.e., every subdeterminant is $0$, $1$ or $-1$
      \cite[Ch.~19, Thm.~9]{Schrijver1986Theory-of-linea}.  Therefore $T_\Co :=
      [\ve r_1,\ldots,\ve r_l,\ve r_{l+1},\ldots,\ve r_p]$ is totally unimodular since
      augmenting $\tilde T_\Co$ by a vector $\ve e_i$ or $-\ve e_j$ preserves this
      subdeterminant property: for any submatrix containing part of a vector
      $\ve e_i$ or $-\ve e_j$ perform the cofactor expansion down the vector $\ve e_i$ or
      $-\ve e_j$ when calculating the determinant.  

Since $T_\Co$ is totally unimodular, each basis of $T_\Co$ generates the entire integer
lattice $\Z^n \cap \lin(\Co)$ and hence every simplicial cone of a triangulation has
normalized volume $1$.
\end{proof}

\begin{lemma} \label{lem:aug}
    Let $\Co \subseteq \R^n$ be a cone generated by $l$ extremal
    rays $\{\ve r_1,\ldots,\ve r_l\} \subseteq R_A$ where $R_A$ is an elementary set of some $A \subseteq [n]$, where $\dim(\Co) < n$. The extremal rays
    $\{\ve r_1,\ldots,\ve r_l\}$ can be augmented by a vector $\ve{\tilde r}$ such that
    $\dim ( \cone \{ \ve r_1,\ldots,\ve r_l,\ve{\tilde r}\}) = \dim(\Co) +1$,
    the vectors
    $\ve r_1,\ldots,\ve r_l,\ve{\tilde r}$ are all extremal,
    and $\ve{\tilde r} \in R_A$.
\end{lemma}
\begin{proof}
    It follows from convexity that at most two of $\ve e_i-\ve e_j$, $\ve e_i$ or
    $-\ve e_j$ are extremal generators of $\Co$ for $i \in [n]$ and $j \in A$.
    There are at least $n$ possible extremal ray generators, considering two of
    $\ve e_i-\ve e_j$, $\ve e_i$ or $-\ve e_j$ for each $i \in [n]$ and $j \in A$.
    Moreover, all these pairs span $\R^n$. Thus by the basis augmentation theorem
    of linear algebra, there exists a vector $\ve{\tilde r}$ such that $\dim (
            \cone \{ \ve r_1,\ldots,\ve r_l,\ve{\tilde r}\} )  = \dim(\Co) +1$ and
    $\ve r_1,\ldots,\ve r_l,\ve{\tilde r}$ are all extremal. 
\end{proof}

\begin{lemma} \label{lem:elsimpbnd}
%%%    Let $\Co \subseteq \R^n$ be a cone generated by $l$ extremal rays $\{\ve
%%%    r_1,\ldots,\ve r_l\} \subseteq R_A$ where $R_A$ is an elementary set of some fixed $A \subseteq [n]$.  Any triangulation of $\conv(\{\ve 0, \ve
%%%            r_1,\ldots,\ve r_l\})$ has at most a polynomial in $n$ number of
%%%    top-dimensional simplices.
    Let $r$ be a fixed integer, $n$ be an integer, $A \subseteq [n]$ with $|A| \leq r$ and
    let $\Co \subseteq \R^n$ be a cone generated by $l$ extremal rays $\{\ve
    r_1,\ldots,\ve r_l\} \subseteq R_A$ where $R_A$ is an elementary set of $A$. Then any triangulation of $\conv(\{\ve 0, \ve
            r_1,\ldots,\ve r_l\})$ has at most a polynomial in $n$ number of
    top-dimensional simplices.
\end{lemma}
\begin{proof}

Assume $\dim(\Co) = n$. Later, we will show how to remove this restriction.
We can see that
$    \conv \{\ve  0, \ve r_1,\ldots,\ve r_l \} \; \subseteq \;  [-1,1]^{A} \times \tilde \Delta_{[n] \setminus A}$
where
\begin{align*}
    [-1,1]^{A} & :=  \{ \, \ve x \in \R^A \mid |\ve x_j| \leq 1 \ \ j \in A \, \}  \subseteq \R^A  \\
    \tilde \Delta_{[n] \setminus A} & := \conv ( \{ \, \ve e_i \; | \; i \in [n] \setminus A \, \} \cup \{\ve 0\} )  \subseteq \R^{[n] \setminus A}.
\end{align*}
The volume of a $d$-simplex $\conv\{\ve v_0,\ldots,\ve v_d\}$ is \cite{1997Handbook-of-dis}
\begin{equation}
    \label{eq:simpvol}
    \frac{1}{d!} \left| \det \left( \begin{array}{ccc} \ve v_0 & \cdots & \ve v_d \\ \ve 1 & \cdots & \ve 1 \end{array} \right) \right| .
\end{equation}
%%%  Thus the $m$-volume of $\conv (\{\ve 0,\ve e_1,\ldots,\ve e_m \}) = \frac{1}{m!}$ and the $|A|$-volume of $[-1,1]^{|A|}$ is $2^{|A|}$.
Thus the $(n-|A|)$-volume of $\tilde \Delta_{[n] \setminus A}$ is $ 
\frac{1}{(n-|A|)!}$ and the $|A|$-volume of $[-1,1]^{A}$ is $2^{|A|}$.
Therefore
\begin{equation*}
    \vol \big( [-1,1]^{A} \times \tilde \Delta_{[n] \setminus A} \big) =
    2^{|A|} \frac{1}{(n-|A|)!} = \frac1{n!} 2^{|A|} n(n-1) \cdots (n-|A| +1) .
\end{equation*}
It is also a fact that any integral $n$-simplex has $n$-volume bounded below by
$\frac{1}{n!}$, using the simplex volume equation \eqref{eq:simpvol}. Therefore
any triangulation of $\conv \{ \ve 0,\ve r_1,\ldots,\ve r_l \}$ has at most 
\begin{equation*}
  2^{|A|} n(n-1) \cdots (n-|A| +1) \leq
  2^{r} n(n-1) \cdots (n-r +1)
\end{equation*}
full-dimensional simplices, a polynomial function in~$n$ of degree~$r$.

Let $d_\Co := n-\dim( \Co )$. If $\dim(\Co) < n$, then by Lemma \ref{lem:aug}, $\{\ve r_1,\ldots,\ve r_l\}$ can
be augmented with vectors $\{\ve{\tilde r}_1, \ldots, \ve{\tilde r}_{d_\Co} \}$ where $ \ve{\tilde
r}_k \in R_A$ for $A$ above, such that $\dim( \cone \{ \ve r_1,\ldots,\ve r_l,\allowbreak \ve{\tilde r}_1, \ldots, \ve{\tilde r}_{d_\Co} \} ) = n$ and $\{
\ve r_1,\ldots,\ve r_l,\ve{\tilde r}_1, \ldots, \ve{\tilde r}_{d_\Co} \}$ are extremal. Moreover,
\begin{equation*}
    \begin{array}{c}
    \dim ( \conv \{\ve 0,\ve r_1,\ldots,\ve r_l\} ) < \dim ( \conv \{\ve 0, \ve r_1,\ldots,\ve r_l,\ve{\tilde r}_1\} ) < \cdots  \\
    < \dim ( \conv \{\ve 0 ,\ve  r_1,\ldots,\ve  r_l,\ve{\tilde r}_1,\ldots,
    \ve{\tilde r}_{d_\Co -1}\} ) 
    < \dim ( \conv \{\ve 0,\ve r_1,\ldots,\ve r_l,\ve{\tilde r}_1,\ldots, \ve{\tilde r}_{d_\Co}\} ),
    \end{array}
\end{equation*}
that is, $\ve{\tilde r}_k \notin \affine \{\ve 0, \ve r_1,\ldots,\ve r_l, \ve{\tilde r}_1, \ldots, \ve{\tilde r}_{k-1}\}$ for $1 \leq k \leq d_\Co$. 

\begin{figure}[!htb]
    \centerline{ \ifpdf
    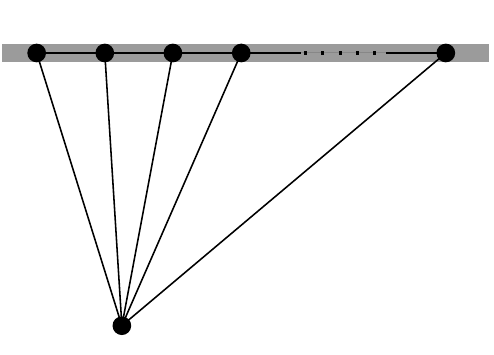
    \else
    \input{lemnsimp.pstex_t}
    \fi }
    \caption{If $\ve{\tilde r}_k$ is not contained in the affine span of
      $\{\ve 0, \ve r_1,\ldots,\ve r_p, \ve{\tilde r}_1,\ldots, \ve{\tilde
        r}_{k-1}\}$ then every full-dimensional simplex must contain
      $\ve{\tilde r}_k$.}
    \label{fig:lemnsimp}
\end{figure}

Since $\ve{\tilde r}_k \notin \affine \{ \ve 0, \ve r_1,\ldots,\ve r_l,
        \ve{\tilde r}_1, \ldots, \ve{\tilde r}_{k-1} \} $, any full-dimensional
simplex in a triangulation of $\conv \{ \ve 0, \ve r_1,\ldots,\ve r_l,
        \ve{\tilde r}_1, \ldots, \ve{\tilde r}_k \} $ must contain $\ve{\tilde
    r}_k$, see Figure \ref{fig:lemnsimp}. If not, then there exists a
    top-dimensional simplex using the points $\{ \ve 0, \ve r_1,\ldots,\ve
    r_l, \allowbreak \ve{\tilde r}_1, \ldots, \ve{\tilde r}_{k-1} \}$, but we know all
    these points lie in a subspace of one less dimension, a contradiction. Therefore, a bound on
    the number of simplices in a triangulation of $\conv \{ \ve 0, \ve
            r_1,\ldots,\ve r_l, \ve{\tilde r}_1, \ldots, \ve{\tilde r}_k \}$ is
    a bound on that of $\conv \{ \ve 0, \ve r_1,\ldots,\ve r_l,\allowbreak \ve{\tilde r}_1,
            \ldots, \ve{\tilde r}_{k-1} \}$.  

Thus, if $\dim (\Co) < n$ we can augment $\Co$ by vectors $\ve{\tilde
r}_1,\ldots, \ve{\tilde r}_{d_\Co}$ so that the cone $\tilde \Co :=
\cone \{ \ve r_1,\ldots,\ve r_p,\ve{\tilde r}_1,\ldots, \ve{\tilde
r}_{d_\Co} \}$ is of dimension $n$ and $ \ve r_l, \ve{\tilde r}_k \;
\in \; R_A$ for $A$ above. We proved any triangulation of $\conv \{
\ve 0, \ve r_1,\ldots,\ve r_p,\ve{\tilde r}_1,\ldots, \ve{\tilde r}_{d_\Co} \}
$ has at most polynomially many full-dimensional $n$-simplices, which
implies that any triangulation of $\conv \{ \ve 0, \ve r_1,\ldots,\ve
r_p \}$ has at most polynomially many top-dimensional
simplices due to the construction of the generators
$\ve{\tilde r}_k$.
\end{proof}

We have shown that for a cone $\Co$ generated by an elementary set of
extremal rays $\{\ve r_1,\ldots,\ve r_l\} \subseteq R_A$ for some $A
\subseteq [n]$, any triangulation of $\conv\{\ve 0,\ve r_1,\ldots,\ve
r_l\}$ has at most polynomially many simplices. What we need next is
an efficient method to compute some triangulation of $\conv\{\ve
0,\ve r_1,\ldots,\ve r_l\}$. We will show that the placing
triangulation is a suitable candidate.

Let $\Po \subseteq \R^n$ be a polytope of dimension $n$ and $\Delta$
be a facet of $\Po$ and $\ve v \in \R^n$. There exists a unique
hyperplane $H$ containing $\Delta$ and $\Po$ is contained in one of
the closed sides of $H$, call it $H^+$. If $\ve v$ is contained in the
interior of $H^-$, the other closed halfspace defined by $H$, then
$\Delta$ is \emph{visible} from $\ve v$ (see chapter 14.2 in
\cite{1997Handbook-of-dis}). The well-known \emph{placing
triangulation} is given by an algorithm where a point is iteratively added to an
intermediate triangulation by determining which facets are visible to
the new point
\cite{1997Handbook-of-dis, De-Loera2006Triangulations}. We recall now how
to determine if a facet is visible to a vertex in polynomial time.

\begin{figure}[!htb]    
    \centering
    \hfil(a)\subfigure{\scalebox{0.8}{\ifpdf
    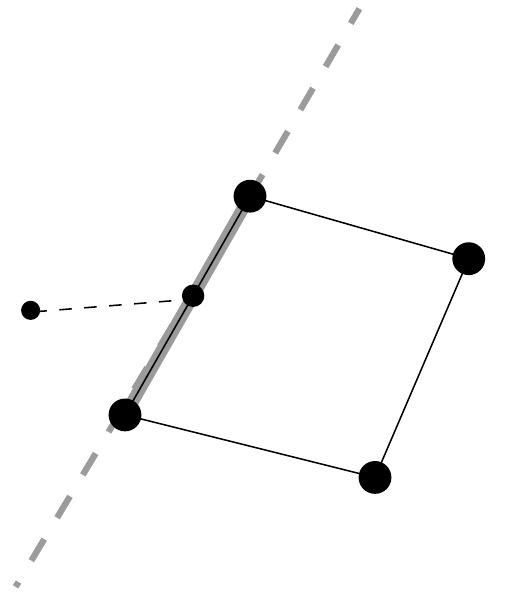
    \else
    \input{visfacet1.pstex_t}
    \fi}}
    \hfil(b)\subfigure{\scalebox{0.8}{\ifpdf
    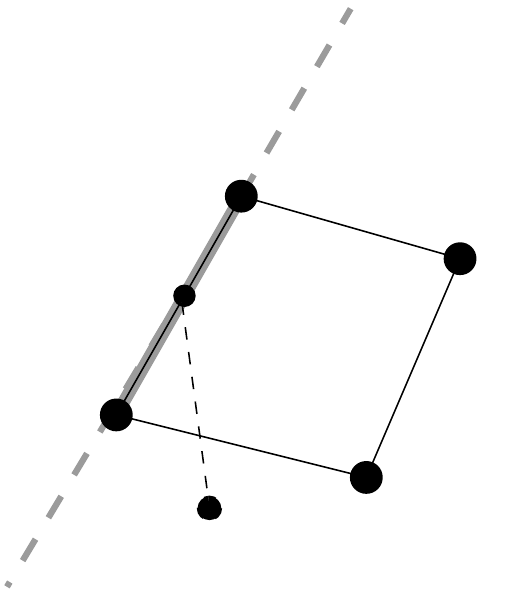
    \else
    \input{visfacet2.pstex_t}
    \fi}}
    %\hfil(a)\subfigure{\inputfig{visfacet1}}
    %\hfil(b)\subfigure{\inputfig{visfacet2}}
    \caption{Figure (a) shows $\Delta$ visible to $\ve v$. Figure (b) shows
        $\Delta$ {\bf not} visible to $\ve v$.}
\end{figure}

\begin{lemma} \label{lem:vislp}
Let $\Po \subseteq \R^n$ be a polytope given by $t$ vertices
$\{\ve v^1,\ldots,\ve v^t\} \subseteq \R^n$ and $\Delta \subseteq \Po$ be a facet of $\Po$
given by $q$ vertices $\{\ve{\tilde v}^1,\ldots, \ve{\tilde v}^q \} \subseteq
\{\ve v^1,\ldots,\ve v^t\} $. If $\ve v \in \R^n$ where $\ve v \notin \Po$ then deciding if
$\Delta$ is visible to $\ve v$ can be done in polynomial time in the input
$\{\ve{\tilde v}^1,\ldots, \ve{\tilde v}^q \}$, $\{\ve v^1,\ldots,\ve v^t\}$ and $\ve v$.
\end{lemma}

\begin{proof}
Let $  \ve  z := \frac{1}{q}
\sum_{i=1}^q  \ve{\tilde v}^i $ so that $\ve z \in \interior ( \Delta )$. We consider the linear program:
\begin{equation} \label{vislp}
    \begin{split} 
    \Bigg\{ \, \left( \begin{array}{c} \ve x \\ \ve y \\ \lambda \end{array} \right) \in \R^{n+t+1} \;\; | \;\;  \ve x = \sum_{i=1}^t \ve v^i \ve y_i, \: & \ve y \geq \ve 0, \: \sum_{i=1}^t \ve y_i=1,    \\
    & \: 0 \leq \lambda < 1, \: \lambda \ve v + (1-\lambda) \ve z = \ve x \, \Bigg\} .  
    \end{split}
\end{equation}

If \eqref{vislp} has a solution then there exists a point $\ve{\bar x} \in \Po$
between the facet $ \Delta$ and $\ve v$, hence $ \Delta$ is not visible from $\ve v$. If
\eqref{vislp} does not have a solution, then there are no points of $
\Po$ between $\ve v$ and $\Delta$, hence $\Delta$ is visible from $\ve v$
(see Lemma 4.2.1 in \cite{De-Loera2006Triangulations}). It is well known that
a strict inequality, such as the one in \eqref{vislp}, can be handled by an
equivalent linear program which has only one additional variable. Determining
if \eqref{vislp} has a solution can be done in polynomial time in the input
\cite{Schrijver1986Theory-of-linea}.
\end{proof}

The placing triangulation is obtained by incrementally adding one point at a time, connecting the new point to the current triangulation. More precisely:
%%%  The table environment created problems with the algorithms placement.
%%%  \begin{table}[htbp]
%%%      \caption{The Polytime Placing Triangulation}
%%%      \label{tab:ppt}
\vskip 10pt
\begin{minipage}{\linewidth}
\begin{algorithm}[The Placing Triangulation \cite{1997Handbook-of-dis, De-Loera2006Triangulations}]
    \mbox{}
    \vskip .25cm
    \label{alg:ept}
    \begin{center}
        \begin{tabular}{l}
            %\toprule
            \begin{minipage}{.85\linewidth}
                \begin{algorithmic}
                    %%%   \INPUT A set of ordered points $\{ v_1,\ldots,v_t \} \in \R^n$. 
                    \item[Input:] A set of ordered points $\{ \ve v_1,\ldots,\ve v_t \} \in \R^n$.
                    %%%   \OUTPUT A triangulation $\Tr$ of $\{ v_1,\ldots,v_t \}$
                    \item[Output:] A triangulation $\Tr$ of $\{ \ve v_1,\ldots,\ve v_t \}$
                    \STATE $\Tr := \{\{\ve v_1\}\}$.
                    \FOR{ each $\ve v_i \in \{\ve v_2,\ldots,\ve v_t\}$} \label{alg:ept:l4}
                        \STATE Let $B \in \Tr$.
                        \STATE $P_i := \{\ve v_1,\ldots,\ve v_{i-1}\}$
                        %\STATE Create and run the linear program \eqref{vislp} with $(P_i,B,v_i)$ to decide visibility of $B$ to $v_i$.
                        \IF{ $\ve v_i \notin \affine(P_i)$} \label{alg:ept:aff} 
                            \STATE $\Tr' := \emptyset$
                            \FOR{ each $D \in \Tr$} \label{alg:ept:l9}
                                \STATE $\Tr' := \Tr' \cup \{D \cup \{\ve v_i\}\}$.
                            \ENDFOR
                        %\IF{ $\Tr' = \emptyset$}
                        \ELSE
                            \FOR{ each $B \in \Tr$ and each $(|B|-1)$-subset $C$ of $B$} \label{alg:ept:l12}
                                \STATE Create and solve the linear program \eqref{vislp} with $(P_i,C,\ve v_i)$ to decide visibility of $C$ to $\ve v_i$. \label{alg:ept:l13}
                                \IF{$C$ is visible to $\ve v_i$}
                                    \STATE $\Tr' := \Tr' \cup \{C \cup \{\ve v_i\}\}$
                                \ENDIF
                            \ENDFOR
                        \ENDIF
                        \STATE $\Tr := \Tr'$
                    \ENDFOR
                    \RETURN $\Tr$
                \end{algorithmic}
            \end{minipage}\\
            %\bottomrule
        \end{tabular}
    \end{center}
\end{algorithm}
\end{minipage}
\vskip 10pt
%%%   \end{table}

Indeed, Algorithm \ref{alg:ept} returns a triangulation
\cite{1997Handbook-of-dis}. We will show that for certain input (a point set
corresponding to a vertex cone of a matroid polytope or independence matroid polytope), it runs in
polynomial time. We remark that there are exponentially, in $n$, many lower
dimensional simplices in any given triangulation. But, it is important to note
that only the highest dimensional simplices are listed in an intermediate
triangulation (and thus the final triangulation) in the placing triangulation
algorithm.

%% The following theorem takes as input the rays generating a cone, considers the
%% rays as a point set and triangulates them.

\begin{theorem} \label{thm:ept}
Let $r$ be a fixed integer, $n$ be an integer, $A \subseteq [n]$ with $|A|
\leq r$, and let 
$\{\ve r_1,\ldots,\ve r_l\} \subseteq R_A$. Then the placing triangulation
(Algorithm \ref{alg:ept}) with input $\{\ve 0, \ve r_1,\ldots,\ve r_l\}$ runs
in polynomial time. 
%%%  Let $\Co \subseteq \R^n$ be a tangent cone of $\Po(M)$ generated by $l$ extremal rays $\{\ve r_1,\ldots,\ve r_l\}$. The placing triangulation (Algorithm \ref{alg:ept}) with input $\{\ve 0, \ve r_1,\ldots,\ve r_l\}$, runs in polynomial time.
\end{theorem}
\begin{proof}
%%%   Since $\Co$ is a tangent cone of $\Po(M)$ then $\{\ve r_1,\ldots,\ve r_l\} \subseteq R_A$ for some $A \subseteq [n]$
By Equation \eqref{eq:genbnd} there is only a polynomial, in $n$, number of
extremal rays $\{\ve r_1,\ldots,\ve r_l\}$. Thus, the {\bf for}
statement on line \ref{alg:ept:l4} repeats a polynomial number of
times. Step \ref{alg:ept:aff} can be done in polynomial time by solving the linear
equation $[\ve r_1,\ldots,\ve r_{i-1}] \ve x = \ve v_i$.  %%%%  \davecomment{(Double check. Size of matrix is nxp(n), where p is a polynomial in n)}

The {\bf for} statement on line \ref{alg:ept:l9} repeats for every simplex $D$ in the
triangulation $\Tr$, and the number of simplices in $\Tr$ is bounded by the
number of simplices in the final triangulation. By Lemma \ref{lem:elsimpbnd}
any triangulation of extremal cone generators in $R_A$ with the origin will use
at most polynomially many top-dimensional simplices. Hence the number of top-dimensional simplices of any
partial triangulation $\Tr$ will be polynomially bounded since it is a subset
of the final triangulation.

The {\bf for} statement on line \ref{alg:ept:l12} repeats for every simplex $B$ and
every $(|B|-1)$-simplex of $B$. As before, the number of simplices $B$ is
polynomially bounded, and there are at most $n$ $(|B|-1)$-simplices of $B$.
Thus the {\bf for} statement will repeat a polynomial number of times.

Finally, by Lemma \ref{lem:vislp}, determining if $C$ is visible to $\ve v_i$ can
be done in polynomial time. Therefore Algorithm \ref{alg:ept} runs in a
polynomial time.
\end{proof}

\begin{corollary} \label{elemconepoly}
%%% Let $M$ be a matroid on $n$ elements and $\Co \subseteq \R^n$ a tangent cone of $\Po(M)$ generated by $l$ extremal rays $\{\ve r_1,\ldots,\ve r_l\}$. 
Let $r$ be a fixed integer, $n$ be an integer, $A \subseteq [n]$ with $|A|
\leq r$, and let $\Co \subseteq \R^n$ be a cone generated by extremal rays
$\{\ve r_1,\ldots,\ve r_l\} \subseteq R_A$.  A triangulation of $\Co$ can be
computed in polynomial time in the input of the extremal ray generators $\{\ve
r_1,\ldots,\ve r_l \}$. 
\end{corollary}

\begin{proof}
%%%   Since $\Co$ is a tangent cone of $\Po(M)$ then $\{\ve r_1,\ldots,\ve r_l\} \subseteq R_A$ for some $A \subseteq [n]$ by Lemma \ref{lem:adj}
Let $\Po_\Co := \conv \{ \ve 0, \ve r_1,\ldots,\ve r_t\} $. We give an algorithm which produces a triangulation of
    $\Po_\Co := \conv \{ \ve 0, \ve r_1,\ldots,\ve r_t\} $ such that each full-dimensional
    simplex has $\ve 0$ as a vertex. Such a triangulation would extend to a triangulation of
    the cone $\Co$. This can be accomplished by applying two placing
    triangulations: one to triangulate the boundary of $\Po_\Co$ not incident to
    $\ve 0$, and another to attach the triangulated boundary faces to $\ve 0$.  The
    algorithm goes as follows:

\begin{itemize}
    \item[1)] {\bf Triangulate $\Po_\Co$ using the placing triangulation algorithm. Call it $\Tr'$.} 

    \item[2)] {\bf Triangulate $\Po_\Co$ using the boundary faces of $\Tr'$ which do not contain $v$}. 
\vskip 10pt
\begin{minipage}{\linewidth}
\begin{algorithm}[Triangulation joining $\ve 0$ to boundary faces]
    \mbox{}
    \vskip .25cm
    \label{alg:jt}
    \begin{center}
        \begin{tabular}{l}
            %\toprule
            \begin{minipage}{.85\linewidth}
                \begin{algorithmic}
                        %%%   \INPUT A triangulation $\Tr'$ of $\Po_\Co$.
                        \item[Input:] A triangulation $\Tr'$ of $\Po_\Co$, given by its vertices.
                        %%%   \OUTPUT A triangulation $\Tr$ of $\Po_\Co$ such that every highest dimension simplex of $\Tr$ is incident to $0$.
                        \item[Output:] A triangulation $\Tr$ of $\Po_\Co$ such that every highest dimension simplex of $\Tr$ is incident to $\ve 0$.
                        \STATE $\Tr := \emptyset$ \label{alg:jt:1}
                        %\STATE $\mathscr{S} := \emptyset$ \label{alg:jt:2}
                        \FOR{each $C$ where $C$ is a $(|B|-1)$-simplex of $B \in \Tr'$} \label{alg:jt:3}
                            \IF{$C$ is not a $(|A|-1)$-simplex of $A \in \Tr'$ where $A \neq B$} \label{alg:jt:4}
                                \STATE $\Tr := \Tr \cup \{ \, C \cup \{\ve 0\} \, \}$ \label{alg:jt:5}
                            \ENDIF
                        \ENDFOR
                        \RETURN $\Tr$

                \end{algorithmic}
            \end{minipage}\\
            %\bottomrule
        \end{tabular}
    \end{center}
\end{algorithm}
\end{minipage}
\vskip 10pt
\end{itemize}

\begin{figure}[!htb]    
    \centering
    \hfil(a)\subfigure{\scalebox{0.6}{\ifpdf
    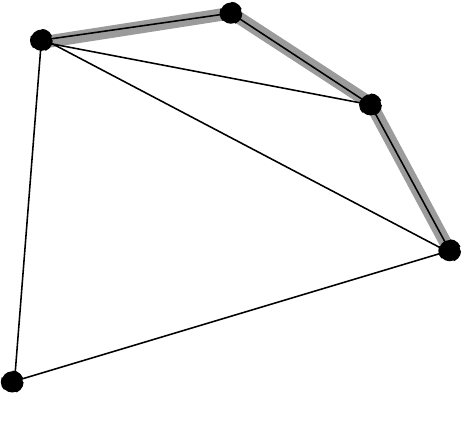
    \else
    \input{bndtri1.pstex_t}
    \fi}}
    \hfil(b)\subfigure{\scalebox{0.6}{\ifpdf
    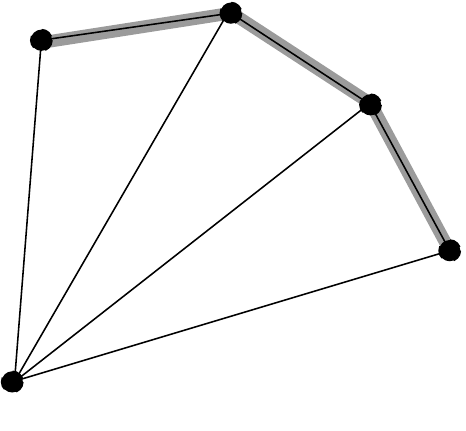
    \else
    \input{bndtri2.pstex_t}
    \fi}}
    \caption{A triangulation $\Tr'$ of $\Po_\Co$ can be used to extend to a triangulation $\Tr$ such that $\ve 0$ is incident to every highest dimensional simplex}
\end{figure}

By Theorem \ref{thm:ept}, triangulating $\Po_\Co$ using Algorithm \ref{alg:ept}
can be done in polynomial time. Algorithm \ref{alg:jt} indeed produces a triangulation of $\Po_\Co$. It covers
$\Po_\Co$ since every extremal ray generator $\ve r_k$ of $\Co$ is on some
$(\dim(\Co)-1)$-simplex. Moreover, $\Tr$ by construction has the property that
the intersection of any two simplices of $\Tr$ is a simplex. Step
\ref{alg:jt:4} checks if $C$ is on the boundary, since if $C$ is on the
boundary it will not be on the intersection of two higher-dimensional
simplices. 

Step \ref{alg:jt:3} repeats a polynomial number of times since any
triangulation of $\Po_\Co$ has at most a polynomial number of simplices, and
each simplex $B$ has at most $n$ $(|B|-1)$-simplices. Step \ref{alg:jt:4} can
be computed in polynomial time since again there are only polynomially many
simplices $B$ in the triangulation $\Tr'$ and at most $n$ $(|B|-1)$-simplices
to check if they are equal to $C$. Hence, Algorithm \ref{alg:jt} runs in
polynomial time.
\end{proof}

\begin{example}
    The tangent cone at the vertex $\ve e_B := \ve e_{\{1,2,3\}}$ on the polytope $\Po(M(K_4))$ can be triangulated as:
\begin{align*}
\Big\{ \, & \left\{ \ve e_{\{2, 3, 5\}}- \ve e_B,\, \ve e_{\{2, 3, 4\}}- \ve e_B,\, \ve e_{\{1, 3, 6\}}- \ve e_B,\, \ve e_{\{1, 3, 4\}}- \ve e_B,\, \ve e_{\{1, 2, 6\}} - \ve e_B \right\},  \\
& \left\{ \ve e_{\{2, 3, 5\}}- \ve e_B,\, \ve e_{\{1, 3, 6\}}- \ve e_B,\, \ve e_{\{1, 3, 4\}}- \ve e_B,\, \ve e_{\{1, 2, 6\}}- \ve e_B,\, \ve e_{\{1, 2, 5\}}- \ve e_B \right\}, \\
& \left\{ \ve e_{\{2, 3, 5\}}- \ve e_B,\, \ve e_{\{2, 3, 4\}}- \ve e_B,\, \ve e_{\{1, 3, 4\}}- \ve e_B,\, \ve e_{\{1, 2, 6\}}- \ve e_B,\, \ve e_{\{1, 2, 5\}}- \ve e_B \right\} \, \Big\}.
\end{align*}
\end{example}

        %%%%%%%%%%%%%%%%%%%%%%%%%%%%%%%%%%%%%%%%%%%%%%%%%
        %%%%%%            POLYMATROIDS             %%%%%%
        %%%%%%%%%%%%%%%%%%%%%%%%%%%%%%%%%%%%%%%%%%%%%%%%%
\section{Polymatroids}
\label{subsec:polymatroids}
We will show that certain lemmas from Subsection
    \ref{subsec:ontangent} also hold for certain polymatroids.
 Recall that the rank of the matroid $M$ is the size of any basis of $M$ which
 equals $\varphi([n])$. Our lemmas from Subsection \ref{subsec:ontangent} rely
 on the fact that $M$ has fixed rank, that is, for some $r \in \Z$, $r \geq 0$,
 $\varphi(A) \leq r$ for all $A \subseteq [n]$. We will show that a similar
 condition on a polymatroid rank function is sufficient for the lemmas of
 Subsection \ref{subsec:ontangent} to hold.  

 \begin{lemma} \label{lem:polyvert}
    Let $\psi \colon 2^{[n]} \longrightarrow \mathbb{N}$ be an integral polymatroid
    rank function where $\psi(A) \leq r$ for all $A \subseteq [n]$, where $r$ is
    a fixed integer. Then the number of vertices of $\Po(\psi)$ is bounded by a
    polynomial in $n$ of degree $r$. 
 \end{lemma}
\begin{proof}
  It is known that if $\psi$ is integral then all vertices of
    $\Po(\psi)$ are integral \cite{Welsh1976Matroid-Theory}. The number of vertices
    of $\Po(\psi)$ can be bounded by the number of non-negative integral solutions to $x_1 +
    \cdots +  x_n \leq r$, which has ${n + r \choose r}$
    solutions, a polynomial in $n$ of degree $r$ \cite{Stanley1997Enumerative-Com}.
\end{proof}

\begin{lemma} \label{lem:polyenum}
    Let $\psi \colon 2^{[n]} \longrightarrow \mathbb{N}$ be an integral polymatroid
    rank function. If $\ve v$ is a vertex of $\Po(\psi)$ then all adjacent
    vertices of $\ve v$ can be enumerated in polynomial time. Moreover if
    $\psi(A) \leq r$ for all $A \subseteq [n]$, where $r$ is a fixed integer,
    then the vertices of $\Po(\psi)$ can be enumerated in polynomial time.
\end{lemma}

\begin{proof}
     If $\ve v$ is a vertex of $\Po(\psi)$
     then generating and listing all adjacent vertices to $\ve v$ can
     be done in polynomial time by \autoref{lem:adj}. If $\psi(A) \leq r$ for all $A
     \subseteq [n]$, where $r$ is a fixed integer, then, by Lemma
     \ref{lem:polyvert}, there is a polynomial number of vertices for
     $\Po(\psi)$.  We know that $\ve 0 \in \R^n$ is a vertex of any
     polymatroid. Therefore, beginning with $\ve 0$, we can perform a
     breadth-first search, which is output-sensitive polynomial time, on the graph of
     $\Po(\psi)$, enumerating all vertices of $\Po(\psi)$.
\end{proof}

What remains to be shown is that these polymatroids have cones like the
ones in Subsection \ref{subsec:ontangent}.  

%We first recall some needed
%definitions from \cite{Topkis1984Adjacency-on-Po}.  Let $\ve v,\ve w \in
%\R^n$ and define $\Delta(\ve v,\ve w) := \{ \, i \in [n] \mid v_i \neq w_i
%\, \}$ and $\cl(\ve v) := \{ \, S \mid S \subseteq [n], \ \sum_{i \in S}
%v_i = \psi(S) \, \}$. Let $F = \{ f_1,\ldots,f_{|F|} \}$ be an ordered
%subset of $[n]$ and $F_i := \{f_1,\ldots,f_i\}$. If $\psi$ is a polymatroid
%rank function then we construct $\ve v \in \R^n$ where $v_i = \psi(F_i) -
%\psi(F_{i-1})$ where $ v_j = 0$ when $j \notin F$ and one says $F$
%\emph{generates} $\ve v$. A classical result of Edmonds \cite{Edmonds2003Submodular-func} says that the set of vectors generated by all ordered subsets of $[n]$
%is exactly the set of vertices of $\Po(\psi)$. Now we can restate an important
%lemma. 
%
%\begin{lemma}[See Theorem 4.1 and Section 2 in \cite{Topkis1984Adjacency-on-Po}.] \label{lem:topadj}
%    Let $\psi$ be a polymatroid rank function. If $\ve v$ and $\ve w$ are vertices
%    of the polymatroid $\Po(\psi)$ then either
%    \begin{itemize}
%        \item[(i)] $|\Delta(\ve v, \ve w)| = 1$ or
%        \item[(ii)] $\cl(\ve v) = \cl(\ve w)$ and $\Delta(\ve v,\ve w) =
%        \{c,d\}$ for some $c,d \in [n]$ where there exists some ordered set $F
%        = \{f_1,\ldots,f_{|F|}\}$ which generates $\ve v$ with $f_{k+1} = d$
%        and $f_k = c$ for some integer $k$, $1 \leq k \leq |F| -1$; moreover
%        the ordered set\\ $\tilde F :=
%        \{f_1,\ldots,f_{k-1},f_{k+1},f_k,f_{k+2},\ldots,f_{|F|}\}$ generates
%        $\ve w$.
%    \end{itemize}
%\end{lemma}

\begin{lemma} \label{lem:polyelem}
    Let $\psi$ be an integral polymatroid rank function and $\Co$ the tangent
    cone of a vertex $\ve v$ of the polymatroid $\Po(\psi)$, translated to the
    origin. Then $\Co$ is generated by extremal ray generators $\{\ve
    r_1,\ldots, \ve r_l\} \subseteq R_{\supp(\ve v)}$, where $R_{\supp(\ve v)}$
    is an elementary set of $\supp(\ve v)$.
\end{lemma}
\begin{proof}
Let $\psi \colon 2^{[n]} \longrightarrow \Z$ be a integral polymatroid rank function.
    Let $\ve v$ and $\ve w$ be adjacent vertices of the polymatroid
    $\Po(\psi)$. Using Lemma \ref{lem:topadj}, if $|\Delta(\ve v, \ve w)| = 1$
    then $\ve w - \ve v = h \ve e_i$ where $h$ is some integer and $\ve e_i$ is the
    standard $i$th elementary vector for some $i \in [n]$. If $h < 0$ then certainly $i
    \in \supp(\ve v)$, else $i \in [n]$. Thus $\ve w - \ve v$, a generator of
    $\Co$, is parallel to a vector in $R_{\supp(\ve v)}$.

Let $\ve v$ and $\ve w$ be adjacent and satisfy ({\bf ii}) in Lemma
\ref{lem:topadj}, where $\Delta(\ve v,\ve w) = \{c,d\}$. Hence there exists an
$F = \{f_1,\ldots,f_{|F|}\}$ which generates $\ve v$ with $f_{k+1} = d$ and
$f_k = c$ for some integer $k$, $1 \leq k \leq |F| -1$; moreover the ordered
set $\tilde F := \{f_1,\ldots,f_{k-1},f_{k+1},f_k,f_{k+2},\ldots,f_{|F|}\}$
generates $\ve w$. First we note that $\psi(F_{k-1}) = \psi(\tilde F_{k-1})$
and $\psi(F_{k+1}) = \psi(\tilde F_{k+1})$. By assumption, we know $ v_{c}
\neq  w_{c}$, $v_{d} \neq w_{d}$ and $v_l = w_l$ for all $l
    \in [n] \setminus \{c,d\}$. Thus
\begin{align*}
    (\ve v - \ve w)_{c} & =  v_{c} -  w_{c} & = \phantom{-} & \psi(F_{k+1}) - \psi(F_k) - \left(  \psi(\tilde F_k) - \psi(\tilde F_{k-1}) \right) \\
                        &                         & = \phantom{-} & \psi(F_{k+1}) - \psi(F_k) -  \psi(\tilde F_k) + \psi(\tilde F_{k-1})  \\
    \intertext{and}
    (\ve v - \ve w)_{d} & =  v_{d} -  w_{d} & = \phantom{-} & \psi(F_{k}) - \psi(F_{k-1}) - \left(  \psi(\tilde F_{k+1}) - \psi(\tilde F_{k}) \right) \\
                        &                         & = \phantom{-} & \psi(F_{k}) - \psi(\tilde F_{k-1}) - \left(  \psi(F_{k+1}) - \psi(\tilde F_{k}) \right) \\
                        &                         & = - & \psi(F_{k+1}) + \psi(F_{k}) +  \psi(\tilde F_{k}) - \psi(\tilde F_{k-1}).    
\end{align*}
Therefore $(\ve v - \ve w)_{c} = - (\ve v - \ve w)_{d}$ and $\ve w - \ve v$ is
    parallel to $\ve e_d - \ve e_c$. Moreover, $c \in \supp(\ve v)$ since $\ve w, \ve v \geq \ve 0$
    by assumption that $\ve v, \ve w \in \Po(\psi)$. Thus $\ve w - \ve v$, a
    generator of $\Co$, is parallel to a vector in $R_{\supp(\ve v)}$.
\end{proof}

\section{The construction of a short multivariate rational generating function}
\label{subsec:short}

From the knowledge of triangulations of tangent cones of matroid
polytopes, independence matroid polytopes, and polymatroids we will now recover short
multivariate generating functions. 
\begin{remark}
  Notice that
  formula~\eqref{eq:ratgenfun-identity-triang} is of exponential size, even
  when the triangulation~$\Tr$ only has polynomially many simplicial cones
  of maximal dimension.
  The reason is that, when the dimension~$n$ is
  allowed to vary, there are exponentially many intersecting proper faces in
  the set~$\hat\Tr$.  Therefore, we cannot
  use~\eqref{eq:ratgenfun-identity-triang} to compute the multivariate
  rational generating function of~$\Co$ in polynomial time for varying dimension.
\end{remark}
%%% \mkoeppesays{Need to explain why Brion's dualization trick does NOT work in this situation.  Can we prove that the volume of the dual of the vertex cones is exponential in our case, so at least our volume-bounding method would provably fail in the dual situation.  Note that the uniform matroids do not form a suitable family, since their the polyopes associated with their normal cones have a volume bounded by $n/n!$.}

\begin{figure}[t]
  \centering
  \ifpdf
    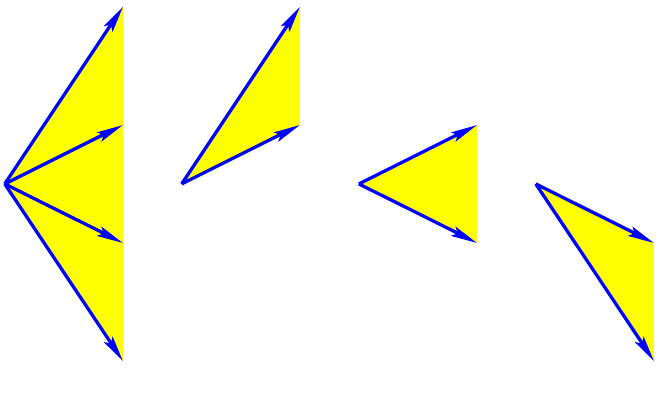
    \else
    \input{2d-triang.pstex_t}
    \fi
  \caption{An identity, valid modulo lower-dimensional cones, corresponding to
    a polyhedral subdivision of a cone}
  \label{fig:2d-triang}
\end{figure}

\begin{figure}[t]
  \ifpdf
    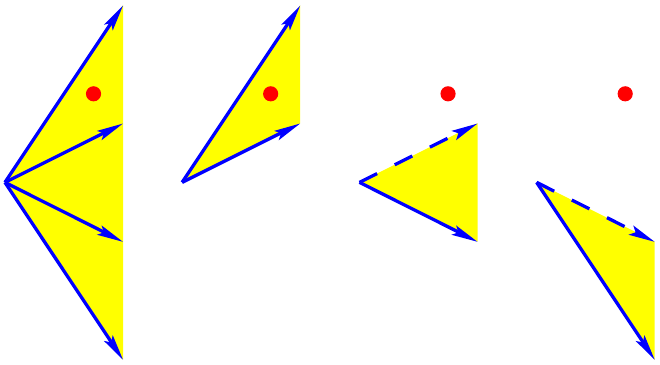
    \else
    \input{2d-triang-halfopen-good.pstex_t}
    \fi
  \caption{The technique of half-open exact decomposition.  The relative
    location of the vector $\ve y$ (represented by a dot) determines which defining
    inequalities are strict (broken lines) and which are weak (solid lines).}
  \label{fig:2d-triang-halfopen-good}
\end{figure}

To obtain a shorter formula, we use the technique of \emph{half-open exact
  decompositions}~\cite{koeppe-verdoolaege:parametric}, which is a
refinement of the method of ``irrational'' perturbations
\cite{beck-sottile:irrational,koeppe:irrational-barvinok}.  We use the
following result; see also \autoref{fig:2d-triang} and
\autoref{fig:2d-triang-halfopen-good}.  

\begin{lemma}~
  \label{th:exactify-identities}
  \begin{enumerate}[\rm(a)]
  \item\label{th:exactify-identities-general-part} Let
    \begin{equation}
      \label{eq:full-source-identity}
      \sum_{i\in I_1} \epsilon_i [\Co_i] + \sum_{i\in I_2} \epsilon_i [\Co_i] = 0
    \end{equation}
    be a linear identity (with rational coefficients~$\epsilon_i$) of
    indicator functions of cones~$\Co_i\subseteq\R^n$, where the cones~$\Co_i$
    are full-dimensional for $i\in I_1$ and lower-dimensional for $i\in I_2$.
    Let each cone be given as
    \begin{align}
      \Co_i &= \bigl\{\, \ve x \in\R^n : \langle \ve b^*_{i,j}, \ve x\rangle
      \leq 0 \text{ for $j\in J_i$}\,\bigr\}.
    \end{align}
    Let $\ve y\in\R^n$ be a vector such that $\langle \ve b^*_{i,j}, \ve
    y\rangle \neq 0$ for all $i\in I_1\cup I_2$, $j\in J_i$.  For $i\in I_1$,
    we define the ``half-open cone''
    \begin{equation}
      \label{eq:half-open-by-y}
      \begin{aligned}
        \tilde \Co_i = \Bigl\{\, \ve x\in\R^d : {}& \langle \ve b^*_{i,j}, \ve
        x\rangle \leq 0
        \text{ for $j\in J_i$ with $\langle \ve b^*_{i,j}, \ve y \rangle < 0$,} \\
        & \langle \ve b^*_{i,j}, \ve x\rangle < 0 \text{ for $j\in J_i$ with
          $\langle \ve b^*_{i,j}, \ve y \rangle > 0$} \,\Bigr\}.
      \end{aligned}
    \end{equation}
    Then
    \begin{equation}
      \label{eq:target-identity}
      \sum_{i\in I_1} \epsilon_i [\tilde \Co_i] = 0.
    \end{equation}
\item\label{th:exactify-identities-triang-part}
  In particular, let 
  \begin{equation}
    \label{eq:source-identity-triang}
    [\Co] = \sum_{\Co_i\in \Tr} [\Co_i] + \sum_{\Co_i\in \hat \Tr} \epsilon_i\, [\Co_i]
  \end{equation}
  be the identity corresponding to a triangulation of the cone~$\Co$, where
  $\Tr$ is the set of simplicial cones of maximal dimension and $\hat \Tr$ is
  the set of intersecting proper faces.
  Then there exists a polynomial-time algorithm to construct a vector~$\ve
  y\in\Q^n$ such that the above construction yields the identity
  \begin{equation}
    \label{eq:half-open-triang}
    [\Co] = \sum_{\Co_i\in \Tr} [\tilde \Co_i],
  \end{equation}
  which describes a partition of $\Co$ into half-open cones of maximal dimension.
  \end{enumerate}
\end{lemma}
\begin{proof}
  Part~(\ref{th:exactify-identities-general-part}) is a slightly less general
  form of Theorem~3 in \cite{koeppe-verdoolaege:parametric}.
  Part~(\ref{th:exactify-identities-triang-part}) follows from the discussion
  in section~2 of  \cite{koeppe-verdoolaege:parametric}.
\end{proof}


%% To apply this result to our situation, we construct a suitable vector~$\ve
%% y\in\R^d$ which allows to decide which of the inequalities of the cones should
%% be strict.  
%% \begin{lemma}
%%   Let $n\in N$ and $A\subseteq[n]$.  
%%   Let $\Co$ be a cone generated by rays $\ve
%%   r_k \in R_A$. %%% \{ \, \ve e_i - \ve e_j, \, \ve e_i, \, -\ve e_j \; | \; i \in [n] \setminus A, \; j \in A \, \}$.
%%   Let 
%%   \begin{displaymath}
%%     \Co = \{\, \ve x : \langle \ve b^*_{k}, \ve x\rangle \leq 0 \text{
%%       for $k\in J$}\,\bigr\}.
%%   \end{displaymath}
%%   Then $\ve y_A := $ \mkoeppesays{Need to find a vector that works -- what do we
%%     know about the facet vectors~$\ve b^*_{k}$?}
%%   satisfies $\langle \ve b^*_{i,j}, \ve y\rangle \neq 0$ for $k\in J$.
%% \end{lemma}

Since the cones in a triangulation~$\Tr$ of all tangent cones~$\Co_\Po(\ve v)$
of our polytopes are unimodular by \autoref{lem:eluni}, we can efficiently
write the multivariate generating functions of their half-open counterparts.
\begin{lemma}[Lemma~9 in \cite{koeppe-verdoolaege:parametric}] \label{th:unimodular-formula}
  Let $\tilde \Co\subseteq\R^n$ be an $N$-dimensional half-open pointed
  simplicial affine cone with an integral apex~$\ve v\in\Z^n$ and the
  ray description 
  \begin{align}
    \label{eq:explicit}
    \tilde \Co&= \Bigl\{\, \ve v + \textstyle\sum_{j=1}^N \lambda_j \ve b_j :{} 
    \lambda_j \geq 0 \text{ for $j\in J_\leq$ and} \ 
    %% \\ &
    \lambda_j > 0 \text{ for $j\in J_<$} 
    \,\Bigr\}
  \end{align}
  where $J_\leq\cup J_< = \{1,\dots,N\}$ and $\ve b_j\in\Z^n\setminus\{\ve0\}$.
  We further assume that $\tilde \Co$ is unimodular, i.e., the vectors~$\ve
  b_j$ form a basis of the lattice $(\R\ve b_1+\dots+\R\ve b_N)\cap\Z^n$.
  Then the unique point in the fundamental parallelepiped of the half-open
  cone~$\tilde \Co$ is
  \begin{equation}
    \ve a = \ve v + \sum_{j\in J_<} \ve b_j,
  \end{equation}
  and the generating function of $\Co$ is given by
  \begin{equation}
    g_\Co(\ve{z}) = \frac{\ve{z}^{\ve a}}{\prod_{j=1}^N (1 - \ve{z}^{\ve b_j})}.
  \end{equation}
\end{lemma}
\smallbreak

Taking all results together, we obtain:
\begin{corollary}
  \label{th:matroid-multiratgenfun-summary}
  Let $r$ be a fixed integer.  There exist algorithms that, given 
  \begin{enumerate}[(a)]
  \item a matroid $M$ on $n$ elements, presented by an
    evaluation oracle for its rank function~$\varphi$, which is bounded above
    by~$r$, or
  \item an evaluation oracle for an integral polymatroid rank function $\psi \colon 2^{[n]}
    \longrightarrow \mathbb{N}$, which is bounded above by~$r$,
  \end{enumerate}
  compute in time polynomial in~$n$ vectors $\ve a_i\in \Z^n$, $\ve
  b_{i,j}\in\Z^n\setminus\{\ve0\}$, and $\ve v_i\in\Z^n$ for $i\in I$ (a
  polynomial-size index set) and $j=1,\dots,N$, where $N\leq n$, such that
  the multivariate generating function of $\Po(M)$, $\Po^\In(M)$ and
  $\Po(\psi)$, respectively, is the sum of rational functions
  \begin{align}
    g_{\Po}(\ve z) &= \sum_{i\in I} \frac{\ve z^{\ve a_i}} {\prod_{j=1}^N
      (1-\ve z^{\ve b_{i,j}})}
  \end{align}
  and  the $k$-th dilation of the polytope has the multivariate rational
  generating function 
  \begin{align}
    g_{k\Po}(\ve z) &= \sum_{i\in I} \frac{\ve z^{\ve a_i+(k-1)\ve v_i}}
    {\prod_{j=1}^N (1-\ve z^{\ve b_{i,j}})}.
  \end{align}
\end{corollary}
\begin{proof}
    Lemma \ref{lem:bl} implies that finding the multivariate generating
    function of $\Po(M)$, $\Po^\In(M)$ or $\Po(\psi)$ can be reduced to
    finding the multivariate generating functions of their tangent cones.
    Moreover, $\Po(M)$, $\Po^\In(M)$ and $\Po(\psi)$ have only polynomially in
    $n$ many vertices as described in subsection \ref{subsec:ontangent} or
    Lemma \ref{lem:polyvert}. Enumerating their vertices can be done in
    polynomial time by Lemma \ref{lem:adj} or \ref{lem:polyenum}. 

    Given a vertex $\ve v$ of $\Po(M)$, $\Po^\In(M)$ or $\Po(\psi)$, its neighbors
    can be computed in polynomial time by Lemma \ref{lem:adj} or
    \ref{lem:polyenum}. The tangent cone $\Co(\ve v)$ at $\ve v$ is generated by
    elements in $R_A$, where $R_A$ is an elementary set for some $A \subseteq [n]$.
    See subsection \ref{subsec:ontangent} or Lemma \ref{lem:polyelem}. We also proved
    in Lemma \ref{lem:eluni} that every triangulation of $\Co(\ve v)$ generated by elements in
    $R_A$ is unimodular and Lemma \ref{lem:elsimpbnd} states that any
    triangulation of $\Co(\ve v)$ has at most a polynomial in $n$ number of top-dimensional
    simplices. Moreover, a triangulation of the cone $\Co(\ve v)$ can be computed
    in polynomial time by Lemma \ref{thm:ept}. Finally, using Lemmas
    \ref{th:exactify-identities} and \ref{th:unimodular-formula} we can write
    the polynomial sized multivariate generating function of $\Co(\ve v)$ in polynomial time.
    Therefore we can write the multivariate generating function of $\Po(M)$,
    $\Po^\In(M)$ or $\Po(\psi)$ in polynomial time.
\end{proof}

\section{Polynomial-time specialization of rational generating functions in
  varying dimension}
\label{subsec:specialize}

We now compute the Ehrhart polynomial $i(\Po, k) = \#(k\Po\cap\Z^n)$
from the multivariate rational generating function $g_{k\Po}(\ve z)$
of \autoref{th:matroid-multiratgenfun-summary}.  This amounts to the
problem of evaluating or \emph{specializing} a rational generating
function $g_{k\Po}(\ve z)$, depending on a parameter~$k$, at the
point~$\ve z=\ve 1$. This is a pole of each of its summands but a
regular point (removable singularity) of the function itself.  From
now on we call this the \emph{specialization problem}. We explain a
very general procedure to solve it which we hope will allow future
applications.

To this end,
let the generating function of a polytope~$\Po\subseteq\R^n$ be given in the form
\begin{equation}
  \label{eq:generating-function-0}
  g_\Po(\ve z) = \sum_{i\in I} \epsilon_i \frac{\ve
    z^{\ve a_i}}{\prod_{j=1}^{s_i} (1-\ve z^{\ve b_{ij}})}
\end{equation}
where $\epsilon_i\in\{\pm1\}$, $\ve a_i\in\Z^n$, and $\ve
b_{ij}\in\Z^n\setminus\{\ve0\}$.  Let $s = \max_{i\in I} s_i$ be the maximum
number of binomials in the denominators.  In general, if $s$ is allowed to
grow, more poles need to be considered for each summand, so the evaluation
will need more computational effort.

In previous literature, the specialization problem has been considered, but not in
sufficient generality for our purpose.
In the original paper by~\citet[Lemma~4.3]{bar}, the dimension~$n$ is
fixed, and each summand has exactly $s_i=n$ binomials in the denominator.  The
same restriction can be found in the survey by~\citet{barvinok:99}.  In the more
general algorithmic theory of monomial substitutions developed by
\citet{barvinok-woods-2003,Woods:thesis}, there is no assumption on the
dimension~$n$, but the number~$s$ of binomials in the denominators is fixed.  The
same restriction appears in the paper by~\citet[Lemma
2.15]{verdoolaege-woods-2005}.  In a recent paper, \citet[section~5]{barvinok-2006-ehrhart-quasipolynomial} 
gives a polynomial-time algorithm for the specialization problem for rational
functions of the form 
\begin{equation}
  \label{eq:generating-function-with-exponents}
  g(\ve z) = \sum_{i\in I} \epsilon_i \frac{\ve
    z^{\ve a_i}}{\prod_{j=1}^{s} {(1-\ve z^{\ve b_{ij}})}^{\gamma_{ij}}}
\end{equation}
where the dimension~$n$ is fixed, the number~$s$ of different binomials in
each denominator equals~$n$, but 
the multiplicity $\gamma_{ij}$ is varying.  

We will show that the technique from \citet[section~5]{barvinok-2006-ehrhart-quasipolynomial} can be implemented
in a way such that we obtain a polynomial-time algorithm even for the case of
a general formula~\eqref{eq:generating-function-0},
when the dimension and the number of binomials are allowed to grow.

\begin{theorem}[Polynomial-time specialization \cite{deloera-haws-koeppe:ehrhart-matroid}]
  \label{th:specialization-polytime-varydim}
  \begin{enumerate}[\rm(a)]
  \item There exists an algorithm for computing the
    specialization of a rational function of the form
    \begin{equation}
      \label{eq:generating-function-01}
      g_\Po(\ve z) = \sum_{i\in I} \epsilon_i \frac{\ve
        z^{\ve a_i}}{\prod_{j=1}^{s_i} (1-\ve z^{\ve b_{ij}})}
    \end{equation}
    at its removable singularity~$\ve z=\ve 1$, 
    which runs in time polynomial in the encoding size of its data
    $\epsilon_i\in\Q$, $\ve a_i\in\Z^n$ for $i\in I$ and 
    $\ve b_{ij}\in\Z^n$ for $i\in I$, $j=1,\dots,s_i$, 
    even when the dimension~$n$ and
    the numbers~$s_i$ of terms in the denominators are not fixed.
  \item In particular, there exists a polynomial-time algorithm that, given 
    data $\epsilon_i\in\Q$, $\ve a_i\in\Z^n$ for $i\in I$ and 
    $\ve b_{ij}\in\Z^n$ for $i\in I$, $j=1,\dots,s_i$ describing a rational
    function in the form~\eqref{eq:generating-function-01},  
    computes 
    a vector $\velambda\in\Q^n$ with $\inner{\velambda, \ve b_{ij}} \neq
    0$ for all $i,j$ 
    and rational weights $w_{i,l}$
    for $i\in I$ and $l=0,\dots,s_i$. 
    Then the number of integer points is given by
    \begin{equation}
      \#(\Po\cap\Z^n) 
      = \sum_{i\in I} \epsilon_i \sum_{l=0}^{s_i} w_{i,l} \inner{\velambda, \ve a_i}^l.
    \end{equation}
  \item Likewise, given a parametric rational function for the dilations of an
    integral polytope~$\Po$, 
    \begin{align}\label{eq:generating-function-dil}
      g_{k\Po}(\ve z) 
      &= \sum_{i\in I} \epsilon_i \frac{\ve z^{\ve a_i+(k-1)\ve v_i}} 
      {\prod_{j=1}^d (1-\ve z^{\ve b_{i,j}})},
    \end{align}
    the Ehrhart polynomial $i(\Po, k) = \#(k\Po\cap\Z^n)$ is given by the explicit formula
    \begin{equation}
      i(\Po, k) = 
      \sum_{m=0}^M \paren{
        \sum_{i\in I} 
        \epsilon_i 
        \inner{\velambda,\ve v_i}^m 
        \sum_{l=m}^{s_i} \binom{l}{m} w_{i,l} 
        \inner{\velambda,\ve a_i - \ve v_i}^{l-m} 
      } k^m,
    \end{equation}
    where $M = \min\{s,\dim\Po\}$.
  \end{enumerate}
\end{theorem}

\begin{proof} [Proof of Theorem \ref{volume}] 
Corollary \ref{th:matroid-multiratgenfun-summary} and
Theorem \ref{th:specialization-polytime-varydim} imply Theorem
(\autoref{volume}) directly. 
\end{proof}

The remainder of this section contains the proof of
\autoref{th:specialization-polytime-varydim}.  We follow \cite{barvinok:99} and
recall the definition of Todd polynomials. We will prove that they can be
efficiently evaluated in rational arithmetic.
\begin{definition}\label{def:todd-poly}
  We consider the function
  \begin{displaymath}
    H(x, \xi_1,\dots,\xi_s) = \prod^s_{i=1} \frac{x\xi_i}{1-\exp\{-x\xi_i\}},
  \end{displaymath}
  a function that is analytic in a neighborhood of~$\ve 0$.
  The $m$-th ($s$-variate) \emph{Todd polynomial} is the coefficient of~$x^m$
  in the Taylor expansion 
  \begin{displaymath}
    H(x, \xi_1,\dots,\xi_s) = \sum_{m=0}^\infty \td_m(\xi_1,\dots,\xi_s) x^m.
  \end{displaymath}
\end{definition}
We remark that, when the numbers $s$ and $m$ are allowed to vary, the Todd
polynomials have an exponential number of monomials.
\begin{theorem}[\cite{deloera-haws-koeppe:ehrhart-matroid}]\label{th:todd-evaluation}
  The Todd polynomial $\td_m(\xi_1,\dots,\xi_s)$ can be evaluated for given
  rational data $\xi_1,\dots,\xi_s$ in time polynomial in~$s$, $m$, and the
  encoding length of $\xi_1,\dots,\xi_s$.
\end{theorem}
The proof makes use of the following lemma.
\begin{lemma}[\cite{deloera-haws-koeppe:ehrhart-matroid}]\label{lemma:toddy-expansion} 
  The function $h(x) = x/(1-\mathrm{e}^{-x})$ is a function that is analytic in a
  neighborhood of~$0$.  Its Taylor series about $x=0$ is of the form 
  \begin{equation}
    \label{eq:toddy-coefficients-format}
    h(x) = \sum_{n=0}^\infty b_n x^n \quad\text{where}\quad
    b_n = \frac1{n!\, (n+1)!} c_n
  \end{equation}
  with integer coefficients $c_n$ that have a binary encoding length of
  $\mathrm O(n^2 \log n)$.
  The coefficients $c_n$ can be computed from the recursion
  \begin{equation}
    \label{eq:recursion-toddy}
    \begin{aligned}
      c_0 &= 1\\
      c_n &= \sum_{j=1}^n (-1)^{j+1} \binom{n+1}{j+1} \frac{n!}{(n-j+1)!} c_{n-j}
      && \text{for $n=1,2,\dots$.}
    \end{aligned}
  \end{equation}
\end{lemma}
\begin{proof}
  The reciprocal function $h^{-1}(x) = (1-\mathrm{e}^{-x})/x$ has the Taylor
  series
  \begin{displaymath}
    h^{-1}(x) = \sum_{i=0}^\infty a_n x^n
    \quad\text{with}\quad
    a_n = \frac{(-1)^n}{(n+1)!}.
  \end{displaymath}
  Using the identity $h^{-1}(x)h(x) = \bigl( \sum_{n=0}^\infty a_n x^n \bigr)
  \bigl( \sum_{n=0}^\infty b_n x^n \bigr) = 1$, we obtain the recursion
  \begin{equation}
    \begin{aligned}
      b_0 &= \tfrac1{a_0} = 1\\
      b_n &= -(a_1 b_{n-1} + a_2 b_{n-2} + \dots + a_n b_0)
      && \text{for $n=1,2,\dots$.}
    \end{aligned}
  \end{equation}
  We prove~\eqref{eq:toddy-coefficients-format} inductively.  Clearly $b_0 =
  c_0 = 1$.  For $n=1,2,\dots$, we have
  \begin{align*}
    c_n &= n! \, (n + 1)! \, b_n \\
    &= -n! \, (n + 1)! \, (a_1 b_{n-1} + a_2 b_{n-2} + \dots + a_n b_0) \\
    &= n! \, (n + 1)! \, \sum_{j=1}^n \frac{(-1)^{j+1}}{(j+1)!} \cdot
    \frac{1}{(n-j)!\, (n-j+1)!} c_{n-j} \\
    &= \sum_{j=1}^n (-1)^{j+1} \frac{(n+1)!}{(j+1)!\, (n-j)!} \cdot 
    \frac{n!}{(n-j+1)!} c_{n-j}.
  \end{align*}
  Thus we obtain the recursion formula~\eqref{eq:recursion-toddy}, which also
  shows that all $c_n$ are integers.  A rough estimate shows that 
  \begin{displaymath}
    \abs{c_n} \leq n (n+1)!\, n!\, \abs{c_{n-1}} \leq \paren[big]{(n+1)!}^2 \abs{c_{n-1}},
  \end{displaymath}
  thus $\abs{c_n}\leq \paren[big]{(n+1)!}^{2n}$, so $c_n$ has a binary
  encoding length of $\mathrm O(n^2 \log n)$. 
\end{proof}

\begin{proof}[Proof of \autoref{th:todd-evaluation}]
  By definition, we have
  \begin{displaymath}
    H(x,\xi_1,\dots,\xi_s) = \sum_{m=0}^\infty \td_m(\xi_1,\dots,\xi_s) x^m 
    = \prod_{j=1}^s h(x\xi_j).
  \end{displaymath}
  From \autoref{lemma:toddy-expansion} we have
  \begin{equation}
    \label{eq:toddy-expansion-trunc}
    h(x\xi_j) = \sum_{n=0}^m \beta_{j,n} x^n  + o(x^m)
    \quad\text{where}\quad
    \beta_{j,n} = \frac{\xi_j^n}{n!\,(n+1)!} c_n 
  \end{equation}
  with integers~$c_n$ given by the recursion~\eqref{eq:recursion-toddy}.
  Thus we can evaluate $\td_m(\xi_1,\dots,\xi_s)$  by summing over all the
  possible compositions $n_1 + \dots + n_{s} = m$ of the order~$m$
  from the orders~$n_j$ of the factors:
  \begin{equation}\label{eq:evaluation-by-compositions}
    \td_m(\xi_1,\dots,\xi_s)
    = \sum_{\substack{(n_1,\dots,n_{s})\in\Z_+^{s}\\
        n_1 + \dots + n_{s} = m}} \!\!\!\! \beta_{1,n_1} \dots
    \beta_{s,n_s}
  \end{equation}
  We remark that the length of the above sum is equal to the number of compositions 
  of~$m$ into $s$~non-negative parts, 
  \begin{align*}
    C'_s(m) %% Note this C is not a cone, but notation for #compositions. --mkoeppe
    &= \binom{m + s - 1}{s - 1}\\
    &= \frac{(m+s-1)(m+s-2)\dots(m+s-(s-1))} {(s-1) (s-2) \dots 2\cdot 1}\\
    &= \Omega\paren{\paren[big]{1 + \tfrac{m}{s-1}}^{s}},
  \end{align*}
  which is \emph{exponential} in~$s$ (whenever $m\geq s$).  Thus we cannot evaluate the
  formula~\eqref{eq:evaluation-by-compositions}
  efficiently when $s$ is allowed to grow.  

  However, we show that we can evaluate $\td_m(\xi_1,\dots,\xi_s)$ more efficiently. 
  To this end, we multiply up the $s$ truncated Taylor
  series~\eqref{eq:toddy-expansion-trunc}, one
  factor at a time, truncating after order~$m$.  Let us denote
  \begin{align*}
    H_1(x) &= h(x\xi_1) \\
    H_2(x) &= H_1(x) \cdot h(x\xi_2) \\
    &\vdots \\
    H_s(x) &= H_{s-1}(x) \cdot h(x\xi_s) = H(x,\xi_1,\dots,\xi_s).
  \end{align*}
  Each multiplication
  can be implemented in $\mathrm O(m^2)$ elementary rational operations.  
  We finally show that all numbers appearing in the calculations have polynomial
  encoding size.  Let $\Xi$ be the largest binary encoding size of any of the rational
  numbers~$\xi_1,\dots,\xi_s$.  Then every $\beta_{j,n}$ given
  by~\eqref{eq:toddy-expansion-trunc} has a binary encoding 
  size $\mathrm O(\Xi n^5 \log^3 n)$.  Let $H_j(x)$ have the truncated Taylor
  series $\sum_{n=0}^m \alpha_{j,n} x^n + o(x^m)$ and let $A_j$ denote the
  largest binary encoding length of any $\alpha_{j,n}$ for $n\leq m$.  Then 
  \begin{displaymath}
    H_{j+1}(x) = \sum_{n=0}^m \alpha_{j+1,n} x^n + o(x^m)
    \quad\text{with}\quad
    \alpha_{j+1,n} =  \sum_{l=0}^n \alpha_{j,l} \beta_{j, n-l} .
  \end{displaymath}
  Thus the binary encoding size of $\alpha_{j+1,n}$ (for $n\leq m$) is bounded
  by $A_j + \mathrm O(\Xi m^5 \log^3 m)$. 
  Thus, after $s$ multiplication steps, the encoding size of the coefficients
  is bounded by $\mathrm O(s \Xi m^5 \log^3 m)$, a polynomial quantity.
\end{proof}

\begin{proof}[Proof of \autoref{th:specialization-polytime-varydim}]
\emph{Parts (a) and (b).}
  We recall the technique of \citet[Lemma~4.3]{bar}, refined by
  \citet[section~5]{barvinok-2006-ehrhart-quasipolynomial}.
  
  We first construct a rational vector $\velambda\in\Z^n$ such that
  $\inner{\velambda, \ve b_{ij}}\neq0$ for all $i,j$.  One such construction
  is to consider the \emph{moment curve} $\velambda(\xi) = (1, \xi, \xi^2,
  \dots, \xi^{n-1})\in\R^n$.  For each exponent vector $\ve b_{ij}$ occuring
  in a denominator of~\eqref{eq:generating-function-0}, the function
  $f_{ij}\colon \xi\mapsto \inner{\velambda(\xi), \ve b_{ij}}$ is a
  polynomial function of degree at most $n-1$.  Since $\ve b_{ij}\neq\ve 0$,
  the function $f_{ij}$ is not identically zero.  Hence $f_{ij}$ has at most
  $n-1$ zeros.  By evaluating all functions $f_{ij}$ for $i\in I$ and
  $j=1,\dots,s_i$ at $M = (n-1)s|I|+1$ different values for~$\xi$, for
  instance at the integers $\xi=0,\dots, M$, we can find one~$\xi=\bar\xi$
  that is not a zero of any~$f_{ij}$.  Clearly this search can be implemented
  in polynomial time, even when the dimension~$n$ and the number~$s$ of terms
  in the denominators are not fixed.  We set $\velambda =
  \velambda(\bar\xi)$.
  
  For $\tau>0$, let us consider the points $\ve z_\tau = \ve e^{\tau
    \velambda} = (\exp\{ \tau\lambda_1 \}, \dots, \exp\{ \tau\lambda_n
  \})$. We have 
  \begin{displaymath}
    \ve z_\tau^{\ve b_{ij}} 
    = \prod_{l=1}^n \exp\{ \tau\lambda_l b_{ijl} \} = \exp\{ \tau
    \inner{\velambda, \ve b_{ij}} \};
  \end{displaymath}
  since $\inner{\velambda, \ve b_{ij}}\neq0$ for all $i,j$, all the
  denominators $1-\ve z_\tau^{\ve b_{ij}}$ are nonzero.  
  Hence for every
  $\tau>0$, the point $\ve z_\tau$ is a regular point not only of $g(\ve z)$
  but also of the individual summands of~\eqref{eq:generating-function-0}. 
  We have
  \begin{align*}
    g(\ve 1) 
    &= \lim_{\tau\to0^+} \sum_{i\in I} \epsilon_i \frac{\ve
      z_\tau^{\ve a_i}}{\prod_{j=1}^{s_i} (1-\ve z_\tau^{\ve b_{ij}})} 
\displaybreak[0]\\
    &= \lim_{\tau\to0^+} \sum_{i\in I} \epsilon_i 
    \frac{\exp\{ \tau \inner{\velambda, \ve a_{i}} \}}
    {\prod_{j=1}^{s_i} (1-\exp\{ \tau \inner{\velambda, \ve b_{ij}} \} ) } 
\displaybreak[0]\\
    &= \lim_{\tau\to0^+} \sum_{i\in I} \epsilon_i \,
    \tau^{-s_i}
    \exp\{ \tau \inner{\velambda, \ve a_{i}} \}
    \prod_{j=1}^{s_i}
    \frac{\tau} 
    {1-\exp\{ \tau \inner{\velambda, \ve b_{ij}} \} } 
\displaybreak[0]\\
    &= \lim_{\tau\to0^+} \sum_{i\in I} \epsilon_i \,
    \tau^{-s_i}
    \exp\{ \tau \inner{\velambda, \ve a_{i}} \}
    \prod_{j=1}^{s_i}
    \frac{-1}{\inner{\velambda, \ve b_{ij}}} h(-\tau\inner{\velambda, \ve
      b_{ij}}) 
\displaybreak[0]\\
    &= \lim_{\tau\to0^+} \sum_{i\in I} \epsilon_i\,  \frac{(-1)^{s_i}}{\prod_{j=1}^{s_i}
      \inner{\velambda, \ve b_{ij}}}  
    \,    \tau^{-s_i}
    \exp\{ \tau \inner{\velambda, \ve a_{i}} \}
    H(\tau, -\inner{\velambda, \ve b_{i1}}, \dots, -\inner{\velambda, \ve b_{is_i}})
  \end{align*}
  where $H(x, \xi_1,\dots,\xi_{s_i})$ is the function
  from \autoref{def:todd-poly}. 
  We will compute the limit 
  by finding the constant term of the Laurent expansion of each summand
  about~$\tau=0$. 
  Now the function $\tau\mapsto \exp\{ \tau \inner{\velambda, \ve a_{i}} \}$
  is holomorphic and has the Taylor series
  \begin{equation}
    \exp\{ \tau \inner{\velambda, \ve a_{i}} \}
    = \sum_{l = 0}^{s_i} \alpha_{i,l} \tau^l + o(\tau^{s_i})
    \quad\text{where}\quad
    \alpha_{i,l} = \frac{\inner{\velambda, \ve a_{i}}^l}{l!}, 
    \label{eq:numerator-series}
  \end{equation}
  and $H(\tau, \xi_1,\dots,\xi_{s_i})$ has the Taylor series
  \begin{displaymath}
    H(\tau, \xi_1,\dots,\xi_s) = \sum_{m=0}^{s_i} \td_m(\xi_1,\dots,\xi_s)
    \tau^m + o(\tau^{s_i}). 
  \end{displaymath}
  Because of
  the factor $\tau^{-s_i}$, which gives rise to a pole of order~$s_i$ in the
  summand, we can compute the constant term of the
  Laurent expansion by summing over all the possible compositions $s_i = l +
  (s_i - l)$ of the order~$s_i$:
  \begin{equation}
    g(\ve 1) 
    = \sum_{i\in I} \epsilon_i 
    \frac{(-1)^{s_i}}{\prod_{j=1}^{s_i} \inner{\velambda, \ve b_{ij}}}
    \sum_{l=0}^{s_i} \frac{\inner{\velambda, \ve a_i}^l}{l!} 
    \td_{s_i-l}(-\inner{\velambda, \ve b_{i1}}, \dots, -\inner{\velambda, \ve b_{is_i}}).
  \end{equation}
  We use the notation
  \begin{displaymath}
    w_{i,l}  =  (-1)^{s_i} \frac{\td_{s_i-l}( -\langle \velambda,
      \ve b_{i,1} \rangle, \dots, -\langle \velambda, \ve b_{i,s_i}\rangle)}
    {l!\cdot \langle \velambda,
      \ve b_{i,1} \rangle \cdots \langle \velambda, \ve b_{i,s_i}\rangle}
    \quad\text{for $i\in I$ and $l=0,\dots,s_i$};
  \end{displaymath}
  these rational numbers can be computed in polynomial time using
  \autoref{th:todd-evaluation}.  
  We now obtain the formula of the claim, 
  \begin{displaymath}
    g(\ve 1) 
    = \sum_{i\in I} \epsilon_i \sum_{l=0}^{s_i} w_{i,l} \inner{\velambda, \ve a_i}^l.
  \end{displaymath}
  \smallbreak

\noindent\emph{Part~(c).}
Applying the same technique to the parametric rational function~\eqref{eq:generating-function-dil}, we obtain
\begin{align*}
  \#(k\Po\cap\Z^n) &= g_{k\Po}(\ve1)\\
  &= \sum_{i\in I} \epsilon_i \sum_{l=0}^{s_i} w_{i,l} \inner{\velambda, \ve a_i+(k-1)\ve v_i}^l, 
\displaybreak[0]\\
  &= \sum_{i\in I} \epsilon_i \sum_{l=0}^{s_i} w_{i,l} \sum_{m=0}^l
  \binom{l}{m} \inner{\velambda,\ve a_i - \ve v_i}^{l-m} k^m \inner{\velambda,\ve
    v_i}^m 
\displaybreak[0]\\
  &= \sum_{m=0}^s \paren{
    \sum_{i\in I} 
    \epsilon_i 
    \inner{\velambda,\ve v_i}^m 
    \sum_{l=m}^{s_i} \binom{l}{m} w_{i,l} 
    \inner{\velambda,\ve a_i - \ve v_i}^{l-m} 
  } k^m,
\end{align*}
an explicit formula for the Ehrhart polynomial.  We remark
that, since the Ehrhart polynomial is of degree equal to the
dimension of~$\Po$, all coefficients of $k^m$ for $m > \dim\Po$ must vanish.
Thus we obtain the formula of the claim, where we sum only up to $\min\{s,
\dim\Po\}$ instead of~$s$. 
\end{proof}

%\section*{Acknowledgments}
%The first author was supported by NSF grant DMS-0608785. The second
%author was supported by VIGRE-GRANT DMS-0135345 and DMS-0636297. The
%third author was supported by a 2006/2007 Feodor Lynen Research
%Fellowship from the Alexander von Humboldt Foundation. We thank Dillon
%Mayhew and Gordon Royle for providing data for the matroids listed in
%Table \ref{tab:hvec}.

%\bibliographystyle{abbrvnat}
%\bibliography{barvinok,matroid}

%\end{document}

   \chapter[% 
     Algebraic Combinatorics of Matroid Polytopes
   ]{% 
     Results on the Algebraic Combinatorics of Matroid Polytopes 
   }%
   \label{ch:3rdChapterLabel}
   \pdfoutput=1

        %%%%%%%%%%%%%%%%%%%%%%%%%%%%%%%%%%%%%%%%%%%%%%%%%
        %%%%%%    PROPERTIES OF THE H^*-VECTORS    %%%%%%
        %%%%%%%%%%%%%%%%%%%%%%%%%%%%%%%%%%%%%%%%%%%%%%%%%

\section[$h^*$-vectors and Ehrhart Polynomials]{Algebraic Properties of $h^*$-vectors and Ehrhart polynomials of Matroid Polytopes} \label{sec:hstar}
%In the second part of the paper, developed in Section \ref{sec:hstar},
Here we investigate algebraic properties of the Ehrhart functions of
matroid polytopes: The \emph{Ehrhart series} of a polytope $\Po$ is
the infinite series $\sum_{k=0}^\infty i(\Po,k)t^k$. We recall the
following classic result about Ehrhart series (see e.g.,
\cite{Hibi1992Algebraic-Combi, Stanley1996Combinatorics-a}). Let 
$\Po \subseteq \R^n$ be an integral convex polytope of dimension
$d$. Then it is known that its Ehrhart series is a rational function
of the form
\begin{equation} \label{lem:ratehr}
\sum_{k=0}^\infty i(\Po,k)t^k = \frac{ h^*_0 + h^*_1 t + \cdots + h^*_{d-1}t^{d-1} + h^*_d t^d }{(1-t)^{d+1}}.
\end{equation}
The numerator is
often called the $\mathit{h}^*$\emph{-polynomial} of $\Po$ (some authors also call it the \emph{Ehrhart h-polynomial}), and we
define the coefficients of the polynomial in the numerator of Lemma
\ref{lem:ratehr}, $h^*_0 + h^*_1 t + \cdots + h^*_{d-1}t^{d-1} + h^*_d
t^d$, as the $h^*$-\emph{vector} of $\Po$, which we write as $\ve
h^*(\Po) := (h^*_0, h^*_1, \dots, h^*_{d-1}, h^*_d)$.

A vector $(c_0,\ldots,c_d)$ is \emph{unimodal} if there exists an index~$p$, $0
\leq p \leq d$, such that $c_{i-1} \leq c_{i}$ for $i \leq p$ and $c_{j}
\geq c_{j+1}$ for $j \geq p$. Due to its algebraic implications,
several authors have studied the unimodality of $h^*$-vectors (see
\cite{Hibi1992Algebraic-Combi} and \cite{Stanley1996Combinatorics-a}
and references therein). 

Suppose, as before, that $\Po \subseteq
\R^n$, and each vertex of $\Po$ has integral (or rational)
vertices. Let $Y_1,Y_2,\dots,Y_n$ and $T$ be indeterminates over a
field $K$. Letting $q \geq 1$ we define $A(\Po)_q$ as the vector space
over $K$ which is spanned by the monomials
$Y_1^{\alpha_1},Y_2^{\alpha_2},\dots,Y_n^{\alpha_n}T^q$ such that
$(\alpha_1,\alpha_2,\dots,\alpha_n) \in q\Po \cap \Z^n$.  Since $\Po$
is convex it follows that $A(\Po)_qA(\Po)_p \subseteq A_{q+p}(\Po)$
for all $p,q$, and thus the \emph{Ehrhart ring} of $\Po$, $A(\Po) =
\bigoplus_{q=0}^\infty A(\Po)_q$, is a graded algebra
\cite{Hibi1992Algebraic-Combi,Stanley1996Combinatorics-a}.

It is well-known that if the Ehrhart ring of an integral polytope
$\Po$, $A(\Po)$, is Gorenstein, then $\ve h^*(\Po)$ is unimodal, and
symmetric \cite{Hibi1992Algebraic-Combi,Stanley1996Combinatorics-a}.
Nevertheless, the vector $\ve h^*(\Po)$ can be unimodal even when the
Ehrhart ring $A(\Po)$ is not Gorenstein. For matroid polytopes, their
Ehrhart ring is indeed often not Gorenstein.  For instance, De Negri and Hibi
\cite{Negri1997Gorenstein-Alge} prove explicitly when the
Ehrhart ring of a uniform matroid polytope is Gorenstein or not.  Two
fascinating facts, uncovered through experimentation, are that all
$h^*$-vectors seen thus far are unimodal, even for the cases when
their Ehrhart rings are not Gorenstein. In addition, when we computed
the explicit Ehrhart polynomials of matroid polytopes we observe their
coefficients are always positive. We conjecture:

\begin{conjecture} \label{unimodalconj} Let $\Po(M)$ be the matroid polytope of a matroid $M$.
    \begin{itemize}
        \item[(A)]  The $h^*$-vector of $\Po(M)$ is unimodal.
        \item[(B)]  The coefficients of the Ehrhart polynomial of $\Po(M)$ are positive. 
    \end{itemize}
\end{conjecture}

We have proved both parts of this conjecture in many instances.  
A class of matroids that we considered are the uniform matroids; 
recall that the \emph{uniform matroid} on $n$ elements of rank $r$ is the
collection of all $r$-subsets of $n$. 
Using computers, we were able to verify \autoref{unimodalconj} for all uniform
matroids up to $75$ elements as 
well as for a wide variety of non-uniform matroids which are collected at
\cite{HawsMatroid-Polytop}.  We include here this information just
for the $28$ famous matroids presented in
\cite{Oxley1992Matroid-Theory}. Results in \cite{Katzman:math0408038},
with some additional careful calculations, imply that 
\autoref{unimodalconj} is true for all rank $2$ uniform matroids.
Regarding part (A) of the conjecture we were also able to prove
\emph{partial unimodality} for uniform matroids of rank $3$. Concretely we obtained:

\begin{theorem} \label{partialuni} \mbox{}
\begin{itemize}
\item[(1)] Conjecture \ref{unimodalconj} is true for all uniform
matroids up to $75$ elements and all uniform matroids of rank 2.
It is also true for all matroids listed in \cite{HawsMatroid-Polytop}.

\item[(2)] Let $\Po(U^{3,n})$ be the matroid polytope of a uniform
            matroid of rank $3$ on $n$ elements, and let $I$ be a non-negative
            integer.  Then there exists $n(I) \in \mathbb{N}$ such that for all
            $n \geq n(I)$ the $h^*$-vector of $\Po(U^{3,n})$, $\left(
                    h^*_0,\ldots, h^*_{n} \right)$, is non-decreasing
            from index 0 to $I$. That is, $ h^*_{0} \leq h^*_{1} \leq \dots
            \leq h^*_{I}$. 
\end{itemize}
\end{theorem}

Using the programs {\tt cdd+} \cite{cdd}, {\tt LattE} \cite{lattesoft} and {\tt LattE macchiato} \cite{latte-macchiato} we explored patterns for the
Ehrhart polynomials of matroid polytopes. Since previous authors
proposed other invariants of a matroid (e.g., the Tutte polynomials
and the invariants of \cite{SpeyerA-matroid-invar, billera-2006}) we wished to know how well does the Ehrhart
polynomial distinguish non-isomorphic matroids. It is natural to
compare it with other known invariants. Some straightforward
properties are immediately evident.  For example, the Ehrhart
polynomial of a matroid and that of its dual are equal. Also the
Ehrhart polynomial of a direct sum of matroids is the product of their
Ehrhart polynomials.

\let\maybemidrule=\relax
\begin{table}[htbp]
\caption{Coefficients of the Ehrhart polynomials of selected matroids in \cite{billera-2006, SpeyerA-matroid-invar}}
 \label{tab:belspeyer}
\begin{centering}
\def\arraystretch{1.2}
\begin{tabular}{ll}
\toprule
& Ehrhart Polynomial \\
\midrule
%%% Note that Speyer1 and 2 are transposed from their files.
Speyer1 & $1   , \frac{21}{5}   , \frac{343}{45}   , \frac{63}{8}   , \frac{91}{18}   , \frac{77}{40}   , \frac{29}{90}  $ \\
\maybemidrule
Speyer2 & $1   , \frac{135}{28}   , \frac{3691}{360}   , \frac{1511}{120}   , \frac{88}{9}   , \frac{39}{8}   , \frac{529}{360}   , \frac{89}{420}  $\\
\maybemidrule
BJR1 & $1  , \frac{109}{30}   , \frac{23}{4}   , \frac{59}{12}   , \frac{9}{4}   , \frac{9}{20}   $\\
\maybemidrule
BJR2 & $ 1  , \frac{211}{60}   , \frac{125}{24}   , \frac{33}{8}   , \frac{43}{24}   , \frac{43}{120}   $\\
\maybemidrule
BJR3 &  $1   , \frac{83}{20}   , \frac{2783}{360}   , \frac{199}{24}   , \frac{391}{72}   , \frac{247}{120}   , \frac{61}{180}  $\\
\maybemidrule
BJR4 & $1   , \frac{25}{6}   , \frac{353}{45}   , \frac{101}{12}   , \frac{193}{36}   , \frac{23}{12}   , \frac{53}{180}  $\\
\bottomrule
\end{tabular}
\end{centering}
\end{table}

We call the last two matroids in Figure $2$ in \cite{SpeyerA-matroid-invar} {\it Speyer1} and {\it Speyer2} and the  matroids of Figure $2$ in \cite{billera-2006} {\it BJR1}, {\it BJR2}, {\it BJR3}, and {\it BJR4}, and list their Ehrhart polynomials in Table \ref{tab:belspeyer}.
We note that {\it BJR3} and {\it BJR4} 
have the same Tutte polynomial, yet their Ehrhart polynomials are
different. This proves that the Ehrhart polynomial cannot be computed
using deletion and contractions, as is the case for the Tutte
polynomial.  Examples {\it BJR1} and {\it BJR2} show that the Ehrhart
polynomial of a matroid may help to distinguish
non-isomorphic matroids: These two matroids are not isomorphic yet
they have the same Tutte polynomials and the same quasi-symmetric
function studied in \cite{billera-2006}. Although they share some
properties, there does not seem to be an obvious relation to Speyer's
univariate polynomials introduced in \cite{SpeyerA-matroid-invar};
examples {\it Speyer1} and {\it Speyer2} show they are relatively
prime with their corresponding Ehrhart polynomials.

Our experiments included, among others, many examples coming from
small graphical matroids, random realizable matroids over fields of
small positive characteristic, and the classical examples listed in
the Appendix of \cite{Oxley1992Matroid-Theory} for which we list the
results in Table \ref{tab:hvec}. For a comprehensive list of all our
calculations visit \cite{HawsMatroid-Polytop}.  Soon it became evident
that all instances verified both parts of our \autoref{unimodalconj}.

\begin{table}[htbp]
\caption{Coefficients of the Ehrhart polynomials and $h^*$-vectors of selected
  matroids in \cite{Oxley1992Matroid-Theory}} 
 \label{tab:hvec}
\begin{centering}
\tiny
\def\arraystretch{1.6}
\begin{tabular}{lccll}
\toprule
& \normalsize$n$ & \normalsize$r$ & \normalsize$h^*$-vector & \normalsize Ehrhart Polynomial \\
\midrule
 %%%   $U^{2,4}$ & $4$ & $2$ & $1,2,1$                             & $  1   , \frac{7}{3}   , 2   , \frac{2}{3}   $ \\
 %%%   \maybemidrule                                                       
 %%%   $U^{2,5}$ & $5$ & $2$ & $5,5,1$                             & $  1   , \frac{35}{12}   , \frac{85}{24}   , \frac{25}{12}   , \frac{11}{24}   $ \\
 %%%   \maybemidrule                                                       
 %%%   $U^{3,5}$ & $5$ & $3$ & $5,5,1$                             & $  1   , \frac{35}{12}   , \frac{85}{24}   , \frac{25}{12}   , \frac{11}{24}   $ \\
 %%%   \maybemidrule                                                       
 $K_4$ & $6$ & $3$ & $1, 10, 20, 10, 1$                            & $  1   , \frac{107}{30}   , \frac{21}{4}   , \frac{49}{12}   , \frac{7}{4}   , \frac{7}{20}   $ \\
 \maybemidrule                                                             
 $W^3$ & $6$ & $3$ & $1, 11, 24, 11, 10$                            & $  1   , \frac{18}{5}   , \frac{11}{2}   , \frac{9}{2}   , 2   , \frac{2}{5}   $ \\
 \maybemidrule                                                             
 $Q_6$ & $6$ & $3$ & $1, 12, 28, 12, 1$                            & $  1   , \frac{109}{30}   , \frac{23}{4}   , \frac{59}{12}   , \frac{9}{4}   , \frac{9}{20}   $ \\
 \maybemidrule                                                             
 $P_6$ & $6$ & $3$ & $1, 13, 32, 13, 1$                            & $  1   , \frac{11}{3}   , 6   , \frac{16}{3}   , \frac{5}{2}   , \frac{1}{2}   $ \\
 \maybemidrule                                                             
 %%%   $U^{3,6}$ & $6$ & $3$ & $1, 14, 36, 14, 1$                  & $  1   , \frac{37}{10}   , \frac{25}{4}   , \frac{23}{4}   , \frac{11}{4}   , \frac{11}{20}   $ \\
 %%%   \maybemidrule                                                       
 $R_6$ & $6$ & $3$ & $1, 12, 28, 12, 1$                            & $  1   , \frac{109}{30}   , \frac{23}{4}   , \frac{59}{12}   , \frac{9}{4}   , \frac{9}{20}   $ \\
 \maybemidrule                                                             
 $F_7$ & $7$ & $3$ & $21, 98, 91, 21, 1$                           & $  1   , \frac{21}{5}   , \frac{343}{45}   , \frac{63}{8}   , \frac{91}{18}   , \frac{77}{40}   , \frac{29}{90}   $ \\
 \maybemidrule                                                             
 %%%   $F_7^$ & $7$ & $4$ & $21, 98, 91, 21, 1$                    & $  1   , \frac{21}{5}   , \frac{343}{45}   , \frac{63}{8}   , \frac{91}{18}   , \frac{77}{40}   , \frac{29}{90}   $ \\
 %%%   \maybemidrule                                                       
 $F_7^-$ & $7$ & $3$ & $21, 101, 97, 22, 1$                        & $  1   , \frac{253}{60}   , \frac{2809}{360}   , \frac{33}{4}   , \frac{193}{36}   , \frac{61}{30}   , \frac{121}{360}   $ \\
 \maybemidrule                                                             
 %%%   $F_7^-^$ & $7$ & $4$ & $21, 101, 97, 22, 1$               & $  1   , \frac{253}{60}   , \frac{2809}{360}   , \frac{33}{4}   , \frac{193}{36}   , \frac{61}{30}   , \frac{121}{360}   $ \\
 %%%   \maybemidrule                                                       
 $P^7$ & $7$ & $3$ & $21, 104, 103, 23, 1$                         & $  1   , \frac{127}{30}   , \frac{479}{60}   , \frac{69}{8}   , \frac{17}{3}   , \frac{257}{120}   , \frac{7}{20}   $ \\
 \maybemidrule                                                             
 %%%   $P^7^$ & $7$ & $4$ & $21, 104, 103, 23, 1$                & $  1   , \frac{127}{30}   , \frac{479}{60}   , \frac{69}{8}   , \frac{17}{3}   , \frac{257}{120}   , \frac{7}{20}   $ \\
 %%%   \maybemidrule                                                       
 $AG(3,2)$ & $8$ & $4$ & $1, 62, 561, 1014, 449, 48, 1$            & $  1   , \frac{209}{42}   , \frac{1981}{180}   , \frac{881}{60}   , \frac{119}{9}   , \frac{95}{12}   , \frac{499}{180}   , \frac{89}{210}   $ \\
 \maybemidrule                                                             
 $AG'(3,2)$ & $8$ & $4$ & $1, 62, 562, 1023, 458, 49, 1$           & $  1   , \frac{299}{60}   , \frac{4007}{360}   , \frac{5401}{360}   , \frac{122}{9}   , \frac{2911}{360}   , \frac{1013}{360}   , \frac{77}{180}   $ \\
 \maybemidrule                                                             
 $R_8$ & $8$ & $4$ & $1, 62, 563, 1032, 467, 50, 1$                & $  1   , \frac{524}{105}   , \frac{1013}{90}   , \frac{1379}{90}   , \frac{125}{9}   , \frac{743}{90}   , \frac{257}{90}   , \frac{136}{315}   $ \\
 \maybemidrule                                                             
 $F_8$ & $8$ & $4$ & $1, 62, 563, 1032, 467, 50, 1$                & $  1   , \frac{524}{105}   , \frac{1013}{90}   , \frac{1379}{90}   , \frac{125}{9}   , \frac{743}{90}   , \frac{257}{90}   , \frac{136}{315}   $ \\
 \maybemidrule                                                             
 $Q_8$ & $8$ & $4$ & $1, 62, 564, 1041, 476, 51, 1$                & $  1   , \frac{2099}{420}   , \frac{4097}{360}   , \frac{1877}{120}   , \frac{128}{9}   , \frac{337}{40}   , \frac{1043}{360}   , \frac{61}{140}   $ \\
 \maybemidrule                                                             
 $S_8$ & $8$ & $4$ & $1, 44, 337, 612, 305, 40, 1$                 & $  1   , \frac{1021}{210}   , \frac{377}{36}   , \frac{475}{36}   , \frac{193}{18}   , \frac{511}{90}   , \frac{65}{36}   , \frac{67}{252}   $ \\
 \maybemidrule                                                             
 $V_8$ & $8$ & $4$ & $1, 62, 570, 1095, 530, 57, 1$                & $  1   , \frac{2117}{420}   , \frac{4367}{360}   , \frac{2107}{120}   , \frac{146}{9}   , \frac{1133}{120}   , \frac{1133}{360}   , \frac{193}{420}   $ \\
 \maybemidrule                                                             
 $T_8$ & $8$ & $4$ & $1, 62, 564, 1041, 476, 51, 1$                & $  1   , \frac{2099}{420}   , \frac{4097}{360}   , \frac{1877}{120}   , \frac{128}{9}   , \frac{337}{40}   , \frac{1043}{360}   , \frac{61}{140}   $ \\
 \maybemidrule                                                             
 $V_8^+$ & $8$ & $4$ & $1, 62, 569, 1086, 521, 56, 1$              & $  1   , \frac{151}{30}   , \frac{2161}{180}   , \frac{3103}{180}   , \frac{143}{9}   , \frac{1669}{180}   , \frac{559}{180}   , \frac{41}{90}   $ \\
 \maybemidrule                                                             
 $L_8$ & $8$ & $4$ & $1, 62, 567, 1068, 503, 54, 1$                & $  1   , \frac{527}{105}   , \frac{529}{45}   , \frac{83}{5}   , \frac{137}{9}   , \frac{134}{15}   , \frac{136}{45}   , \frac{47}{105}   $ \\
 \maybemidrule                                                             
 $J$ & $8$ & $4$ & $1, 44, 339, 630, 323, 42, 1$                   & $  1   , \frac{512}{105}   , \frac{193}{18}   , \frac{83}{6}   , \frac{205}{18}   , \frac{361}{60}   , \frac{17}{9}   , \frac{23}{84}   $ \\
 \maybemidrule                                                             
 $P_8$ & $8$ & $4$ & $1, 62, 565, 1050, 485, 52, 1$                & $  1   , \frac{1051}{210}   , \frac{2071}{180}   , \frac{2873}{180}   , \frac{131}{9}   , \frac{1547}{180}   , \frac{529}{180}   , \frac{277}{630}   $ \\
 \maybemidrule                                                             
 $W_4$ & $8$ & $4$ & $1, 38, 262, 475, 254, 37, 1$                 & $  1   , \frac{135}{28}   , \frac{3691}{360}   , \frac{1511}{120}   , \frac{88}{9}   , \frac{39}{8}   , \frac{529}{360}   , \frac{89}{420}   $ \\
 \maybemidrule                                                             
 $W^4$ & $8$ & $4$ & $1, 38, 263, 484, 263, 38, 1$                 & $  1   , \frac{169}{35}   , \frac{467}{45}   , \frac{581}{45}   , \frac{91}{9}   , \frac{227}{45}   , \frac{68}{45}   , \frac{68}{315}   $ \\
 \maybemidrule                                                             
 %%%   $_{3,3}^$ & $9$ & $4$ & $78, 1116, 3492, 3237, 927, 72, 1$ & $  1   , \frac{307}{56}   , \frac{137141}{10080}   , \frac{3223}{160}   , \frac{37807}{1920}   , \frac{211}{16}   , \frac{5743}{960}   , \frac{1889}{1120}   , \frac{8923}{40320}   $ \\
 %%%   \maybemidrule                                                       
 $K_{3,3}$ & $9$ & $5$ & $78, 1116, 3492, 3237, 927, 72, 1$        & $  1   , \frac{307}{56}   , \frac{137141}{10080}   , \frac{3223}{160}   , \frac{37807}{1920}   , \frac{211}{16}   , \frac{5743}{960}   , \frac{1889}{1120}   , \frac{8923}{40320}   $ \\
 \maybemidrule                                                             
 $AG(2,3)$ & $9$ & $3$ & $1, 147, 1230, 1885, 714, 63, 1$          & $  1   , \frac{1453}{280}   , \frac{41749}{3360}   , \frac{581}{32}   , \frac{34069}{1920}   , \frac{927}{80}   , \frac{4541}{960}   , \frac{239}{224}   , \frac{449}{4480}   $ \\
 \maybemidrule                                                             
 Pappus & $9$ & $3$ & $1, 147, 1230, 1915, 744, 66, 1$             & $  1   , \frac{729}{140}   , \frac{3573}{280}   , \frac{381}{20}   , \frac{1499}{80}   , \frac{243}{20}   , \frac{49}{10}   , \frac{153}{140}   , \frac{57}{560}   $ \\
 \maybemidrule                                                             
 Non-Pappus & $9$ & $3$ & $1, 147, 1230, 1925, 754, 67, 1$         & $  1   , \frac{4379}{840}   , \frac{25951}{2016}   , \frac{9287}{480}   , \frac{21967}{1152}   , \frac{987}{80}   , \frac{2855}{576}   , \frac{3701}{3360}   , \frac{275}{2688}   $ \\
 \maybemidrule                                                             
 $Q_3(\mathit{GF}(3)^*)$ & $9$ & $3$ & $1, 147, 1098, 1638, 632, 59, 1$      & $  1   , \frac{433}{84}   , \frac{3079}{252}   , \frac{4193}{240}   , \frac{5947}{360}   , \frac{167}{16}   , \frac{601}{144}   , \frac{787}{840}   , \frac{149}{1680}   $ \\
 \maybemidrule                                                             
 $R_9$ & $9$ & $3$ & $1, 147, 1142, 1717, 656, 60, 1$              & $  1   , \frac{723}{140}   , \frac{49}{4}   , \frac{88}{5}   , \frac{24217}{1440}   , \frac{1291}{120}   , \frac{625}{144}   , \frac{821}{840}   , \frac{133}{1440}   $ \\
\bottomrule
\end{tabular}
\end{centering}
\end{table}

By far the most comprehensive study we made was for the family of
uniform matroids. In this case we based our computations on the theory
of Veronese algebras as developed by M. Katzman in
\cite{Katzman:math0408038}. There, Katzman gives an explicit
equation for the $h^*$-vector of uniform matroid polytopes (again,
using the language of Veronese algebras). For this family we were
able to verify computationally the conjecture is true for all uniform
matroids up to $75$ elements and to prove partial unimodality as explained
in the introduction.

\begin{lemma} \label{lem:conjb}
The coefficients of the Ehrhart polynomial of the matroid polytope of the uniform matroid $U^{2,n}$ are positive.
\end{lemma}
\begin{proof}
We begin with the expression in Corollary 2.2 in \cite{Katzman:math0408038} which explicitly gives the Ehrhart polynomial of $\Po(U^{r,n})$ as
\begin{equation} \label{eq:uniehrhart}
 i(\Po(U^{r,n}), k) = \sum_{s=0}^{r-1} (-1)^s {n \choose s}{k(r-s) -s +n -1 \choose n-1}.
\end{equation}
Letting $r=2$, Equation \eqref{eq:uniehrhart} becomes
\begin{equation*}
\frac{(2k+n-1)(2k+n-2) \cdots (2k+1)}{(n-1)!} - n \frac{(k+n-2)(k+n-3) \cdots (k)}{(n-1)!}
\end{equation*}
We next consider the coefficient of $k^{n-p-1}$ for $0 \leq p \leq n -1$, which can be written as
%%%   \begin{multline*}
%%%   \sum_{\substack{1 \leq j_1 < \cdots < j_p \leq n-1 \\ j_q \in \Z }} \frac{2^{n-p-1}(n-1)! \, k^{n-p-1}}{(n-j_1)\cdot \ldots \cdot(n-j_p)}\frac{1}{(n-1)!} \\
%%%   - \sum_{\substack{1 \leq j_1 < \cdots < j_{p-1} \leq n-2 \\ j_q \in \Z }} \frac{n(n-2)! \, k^{n-p-1}}{(n-j_1)\cdot \ldots \cdot(n-j_{p-1})}\frac{1}{(n-1)!}
%%%   \end{multline*}
\begin{align} 
&2^{n-p-1} \sum_{\substack{1 \leq j_1 < \cdots < j_p \leq n-1 \\ j_q \in \Z }} \frac{1}{(n-j_1) \cdots (n-j_p)} \notag\\
&\qquad - \frac{n}{n-1} \sum_{\substack{1 \leq j_1 < \cdots < j_{p-1} \leq n-2 \\ j_q \in \Z }} \frac{1}{(n-j_1) \cdots(n-j_{p-1})}
\notag
\displaybreak[0]\\
&= 2^{n-p-1} \sum_{\substack{1 \leq j_1 < \cdots < j_{p-1} \leq n-2 \\ j_q \in \Z }} \biggl( \frac{1}{(n-j_1) \cdots (n-j_{p-1})}  \sum_{\substack{j_p=1+j_{p-1}}}^{n-1} \frac{1}{n-j_p} \biggr) \notag\\
&\qquad - \frac{n}{n-1} \sum_{\substack{1 \leq j_1 < \cdots < j_{p-1} \leq n-2 \\ j_q \in \Z }} \frac{1}{(n-j_1) \cdots(n-j_{p-1})}\notag
\displaybreak[0]\\
\label{eq:d2e3}
&=  \sum_{\substack{1 \leq j_1 < \cdots < j_{p-1} \leq n-2 \\ j_q \in \Z }} \biggl( \frac{1}{(n-j_1) \cdots (n-j_{p-1})} \biggl[ 2^{n-p-1} \sum_{\substack{j_p=1 + j_{p-1}}}^{n-1} \frac{1}{n-j_p} - \frac{n}{n-1} \biggr] \biggr).
\end{align}
It is known that the constant in any Ehrhart polynomial is $1$
\cite{Stanley1996Combinatorics-a}, thus we only need to show that Equation
\eqref{eq:d2e3} is positive for $0 \leq p \leq n-2$.  It is sufficient to show
that the square-bracketed term of~\eqref{eq:d2e3},  
\begin{equation} \label{eq:d2e4}
2^{n-p-1} \sum_{\substack{j_p=1 + j_{p-1}}}^{n-1} \frac{1}{n-j_p} - \left( 1 + \frac{1}{n-1} \right),
\end{equation}
is positive for $0 \leq p \leq n-2$ and all $1 \leq j_1 \leq \cdots \leq j_{p-1} \leq n-2$. We can see that $\sum_{\substack{j_p=1 + j_{p-1}}}^{n-1} \frac{1}{n-j_p} \geq 1$. Moreover, since $0 \leq p \leq n-2$ then $2^{n-p-1} \geq 2$ and $1 + \frac{1}{n-1} \leq 2$ since $n \geq 2$,  proving the result.
%%%  Since $1 \leq j_1 < \cdots < j_k \leq n-1$ we can see that $\frac{1}{n-j_k} \geq \frac{1}{n-k}$ and $1 \leq j_1 < \cdots < j_{k-1} \leq n-2$. Therefore Equation \ref{eq:d2e1} is 
%%%  \begin{multline*} 
%%%  \geq \sum_{\substack{1 \leq j_1 < \cdots < j_{k-1} \leq n-2 \\ j_p \in \Z }} \frac{2^{n-k-1}i^{n-k-1}}{(n-j_1)\cdot \ldots \cdot (n-j_{k-1})}\frac{(n-k+1)}{n-k} \\
%%%  - \sum_{\substack{1 \leq j_1 < \cdots < j_{k-1} \leq n-2 \\ j_p \in \Z }} \frac{ni^{n-k-1}}{(n-j_1)\cdot \ldots \cdot(n-j_{k-1})}\frac{1}{(n-1)}
%%%  \end{multline*}
%%%  \begin{equation*} 
%%%  = \left( 2^{n-k-1}\frac{(n-k+1)}{n-k}-\frac{n}{n-1} \right) i^{n-k-1} \left( \sum_{\substack{1 \leq j_1 < \cdots < j_{k-1} \leq n-2 \\ j_p \in \Z }} \frac{1}{(n-j_1)\cdot \ldots \cdot (n-j_{k-1})} \right)
%%%  \end{equation*}
%%%  where
%%%  \begin{equation} \label{eq:d2e2}
%%%   2^{n-k-1}\frac{(n-k+1)}{n-k}-\frac{n}{n-1}  =  2^{n-k-1} \left( 1 + \frac{1}{n-k} \right) - \left( 1 + \frac{1}{n-1} \right).
%%%  \end{equation}
%%%  Now, proving that the coefficient of $i^{n-k-1}$ is positive is
%%%  reduced to proving that Equation \ref{eq:d2e2} is positive for $0 \leq
%%%  k \leq n-1$. When $k=n-1$ Equation \ref{eq:d2e2} becomes
%%%  \begin{equation*}
%%%  2^0 \left( 1 + 1 \right) - \left(1 + \frac{1}{n-1} \right) > 0
%%%  \end{equation*}
%%%  and we can see for all other $k$ that $2^{n-k-1} \geq 2$, $\left( 1 +
%%%  \frac{1}{n-k} \right) \geq 1$, and $\left( 1 + \frac{1}{n-1} \right) <
%%%  2$ where $n \geq 2$, proving the result.
\end{proof}

To present our results about the $h^*$-vector we begin explaining the
details with the following numbers introduced in
\cite{Katzman:math0408038}, which we refer to as the \emph{Katzman
coefficients}:
\begin{definition} \label{def:kat}
    For any positive integers $n$ and $r$ define the coefficients $A_i^{n,r}$ by
    $$ \sum_{i=0}^{n(r-1)} A_i^{n,r}T^i = (1 + T + \cdots + T^{r-1})^n.$$
    We also define the vector $\ve A^{n,r}$ as 
    $\left( A_0^{n,r},A_1^{n,r},\ldots ,A_{n(r-1)}^{n,r} \right)$.
\end{definition}

Looking at the definition of the Katzman coefficient, we see that
$A_j^{n,2} = {n \choose j}$ and $A_j^{n,1} = 0$ unless $j=0$, in which
case we have $A_0^{n,1} = 1$. \\ \indent Below we derive some new and
useful equalities and prove symmetry and unimodality for the Katzman
coefficients. Katzman \cite{Katzman:math0408038} gave an explicit
equation for the $h^*$-vector of uniform matroid polytopes and the
coefficients of their Ehrhart polynomials, although he did not use the same
language. We restate it here for our purposes:

%we state a known lemma (see Corollary 2.9 in \cite{Katzman:math0408038}):
\begin{lemma}[See Corollary 2.9 in \cite{Katzman:math0408038}] \label{katlemma}
    Let $\Po(U^{r,n})$ be the matroid polytope 
    of the uniform matroid of rank $r$ on $n$ elements. Then 
    the $h^*$-polynomial of $\Po(U^{r,n})$ is given by

    \begin{equation} \label{eq:hvec2}
        \sum_{s=0}^{r-1} \sum_{j=0}^s \sum_{k=0}^j \sum_{l \geq k} (-1)^{s +j +k} \left[ {n \choose s} {s \choose j} {j \choose k} A_{(l-k)(r-s)}^{n-j,r-s} \right] T^{l}.
    \end{equation}
    That is, for $\ve h^*(\Po(U^{n,r)})) = (h^*_0,\ldots, h^*_r)$, 
    \begin{equation*} \label{eq:hvec2coeff}
        h^*_l = \sum_{s=0}^{r-1} \sum_{j=0}^s \sum_{k=0}^j(-1)^{s +j +k} {n \choose s} {s \choose j} {j \choose k} A_{(l-k)(r-s)}^{n-j,r-s}.
    \end{equation*} 
    For $r=2$ the $h^*$-polynomial of $\Po(U)^{n,2}$ is
    \begin{equation} \label{eq:dim2}
        \left( \sum_{l \geq 0} {n \choose 2l} T^l \right) - nT .
    \end{equation}
\end{lemma}

%%% TOO MUCH\begin{proof}
%%% TOO MUCH    From Corollary 2.9 in \cite{Katzman:math0408038} we get that the 
%%% TOO MUCH    $h^*$-polynomial of $\Po(U^{d,n})$ is
%%% TOO MUCH        \begin{equation} \label{eq:hvec1}
%%% TOO MUCH            \sum_{n=0}^{r-1} (-1)^s {n \choose s} \sum_{j=0}^s (-1)^j {s \choose j} (1-T)^j \sum_{l\geq 0} A_{l(r-s)}^{n-j,r-s} t^l.
%%% TOO MUCH        \end{equation}
%%% TOO MUCH    We expand the term $(1-T)^j = \sum_{k=0}^{j} {j \choose k} (-1)^kT^k$. 
%%% TOO MUCH    Then equation \eqref{eq:hvec1} becomes
%%% TOO MUCH    %%%\begin{equation*}
%%% TOO MUCH    \begin{align*}
%%% TOO MUCH        & \sum_{n=0}^{r-1} (-1)^s {n \choose s} \sum_{j=0}^s (-1)^j {s \choose j} \sum_{k=0}^{j} {j \choose k} (-1)^kT^k \sum_{l \geq 0} A_{l(r-s)}^{n-j,r-s} T^l \\
%%% TOO MUCH         = & \sum_{s=0}^{r-1} \sum_{j=0}^s \sum_{k=0}^j \sum_{l\geq 0} (-1)^{s +j +k} {n \choose s} {s \choose j} {j \choose k} A_{l(r-s)}^{n-j,r-s} T^{l +k} \\
%%% TOO MUCH         = &\sum_{s=0}^{r-1} \sum_{j=0}^s \sum_{k=0}^j \sum_{l \geq k} (-1)^{s +j +k} {n \choose s} {s \choose j} {j \choose k} A_{(l-k)(r-s)}^{n-j,r-s} T^{l}.
%%% TOO MUCH    \end{align*}
%%% TOO MUCH    %%%\end{equation*}
%%% TOO MUCH\end{proof}

The
following lemma is a direct consequence of Corollary 2.9 in
\cite{Katzman:math0408038}:

\begin{lemma}%[See Corollary 2.9 in \cite{Katzman:math0408038}]
    Let $\Po(U^{2,n})$ be the matroid polytope 
    of the uniform matroid of rank $2$ on $n$ elements. 
    The $h^*$-vector of $\Po(U^{2,n})$ is unimodal.
\end{lemma}

%%%  TOO MUCH %%% \begin{proof}
%%%  TOO MUCH %%%     Looking at equation \eqref{eq:dim2} and letting $n=2$ we simply get 
%%%  TOO MUCH %%%     $1$ as the $h^*$-polynomial. When $n=3$ we get $1 + 2t$ 
%%%  TOO MUCH %%%     as the $h^*$-polynomial. If $n > 3$, then 
%%%  TOO MUCH %%%     \begin{equation*}
%%%  TOO MUCH %%%         \left({n \choose 0}, {n \choose 2},\dots, {n \choose n - (n \mod 2)}\right)
%%%  TOO MUCH %%%     \end{equation*}
%%%  TOO MUCH %%%     is unimodal and the difference between indices is greater 
%%%  TOO MUCH %%%     than 1. Hence the following is unimodal:
%%%  TOO MUCH %%%     \begin{equation*}
%%%  TOO MUCH %%%         \left({n \choose 0}, {n \choose 2} - 1,{n \choose 4}, {n \choose 6}, \dots, {n \choose n - (n \mod 2)}\right) 
%%%  TOO MUCH %%%     \end{equation*}
%%%  TOO MUCH %%% \end{proof}

The rank 2 case is an interesting example already. Although the 
$h^*$-vector is unimodal, it is not always symmetric. Next we
present some useful lemmas, the first a combinatorial
description of the Katzman coefficients.

\begin{lemma}
  For $i=0,\dots,n(r-1)$ we have
    \begin{equation} \label{eq:katcomb}
        A_i^{n,r} =\sum_{\substack{0a_0 + 1a_1 + \dots + (r-1)a_{r-1} = i\\ a_0 + a_1 + \dots + a_{r-1} = n }} { n \choose a_0,a_1, \cdots ,a_{r-1}} 
    \end{equation}
    where $a_0,\dots,a_{r-1}$ run through non-negative integers.
\end{lemma}

\begin{proof}
    Using the multinomial formula \cite{Stanley1997Enumerative-Com} we have

    \begin{align*}
     \sum_{i=0}^{n(r-1)} A_i^{n,r}T^i = & (1 + T + \cdots + T^{r-1})^n  \\
    = &\sum_{\substack{ a_0 + a_1 + \dots + a_{r-1} = n}} { n \choose a_0,a_1, \cdots ,a_{r-1}} 1^{a_0}T^{a_1}T^{2a_2}\dots T^{(r-1)a_{r-1}} \\
    = &\sum_{\substack{ a_0 + a_1 + \dots + a_{r-1} = n}} { n \choose a_0,a_1, \cdots ,a_{r-1}} T^{0a_0 + 1a_1 + \cdots + (r-1)a_{r-1}}.
    \end{align*}
    By grouping the powers of $T$ we get equation \eqref{eq:katcomb}.
\end{proof}

Next we present a generalization of a property of the binomial
coefficients. The following lemma relates the Katzman coefficients to
Katzman coefficients with one less element.
\begin{lemma} \label{lem:dimrel}
    \begin{equation} \label{eq:dimrel}
        A_i^{n,r} = \sum_{k=i-r+1}^i A_k^{n-1,r}
    \end{equation}
    where we define $A_p^{n-1,r} := 0$ when $p < 0$ or $p > (n-1)(r-1)$.
\end{lemma}

\begin{proof}
    %%% First define $A_p^{n-1,r} := 0$ when $p < 0$ and $p > (n-1)(r-1)$. Then
    \begin{equation*}
    \begin{split}
    & \sum_{i = 0}^{n(r-1)} A_i^{n,r} T^i  \\
    & = (1 + T + \dots + T^{r-1})^n = (1 + T + \dots + T^{r-1})^{n-1}(1 + T + \dots + T^{r-1})   \\
    & = \left( \sum_{i = 0}^{(n-1)(r-1)} A_i^{n-1,r} T^i \right)  (1 + T + \dots + T^{r-1})  \\
    & = \sum_{i = 0}^{(n-1)(r-1)} A_i^{n-1,r} T^i + \sum_{i = 0}^{(n-1)(r-1)} A_i^{n-1,r} T^{i+1} + \cdots + \sum_{i = 0}^{(n-1)(r-1)} A_i^{n-1,r} T^{i + r-1} \\
    & = \sum_{i = 0}^{(n-1)(r-1)} A_i^{n-1,r} T^i + \sum_{i = 1}^{(n-1)(r-1) +1} A_{i-1}^{n-1,r} T^{i} + \cdots + \sum_{i = r-1}^{(n-1)(r-1) + r-1} A_{i-r+1}^{n-1,r} T^{i} \\
    & = \sum_{i=0}^{(n-1)(r-1)+r-1} \left( A_i^{n-1,r} + A_{i-1}^{n-1,r} + \cdots + A_{i-r+1}^{n-1,r} \right) T^i.
    \end{split}
    \end{equation*}
    Thus we get equation \eqref{eq:dimrel}.
    %%%   \begin{equation*}  
    %%%   A_i^{n,r} = \sum_{k = i-r+1}^i A_k^{n-1,r}.
    %%%   \end{equation*}
\end{proof}

%%%  The following is a relation between the Katzman coefficients.
The following lemma relates the Katzman coefficients of rank~$r$
with those of rank~$r-1$.
\begin{lemma} \label{lem:rankrel}
    \begin{displaymath}
    \sum A_i^{n,r}T^i = \sum_{k=0}^n { n \choose k} T^k \left( \sum_{l=0}^{k(r-2)} A_l^{k,r-1}T^l \right)
    \end{displaymath}
    or in other words
    \begin{equation*} \label{eq:rankrel}
    A_i^{n,r} \; = \; \sum_{\substack{k+l=i \\ 0 \leq k \leq n \\ 0 \leq l \leq k(r-2)}} {n \choose k}  A_l^{k,r-1}.
    \end{equation*}
\end{lemma}

\begin{proof}
    From Definition \ref{def:kat}
    \begin{align*}
        \sum_{i=0}^{n(r-1)} A_i^{n,r}T^i = & (1 + T + \cdots + T^{r-1})^n = (1 + \big[T + \cdots + T^{r-1}\big])^n \\
        = & \sum_{k=0}^n { n \choose k} \big[T + \cdots + T^{r-1}\big]^k = \sum_{k=0}^n { n \choose k} T^k\big[1 + \cdots + T^{r-2}\big]^k \\
        = & \sum_{k=0}^n { n \choose k} T^k \left( \sum_{l=0}^{k(r-2)} A_l^{k,r-1}T^l \right).
    \end{align*}
\end{proof}

\begin{lemma} \label{lem:katsymuni}
    The Katzman coefficients are unimodal and symmetric in the index $i$. That is, the vector $\left( A_0^{n,r}, A_1^{n,r}, \dots, A_{n(r-1)}^{n,r} \right)$ is unimodal and symmetric.
\end{lemma}

\begin{proof}
    We first prove symmetry. Considering equation \eqref{eq:katcomb} we assume that
    \begin{equation*} 
        0a_0 + 1a_1 + \dots + (r-1)a_{r-1} = i \qquad \text{and} \qquad a_0 + a_1 + \dots + a_{r-1} = n. 
    \end{equation*}
    These two assumptions imply that
    \begin{align*}
    &\big((r-1)-0\big)a_0 + \big((r-1)-1\big)a_1 + \dots + \big((r-1)-(r-1)\big)a_{r-1}  \\
     = & (r-1)a_0 + (r-1)a_1 + \dots + (r-1)a_n -0a_0 - 1a_1 - \dots - (r-1)a_{r-1} \\ 
     = & (r-1)n - i.
    \end{align*}
    Therefore
    \begin{align*}
      A_i^{n,r} = & \sum_{\substack{0a_0 + 1a_1 + \dots + (r-1)a_{r-1} = i\\ a_0 + a_1 + \dots + a_{r-1} = n }} { n \choose a_0,a_1, \cdots ,a_{r-1}} & \\
     = & \sum_{\substack{(r-1)a_0 + (r-2)a_1 + \dots + 0a_{r-1} = i\\ a_0 + a_1 + \dots + a_{r-1} = n }} { n \choose a_0,a_1, \cdots ,a_{r-1}} = A_{n(r-1)-i}^{n,r}.
    \end{align*}

    %%%Hence we have symmetry. \\

    To prove unimodality we proceed by induction on $n$, where $r$ is fixed.
    First, $\ve A^{1,r}$ is unimodal. Assume for $n-1$ that $\ve
    A^{n-1,r} = \left(
            A_0^{n-1,r},A_1^{n-1,r},\cdots,A_{(n-1)(r-1)}^{n-1,r}\right)$ is
    unimodal. Using equation \eqref{eq:dimrel} and the fact that $\ve A^{n-1,r}$
    is symmetric we get that $\ve A^{n,r}$ is unimodal. To help see this,
    one can view equation \eqref{eq:dimrel} as a sliding window over $r$ elements of
    the vector $\ve A^{n-1,r}$, that is, $A_i^{n,r}$ is equal to the sum
    of the $r$ elements in a window over the vector $\ve A^{n-1,r}$. As
    the window slides up the vector $\ve A^{n-1,r}$, the sum will
    increase.  When the window is on the center of $\ve A^{n-1,r}$
    symmetry and unimodality of $\ve A^{n-1,r}$ will imply unimodality of
    $\ve A^{n,r}$.
%%%       To prove unimodality, it is helpful to view equation \eqref{eq:dimrel} in the following way, that $\ve A^{n,r}$ is the sum of a sliding window over the vector $\ve A^{n-1,r}$,
%%%       $$ \undersum{A_0^{n-1,r}}{A_0^{n,r}}{A_1^{n-1,r}}{A_{(n-1)(r-1)}^{n-1,r}} $$
%%%   
%%%       $$ \undersum{ A_0^{n-1,r}, A_1^{n-1,r}}{A_1^{n,r}}{ A_2^{n-1,r}}{A_{(n-1)(r-1)}^{n-1,r} } $$
%%%   
%%%       $$ \undersum{ A_0^{n-1,r}, A_1^{n-1,r}, A_2^{n-1,r}}{ A_2^{n,r}}{ A_3^{n-1,r}}{A_{(n-1)(r-1)}^{n-1,r} } $$
%%%   
%%%       $$ \vdots $$
%%%       $$ \left( A_0^{n-1,r}, \lrots , A_{i-r}^{n-1,r}, \right. \underbrace{A_{i-r+1}^{n-1,r}, A_{i-r+2}^{n-1,r}, \ldots , A_{i}^{n-1,r}}_{ \begin{array}{c} + \\ =  A_i^{n,r} \end{array} } \left. , A_{i+1}^{n-1,r} , \ldots , A_{(n-1)(r-1)}^{n-1,r}\right) $$
%%%   
%%%       $$ \left( A_0^{n-1,r}, \ldots , A_{i-r}^{n-1,r} , A_{i-r+1}^{n-1,r} , \right. \underbrace{ A_{i-r+2}^{n-1,r} , \ldots , A_{i}^{n-1,r} , A_{i+1}^{n-1,r} }_{\begin{array}{c} + \\  = A_{i+1}^{n,r} \end{array}} \left. , \ldots , A_{(n-1)(r-1)}^{n-1,r}\right) $$
%%%   
%%%       $$ \vdots $$
%%%       \indent Let $r$ be fixed. We proceed by induction on $n$. 
%%%       Firstly $A_i^{1,r}$ is unimodal. Assume for $n-1$ that 
%%%       $\left( A_0^{n-1,r},A_1^{n-1,r},\cdots,A_{(n-1)(r-1)}^{n-1,r}\right)$
%%%       is unimodal, which we already proved is symmetric. 
%%%       Then, symmetry and equation \eqref{eq:dimrel} imply unimodality of 
%%%       $A^{n,r}$. \davecomment{(Try to consolidate/shrink this proof)}
\end{proof}

Now we use the explicit equation for the $h^*$-vector of uniform
matroid polytopes to prove partial unimodality of rank $3$ uniform
matroids. First we note that the coefficient of $T^l$, $h^*_l$, in
equation \eqref{eq:hvec2} is
\begin{displaymath}
h^*_l = \sum_{s=0}^{r-1} \sum_{j=0}^s \sum_{k=0}^j(-1)^{s +j +k} {n \choose s} {s \choose j} {j \choose k} A_{(p-k)(r-s)}^{n-j,r-s}.
\end{displaymath}

Letting the rank $r = 3$, and using equation \eqref{eq:hvec2}, we 
get the $h^*$-polynomial (which is grouped by values of $s$ 
from \eqref{eq:hvec2}),
\begin{multline*}
    \sum_{l \geq 0} \left[ {n \choose 0 } A_{3l}^{n,3} + {n \choose 1} \left( -A_{2l}^{n,2} + A_{2l}^{n-1,2} - A_{2(l-1)}^{n-1,2} \right) + \right.  \\
  + \left.  {n \choose 2} \left( A_{l}^{n,1} - 2A_{l}^{n-1,1} + 2A_{l-1}^{n-1,1} + A_{l}^{n-2,1} - 2A_{l-1}^{n-2,1} + A_{l-2}^{n-2,1} \right) \right] T^l.
\end{multline*}
Now using that $A_{i}^{n,2} = { n \choose i}$ and $A_{i}^{n,1} = \delta_0(i)$, where
$\delta_j(p) = \left\{
\begin{array}{lll}
1 & \text{if } p = j \\
0 & \text{else}
\end{array}
\right.$ ,
 we get
\begin{multline*}
    \sum_{l \geq 0}  \left[ A_{3l}^{n,3} + n \left( - { n \choose 2l} \right. \left. + {n-1 \choose 2l} -  {n-1 \choose 2l - 2} \right) + \right.  \\
    \left. + {n \choose 2} \left(  \delta_0(l) - 2 \delta_0(l) + 2 \delta_{1}(l) + \delta_0(l) - 2 \delta_{1}(l) + \delta_2(l)   \right) \right] T^l
\end{multline*}
\begin{equation*}
    = \sum_{l \geq 0}  \left[ A_{3l}^{n,3} + n \left( - { n \choose 2l} + {n-1 \choose 2l} - {n-1 \choose 2l - 2} \right) + \delta_2(l) {n \choose 2}  \right] T^l.
\end{equation*}
Using properties of the binomial coefficients, we see that
\begin{align*}
     \left[ -{n \choose 2l} \right]  + {n-1 \choose 2l} - {n-1 \choose 2l - 2} = & \left[ -{n-1 \choose 2l} - {n-1 \choose 2l -1} \right] + {n-1 \choose 2l} - {n-1 \choose 2l - 2} \\
= & - {n-1 \choose 2l -1} - {n-1 \choose 2l - 2} \\
= & -{n \choose 2l - 1}.
\end{align*}
So the $h^*$-polynomial of rank three uniform matroid polytopes is
\begin{equation} \label{hvecrank3}
     \sum_{l \geq 0}  \left[ A_{3l}^{n,3} - n{n \choose 2l - 1}  + \delta_2(l) {n \choose 2}  \right] T^l.
\end{equation}
Using Lemma \ref{eq:dimrel}, the coefficient of $T^l$, if $3l \leq n$, is
\begin{equation} \label{hvecrank3eq2}
    h^*_l = \sum_{\substack{k+p = 3l \\ 0 \leq p \leq k \leq n}} {n \choose k}{k \choose p} - n{n \choose 2l-1} + \delta_2(l){n \choose 2}
\end{equation}
\begin{multline*}
= \left[ {n \choose 3l}{3l \choose 0} + {n \choose 3l-1}{ 3l-1 \choose 1} + \cdots+ {n \choose 3l - \lfloor 3l/2 \rfloor}   {3l - \lfloor 3l/2 \rfloor \choose \lfloor 3l/2 \rfloor} \right] + \\
   - n{n \choose 2l-1} + \delta_2(l){n \choose 2}.
\end{multline*}

Next we show that when $g$ is fixed $A_g^{n,3}$ is a polynomial of degree $g$
in the indeterminate $n$, with positive leading coefficient. Assume  $g \leq
n$. Considering Lemma \ref{lem:rankrel} and when $g \leq n$, 
\begin{equation}
   A_{g}^{n,3} = \sum_{\substack{k+p = g \\ 0 \leq p \leq k \leq n}} {n \choose k}{k \choose p}
   \label{katdim3}
\end{equation}
where $n \choose q$ is a polynomial
of degree $q$ with positive leading coefficient. The highest
degree polynomial in the sum is $n \choose g$, a degree $g$ 
polynomial. Hence $A_g^{n,3}$ is a polynomial of degree $g$ 
in the indeterminate $n$, with positive leading coefficient. 
If $g \geq n$, then $A_g^{n,3} = A_{n-g}^{n,3}$ since
the Katzman coefficients are symmetric by Lemma \ref{lem:katsymuni}.

\begin{proof}[Proof of Theorem \ref{partialuni} Part (2)]
     Let $I$ be a non-negative integer. From above we see that $A_g^{n,3}$ is a
     degree $g$ polynomial in the indeterminate $n$, with positive leading
     coefficient.  Equation \eqref{hvecrank3} is the $h^*$-polynomial of
     $U^{3,n}$, which is a sum of polynomials in $n$, the highest degree
     polynomial being $A_{3l}^{n,3}$, a polynomial of degree $3l$. So, $h^*_l -
     h^*_{l-1}$ is the difference of a degree $3l$ and $3(l-1)$ polynomial,
     hence $h^*_l - h^*_{l-1}$ is a degree $3l$ polynomial with positive
     leading coefficient.  For sufficiently large $n$, call it $n(I)$, $h^*_l
     - h^*_{l-1}$ is positive for $0 \leq l \leq I$. Hence, the $h^*$-vector of
     $U^{3,n}$ is non-decreasing up to the index $I$ for $n \geq n(I)$.
\end{proof}

One might ask if equation \eqref{katdim3} has a simpler form. We ran the
\emph{WZ} algorithm on our expression, which proved that equation
\eqref{katdim3} can not be written as a linear combination of a fixed number of
hypergeometric terms (\emph{closed form}) \cite{Petkovsek1996AB}. There is
still the possibility that this expression has a simpler form, though not a
closed form as described above.

\section{Unimodular Simplices of Matroid Polytopes}
A $k$-\emph{simplex} $S =\{s_1,\ldots,s_{k+1}\}$ is the collection of $k+1$ affinely independent points in $\R^n$ and we say $\conv(S) := \conv(\ve s_1,\ldots,\ve s_{k+1})$. When we refer to a simplex $S$ we consider it both as a collection of points and as a convex polytope when the context is clear. Let $\Po \subseteq \R^n$ be a polytope with vertices $\ve v_1,\ldots,\ve v_p$. A \emph{triangulation} $\Tr$ of $\Po$ is a collection of simplices on the vertices $\ve v_1,\ldots, \ve v_p$ of $\Po$ such that 
\begin{itemize}
    \item[(i)]   If $A \in \Tr$ then all faces of $A$ are in $\Tr$,
    \item[(ii)]  If $A,B \in \Tr$ then $A \cap B$ is a face of both $A$ and $B$,
    \item[(iii)] $\bigcup_{A \in \Tr} \conv(A) = \Po$.
\end{itemize}

 It is sufficient to list only the highest dimensional simplices of a triangulation. A triangulation $\Tr$ is \emph{unimodular} if the volume of every highest dimensional simplex of $\Tr$ are the same. We also say a simplex is unimodular if its normalized volume is one \cite{De-Loera2006Triangulations}. A vector in $R^n$ is a $\pm1${\it -vector} if it contains $+1$ and $-1$ and the rest are zeros.

Recall that two elements $s,t \in S$ of a matroid $\M=(S,\In)$ are \emph{connected} if there exists a circuit that contains $s$ and $t$. This defines an equivalence relation on $S$ and moreover $\M$ can be written as a direct sum of restrictions to the connected components of $\M$ \cite{Oxley1992Matroid-Theory}. This implies a matroid polytope $\Po_\M$ can be written as a direct product of connected matroid polytopes. We note that if a polytope $\Po$ is a direct product of polytopes $\mathcal Q_1, \ldots, \mathcal Q_k$ which each have a unimodular triangulation, then $\Po$ has a unimodular triangulation \cite{De-Loera2006Triangulations}. Hence, we need only consider connected matroids when investigating unimodular triangulations.
 
We wish to find sufficient conditions such that a $(n-1)$-simplex of a connected matroid polytope $\Po_\M$ is unimodular. We will see later that a simplex of $\Po_\M$ is unimodular if and only if the determinant of the vertices is $\rank(\M)$.

Neil White proposed an algebraic conjecture on the toric ideal determined by a matroid \cite{White1980A-unique-exchan}. Let $\F$ be a field and define the polynomial ring $S_\M := \F[\ve y_B \mid B \in \B_\M]$. Let the toric ideal $I_\M$ be defined as the kernel of the homomorphism $S_\M \rightarrow \F[x_1,\ldots,x_n]$ given by $\ve y_B \mapsto \prod_{i\in B} x_i$. Given any pair of bases $B_1, B_2 \in \B$ the \emph{symmetric exchange property} states that for every $x \in B_1$ there exists $y \in B_2$ such that $B_1 - x + y$ and $B_2 - y + x$ are bases and it is said that $x \in B_1$ \emph{double swaps} into $B_2$. If $B_1 - x + y$ and $B_2 - y + x$ are bases they are called a \emph{double swap}.  White proposed the following conjecture:
\begin{conjecture}[\cite{White1980A-unique-exchan}]
For any matroid $\M$, the toric ideal $I_\M$ is generated by the quadratic binomials $\ve y_{B_1}\ve y_{B_2} - \ve y_{D_1}\ve y_{D_2}$ such that the pair of bases $D_1, D_2$ can be obtained from the pair $B_1, B_2$ by a double swap.
\end{conjecture}

A stronger conjecture is that the previous quadratic binomials are in fact a Gr\"obner bases for some term ordering. This implies that there exists a regular unimodular triangulation of any matroid polytope.
The following is a geometric variation of White's famous conjecture:
\begin{conjecture} \label{conj:white}
Let $\M$ be a matroid of rank $r$ on $n$ elements. There exists a unimodular triangulation of $\Po_\M$.
\end{conjecture}
In \cite{Sturmfels1996Grobner-Bases-a} it was shown that all uniform matroids have a regular unimodular triangulation. A partial proof of White's original conjecture was given in \cite{Blasiak2005The-Toric-Ideal} for graphical matroids, but \autoref{conj:white} has yet to be proved in general.

A \emph{covering} of $\Po$ is a triangulation without condition (ii). That is, the simplices of a covering can intersect in any fashion, as long as the polytope is covered. See \autoref{fig:covering}. If asking for a unimodular triangulation is too much, we propose a weaker conjecture:  
\begin{figure}[!htb]
    \scalebox{1}{\ifpdf
    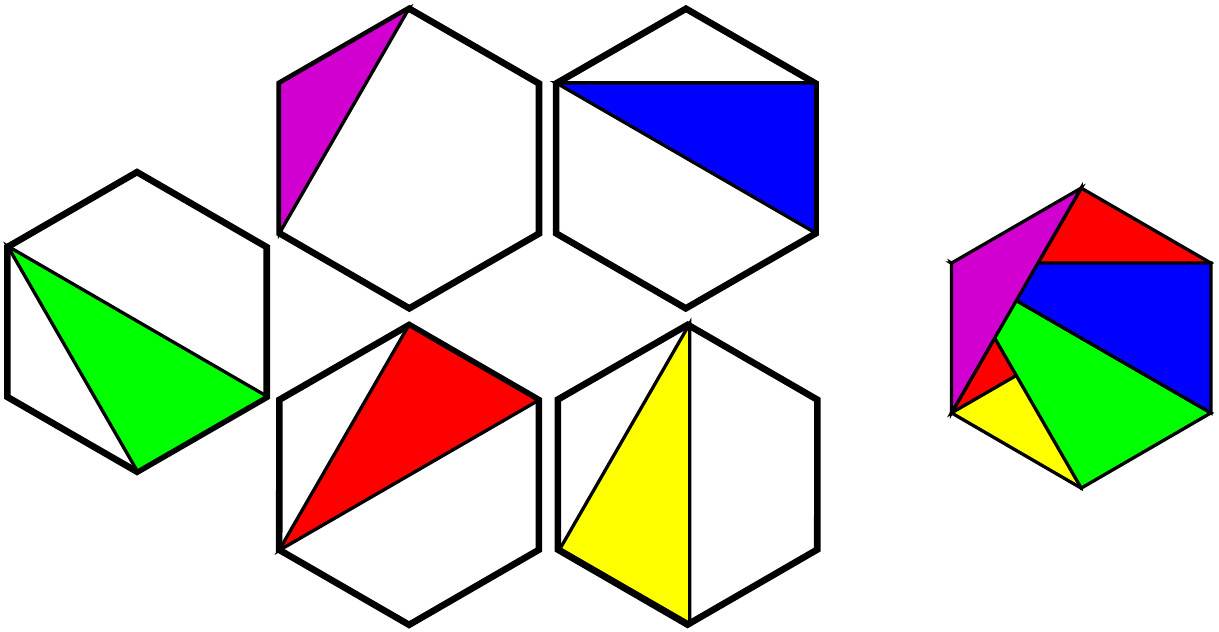
    \else
    \input{covering.pstex_t}
    \fi}
    \caption{A covering of a hexagon.} \label{fig:covering}
\end{figure}
\begin{conjecture} \label{conj:unicovering}
Let $\M$ be a matroid of rank $r$ on $n$ elements. There exists a unimodular covering of $\Po_\M$.
\end{conjecture}
Coverings are important to matroid polytopes. If a unimodular covering exists for a connected matroid polytope $\Po_\M$ it would imply its Caratheodory rank is $n$. In other words, every integer point in the cone of $\Po_\M$ can be written as an integral linear combination of vertices of $\Po_\M$. That is, the vertices of $\Po_\M$ are a Hilbert basis for the cone of $\Po_\M$. See \cite{Cook1986An-integer-anal,Sebo1990Hilbert-bases-C,Pina2003Improved-bound-}.

We now focus on sufficient conditions for the existence unimodular simplices of matroid polytopes and begin by defining two graphs on any collection of incidence vectors of a matroid polytope.

\begin{definition}
Let $\M$ be a matroid on $n$ elements. Given a subset $X \subseteq \Inc(\B_\M) \subseteq \R^n$ we define two graphs $G$ and $g$ on $n$ nodes. The graph $G(X)$ has an edge $(i,j)$ for every pair of bases $\ve X_i,\ve X_j \in X$ that are adjacent, i.e. $\ve X_i - \ve X_j = \ve e_s - \ve e_t$ for some $s,t$. The graph $g(X)$ has an edge $(i,j)$ every pair of adjacent bases $\ve X_s, \ve X_t \in X$ such that they differ in position $i$ and $j$, i.e. $\ve X_s - \ve X_j = \ve e_i - \ve e_j$. Succinctly we define $G(X)$ and $g(X)$ as:
\begin{equation*}
    G(X) := ([n],\, \{\, (i,j) \mid \ve X_i,\ve X_j \in X,\, \ve X_i - \ve X_j = \ve e_s - \ve e_t \, \})
\end{equation*}
\begin{equation*}
    g(X) := ([n],\, \{\, (i,j) \mid \ve X_s,\ve X_t \in X,\, \ve X_s - \ve X_t = \ve e_i - \ve e_j \, \}) 
\end{equation*}
where $\ve e_i,\ve e_j, \ve e_s, \ve e_t \in \R^n$ are standard unit vectors.
\end{definition}
%We see that $G(X)$ and $g(X)$ are defined on the same number of nodes. Yet $G(X)$ contains an edge for each pair of incidence vectors of $X$ that are adjacent on the $1$-skeleton of $\Po_\M$ and $g(X)$ contains an edge for the elements of $[n]$ where a pair of adjacent incidence vectors differ. 
Throughout this section we will take $X$ to mean a set of incidence vectors as well as a $0/1$ matrix with rows as the incidence vectors.  The following proposition is easy to see using elementary row sums on $X$ and cofactor expansion on the resulting last row.
\begin{proposition}
Let $\M$ be a connected matroid of rank $r$ on $n$ elements and $X \subset \Inc(\B_\M)$ be $n$ linearly independent incidence vectors of $\Po_\M$. Then $|det(X)| = kr$, $k \in \N$. 
\end{proposition}
\begin{example}
Consider the following incidence vectors $X$ in row form. We label the rows $\mathbf 1,\ldots,\mathbf 6$ and columns $\mathbf a,\ldots, \mathbf f$.\\
\begin{center}                
\begin{minipage}{0.4\linewidth}
%$$ \left( \begin{array}{cccccc}
%1 & 1 & 0 & 0 & 1 & 0\\ 
%1 & 1 & 0 & 0 & 0 & 1\\
%1 & 0 & 1 & 1 & 0 & 0\\
%0 & 1 & 1 & 1 & 0 & 0\\
%0 & 0 & 1 & 0 & 1 & 1\\
%0 & 0 & 0 & 1 & 1 & 1\\
$$
\bordermatrix{ & \mathbf{a} & \mathbf{b} & \mathbf{c} & \mathbf{d} & \mathbf{e} & \mathbf{f} \cr
\mathbf{1} & 1 & 1 & 0 & 0 & 1 & 0\cr 
\mathbf{2} & 1 & 1 & 0 & 0 & 0 & 1\cr
\mathbf{3} & 1 & 0 & 1 & 1 & 0 & 0\cr
\mathbf{4} & 0 & 1 & 1 & 1 & 0 & 0\cr
\mathbf{5} & 0 & 0 & 1 & 0 & 1 & 1\cr
\mathbf{6} & 0 & 0 & 0 & 1 & 1 & 1}
$$
%\end{array}\right) $$
\end{minipage}
\hskip 8pt
\begin{minipage}{0.3\linewidth}
$$ G(X) = ([6],\{(\mathbf 1,\mathbf 2),(\mathbf 3,\mathbf 4),(\mathbf 5,\mathbf 6)\})$$
$$ g(X) = ([6],\{(\mathbf a,\mathbf b),(\mathbf c,\mathbf d),(\mathbf e,\mathbf f)\})$$
\end{minipage}
\end{center}
\end{example}
We can equivalently define $g(X)$ as the graph whose incidence matrix is given by
\begin{equation*}
   \left\{\, \ve X_s - X_t \mid \ve X_s, \ve X_t \in X,\, \ve X_s - \ve X_t = \ve e_i - \ve e_j \, \right\} 
\end{equation*}
where $\ve e_i, \ve e_j \in \R^n$ are standard unit vectors. We will use this fact for the following theorem.

\begin{theorem} \label{thm:genconn}
Let $\M$ be a connected matroid of rank $r$ on $n$ elements and $X \subset \Inc(\B_\M)$ be $n$ linearly independent incidence vectors of $B_\M$. Then
\begin{itemize}
    \item[(i)] $G(X)$ and $g(X)$ have the same number of connected components $c$,
    \item[(ii)] Let $X = \{\ve v^1,\ldots, \ve v^n \}$ and without loss of generality assume $\{1, \ldots,l_1-1\}$, $\{l_1, \ldots, l_2-1\}$, $\cdots$, $\{l_{c-1}, \ldots, n\}$,  are the connected components of $g(X)$. Then

\begin{equation} \label{eq:genconn}
    |\det(X)| = \left| \det \left (\begin{array}{cccc} 
    \sum_{i=1}^{l_1-1}  v^{m_1}_i & \sum_{i=l_1}^{l_2-1}  v^{m_1}_i & \cdots & \sum_{i=l_{c-1}}^{n}  v^{m_1}_i \\
    \sum_{i=1}^{l_1-1}  v^{m_2}_i & \sum_{i=l_1}^{l_2-1}  v^{m_2}_i & \cdots & \sum_{i=l_{c-1}}^{n}  v^{m_2}_i \\
    \vdots & \vdots & \ddots & \vdots \\
    \sum_{i=1}^{l_1-1}  v^{m_c}_i & \sum_{i=l_1}^{l_2-1}  v^{m_c}_i & \cdots & \sum_{i=l_{c-1}}^{n}  v^{m_c}_i \\
    \end{array} \right) \right| %\in \R^{c \times c}
\end{equation}
where $\ve v^{m_1},\ldots, \ve v^{m_c}$ are representative incidence vectors of the connected components of $G(X)$.
\end{itemize}
%where the second matrix is a $(c) \times (c)$ matrix.
\end{theorem}
\begin{proof}
%{\it \bf (i)} hmmj
%{\it \bf (ii)} I don't know 
We offer a constructive proof of the theorem. Assume $G(X)$ has $p$ connected component and $g(X)$ has $c$ connected components. First let $\ve v^{m_1},\ldots, \ve v^{m_p}$ be representative incidence vectors of the connected components of $G(X)$.  Performing elementary row operations such as row addition, difference or reordering will only change the sign of the determinant. First we organize the rows of $X$ by the $p$ connected components of $G(X)$.
\begin{equation*} \tiny
X = \left(
    \begin{array}{ccc}
    v_1^{m_1} & \cdots & v_{n}^{m_1}  \\
    \vdots & \ddots & \vdots \\
    \hline
    v_1^{m_2} & \cdots & v_{n}^{m_2}  \\
    \vdots & \ddots & \vdots \\
    \hline
    \vdots & \ddots & \vdots \\
    \hline
    v_1^{m_p} & \cdots & v_{n}^{m_p}  \\
    \vdots & \ddots & \vdots \\
    \end{array}
\right)   
\end{equation*}
Consider a connected component of $G(X)$ with representative incidence vector $\ve v^{m_i}$. Let $\ve v^{m_i}$ be the root of a spanning tree $T$ of the corresponding component of $G(X)$. For all non-root vectors in the component determined by $\ve v^{m_i}$ perform the row operation $\ve v^g := \ve v^g - \ve v^h$ where $h$ is the parent of $g$. Now all the rows except $\ve v^{m_i}$ of the component of $G(X)$ are $\pm 1$-vectors. Repeat for each connected component of $G(X)$. (By $\cdots \pm 1 \cdots$ we mean some $\pm 1$ vector).
\begin{center}
\begin{minipage}{0.4\linewidth}
\begin{equation*} \tiny
 \left(
    \begin{array}{ccc}
    v_1^{m_1} & \cdots & v_{n}^{m_1}  \\
     & \cdots    \pm 1   \cdots       \\
%     \cdots   & \pm 1  & \cdots       \\
              & \vdots &              \\
     & \cdots    \pm 1   \cdots       \\
    \hline
    v_1^{m_2} & \cdots & v_{n}^{m_2}  \\
    \vdots & \ddots & \vdots \\
    \hline
    \vdots & \ddots & \vdots \\
    \hline
    v_1^{m_p} & \cdots & v_{n}^{m_p}  \\
    \vdots & \ddots & \vdots \\
    \end{array}
\right)   \rightarrow
\end{equation*}
\end{minipage}
\hskip 5pt
\begin{minipage}{0.4\linewidth}
\begin{equation*} \tiny
\left(
    \begin{array}{ccc}
    v_1^{m_1} & \cdots & v_{n}^{m_1}  \\
     & \cdots    \pm 1   \cdots       \\
%     \ldots   & \pm 1  & \ldots       \\
              & \vdots &              \\
     & \cdots    \pm 1   \cdots       \\
    \hline
    v_1^{m_2} & \cdots & v_{n}^{m_2}  \\
     & \cdots    \pm 1   \cdots       \\
%     \ldots   & \pm 1  & \ldots       \\
              & \vdots &              \\
     & \cdots    \pm 1   \cdots       \\
    \hline
    \vdots & \ddots & \vdots \\
    \hline
    v_1^{m_p} & \cdots & v_{n}^{m_p}  \\
     & \cdots    \pm 1   \cdots       \\
%     \ldots   & \pm 1  & \ldots       \\
              & \vdots &              \\
     & \cdots    \pm 1   \cdots       \\
    \end{array}
\right)   
\end{equation*}
\end{minipage}
\end{center}
Next, due to the assumption on the ordering of the connected components of $g(X)$ we partition the columns of matrix $X$ by the components of $g(X)$. We also reorder the rows with the representatives vectors $\ve v^{m_1},\ldots,\ve v^{m_p}$ on top and the $\pm 1$-vectors by their corresponding components of $g(X)$.
\begin{equation*} \tiny
 \left(
    \begin{array}{c|c|c|c}
    v_1^{m_1} \cdots v_{l_1-1}^{m_1} & v_{l_1}^{m_1} \cdots v_{l_2-1}^{m_1} & \cdots & v_{l_{c-1}}^{m_1} \cdots v_n^{m_1} \\
    v_1^{m_2} \cdots v_{l_1-1}^{m_2} & v_{l_1}^{m_2} \cdots v_{l_2-1}^{m_2} & \cdots & v_{l_{c-1}}^{m_2} \cdots v_n^{m_2} \\
    \vdots          &  \vdots         & \ddots &  \vdots            \\
    v_1^{m_p} \cdots v_{l_1-1}^{m_n} & v_{l_1}^{m_p} \cdots v_{l_2-1}^{m_n} & \cdots & v_{l_{c-1}}^{m_p} \cdots v_n^{m_n} \\
    \hline
    \cdots \pm 1 \cdots &   0 \cdots 0      & \ldots        & 0 \cdots 0          \\
    \vdots              &    \vdots         & \ddots        &   \vdots            \\
    \cdots \pm 1 \cdots &  0 \cdots 0       & \ldots        & 0 \cdots 0          \\
    \hline
       0 \cdots 0       & \cdots \pm 1 \cdots  & \ldots        & 0 \cdots 0          \\
        \vdots          & \vdots               & \ddots        &   \vdots            \\
      0 \cdots 0        & \cdots \pm 1 \cdots  & \ldots        & 0 \cdots 0          \\
    \hline
    \ldots              &   \ldots          & \ddots        & \ldots              \\
    \hline
    0 \cdots 0 &   0 \cdots 0      & \ldots        & \cdots \pm 1 \cdots          \\
      \vdots   &    \vdots         & \ddots        & \vdots                       \\
    0 \cdots 0 &  0 \cdots 0       & \ldots        & \cdots \pm 1 \cdots          \\
    \end{array}
\right)   
\end{equation*}

We note that for isolated nodes of $g(X)$ there will be no corresponding $\pm 1$-vectors, but this will not interfere with this proof. Consider any submatrix of $\pm 1$-vectors determined by a connected component $C$ of $g(X)$. First, if there are $k$ many nodes in $C$, then there will be $k-1$ $\pm 1$-vectors since $X$ is linearly independent and these vectors are a spanning tree for an orientation of the component $C$. Using the connectedness of the component $C$, elementary row addition and the fact that $\pm 1$-vectors sum to $0$ we can rewrite the parts of each representative incidence vector $\ve v_{m_1},\ldots,\ve v_{m_p}$ contained in the columns corresponding to $C$. We sum the elements of each representative incidence vector over the component $C$ and place it in the first position, with $0$'s in all other position of $C$.
%\begin{center}
%\begin{minipage}{0.4\linewidth}
\begin{multline*} \tiny
 \left(
    \begin{array}{c|c|c|c}
    \sum_{i=1}^{l_1 - 1} v_i^{m_1} \; 0 \cdots  0 & v_{l_1}^{m_1} \cdots v_{l_2-1}^{m_1} & \cdots & v_{l_{c-1}}^{m_1} \cdots v_n^{m_1} \\
    \sum_{i=1}^{l_1 - 1} v_i^{m_2} \; 0 \cdots  0 & v_{l_1}^{m_2} \cdots v_{l_2-1}^{m_2} & \cdots & v_{l_{c-1}}^{m_2} \cdots v_n^{m_2} \\
    \vdots          &  \vdots         & \ddots &  \vdots            \\
    \sum_{i=1}^{l_1 - 1} v_i^{m_p} \; 0 \cdots 0 & v_{l_1}^{m_p} \cdots v_{l_2-1}^{m_n} & \cdots & v_{l_{c-1}}^{m_p} \cdots v_n^{m_n} \\
    \hline
    \cdots \pm 1 \cdots &   0 \cdots 0      & \ldots        & 0 \cdots 0          \\
    \vdots              &    \vdots         & \ddots        &   \vdots            \\
    \cdots \pm 1 \cdots &  0 \cdots 0       & \ldots        & 0 \cdots 0          \\
    \hline
       0 \cdots 0       & \cdots \pm 1 \cdots  & \ldots        & 0 \cdots 0          \\
        \vdots          & \vdots               & \ddots        &   \vdots            \\
      0 \cdots 0        & \cdots \pm 1 \cdots  & \ldots        & 0 \cdots 0          \\
    \hline
    \ldots              &   \ldots          & \ddots        & \ldots              \\
    \hline
    0 \cdots 0 &   0 \cdots 0      & \ldots        & \cdots \pm 1 \cdots          \\
      \vdots   &    \vdots         & \ddots        & \vdots                       \\
    0 \cdots 0 &  0 \cdots 0       & \ldots        & \cdots \pm 1 \cdots          \\
    \end{array}
\right) \rightarrow \\
%\end{minipage}
%\hskip 5pt
%\begin{minipage}{0.4\linewidth}
 \tiny \left(
    \begin{array}{c|c|c|c}
    \sum_{i=1}^{l_1 - 1} v_i^{m_1} \; 0 \cdots  0 & \sum_{i=l_1}^{l_2 - 1} v_i^{m_1} \; 0 \cdots  0 & \cdots & \sum_{i=l_{c-1}}^{n} v_i^{m_1} \; 0 \cdots  0  \\
    \sum_{i=1}^{l_1 - 1} v_i^{m_2} \; 0 \cdots  0 & \sum_{i=l_1}^{l_2 - 1} v_i^{m_2} \; 0 \cdots  0 & \cdots & \sum_{i=l_{c-1}}^{n} v_i^{m_2} \; 0 \cdots  0  \\
    \vdots          &  \vdots         & \ddots &  \vdots            \\
    \sum_{i=1}^{l_1 - 1} v_i^{m_p} \; 0 \cdots 0 & \sum_{i=l_1}^{l_2 - 1} v_i^{m_n} \; 0 \cdots  0  & \cdots & \sum_{i=l_{c-1}}^{n} v_i^{m_n} \; 0 \cdots  0 \\
    \hline
    \cdots \pm 1 \cdots &   0 \cdots 0      & \ldots        & 0 \cdots 0          \\
    \vdots              &    \vdots         & \ddots        &   \vdots            \\
    \cdots \pm 1 \cdots &  0 \cdots 0       & \ldots        & 0 \cdots 0          \\
    \hline
       0 \cdots 0       & \cdots \pm 1 \cdots  & \ldots        & 0 \cdots 0          \\
        \vdots          & \vdots               & \ddots        &   \vdots            \\
      0 \cdots 0        & \cdots \pm 1 \cdots  & \ldots        & 0 \cdots 0          \\
    \hline
    \ldots              &   \ldots          & \ddots        & \ldots              \\
    \hline
    0 \cdots 0 &   0 \cdots 0      & \ldots        & \cdots \pm 1 \cdots          \\
      \vdots   &    \vdots         & \ddots        & \vdots                       \\
    0 \cdots 0 &  0 \cdots 0       & \ldots        & \cdots \pm 1 \cdots          \\
    \end{array}
\right) 
\end{multline*}
%\end{minipage}
%\end{center}
For each component $C$ of $g(X)$ there exists a column corresponding to a leaf of the spanning tree determined by the $\pm 1$ vectors in $C$, rooted at the first column. This column contains all $0$'s and one $1$. Performing cofactor expansion down this column eliminates a row and column. This can be repeated until all $\pm 1$-vectors of $C$ are eliminated, preserving the summed parts of the incidence vectors. Repeating for all connected components, all $\pm 1$-vectors can be eliminated while preserving the determinant of $X$ (under a sign change). We get the desired form and since each operation preserves the matrix being square we see that $c=p$.
\begin{equation*} \tiny 
     \left( \begin{array}{cccc} 
    \sum_{i=1}^{l_1-1}  v^{m_1}_i & \sum_{i=l_1}^{l_2-1}  v^{m_1}_i & \cdots & \sum_{i=l_{c-1}}^{n}  v^{m_1}_i \\
    \sum_{i=1}^{l_1-1}  v^{m_2}_i & \sum_{i=l_1}^{l_2-1}  v^{m_2}_i & \cdots & \sum_{i=l_{c-1}}^{n}  v^{m_2}_i \\
    \vdots & \vdots & \ddots & \vdots \\
    \sum_{i=1}^{l_1-1}  v^{m_c}_i & \sum_{i=l_1}^{l_2-1}  v^{m_c}_i & \cdots & \sum_{i=l_{c-1}}^{n}  v^{m_c}_i \\
    \end{array} \right) \in \R^{c \times c}
\end{equation*}

\end{proof}

\begin{corollary} \label{cor:connuni}
Let $\M$ be a connected matroid of rank $r$ on $n$ elements and $X \subset \Inc(\B_\M)$ be $n$ linearly independent incidence vectors of $\B_\M$. If $G(X)$ is connected then $|det(X)| = r$. 
\end{corollary}

\begin{example}
We offer an example of \autoref{thm:genconn} in action. Consider the uniform matroid $\U^{2,6}$ and the following collection of incidence vectors of bases in row form.\\
Partition into components of $G(X)$ and $g(X)$.
\begin{equation*} \tiny
\left| \begin{array}{cccccc}
1 & 1 & 0 & 0 & 1 & 0\\ 
1 & 1 & 0 & 0 & 0 & 1\\
1 & 0 & 1 & 1 & 0 & 0\\
0 & 1 & 1 & 1 & 0 & 0\\
0 & 0 & 1 & 0 & 1 & 1\\
0 & 0 & 0 & 1 & 1 & 1\\
\end{array}\right| \Longrightarrow 
\left| \begin{array}{cc|cc|cc}
1 & 1 & 0 & 0 & 1 & 0\\ 
1 & 1 & 0 & 0 & 0 & 1\\
\hline
1 & 0 & 1 & 1 & 0 & 0\\
0 & 1 & 1 & 1 & 0 & 0\\
\hline
0 & 0 & 1 & 0 & 1 & 1\\
0 & 0 & 0 & 1 & 1 & 1\\
\end{array}\right| %\quad \text{Partition into component of } G(X) \text{ and } g(X)
\end{equation*}
Find rooted spanning trees of each component of $G(X)$ and take differences to get $\pm 1$-vectors and the representative vectors.
\begin{equation*}  \tiny
\left| \begin{array}{cc|cc|cc}
1 & 1 & 0 & 0 & 1 & 0\\ 
0 & 0 & 0 & 0 &-1 & 1\\
\hline
1 & 0 & 1 & 1 & 0 & 0\\
-1& 1 & 0 & 0 & 0 & 0\\
\hline
0 & 0 & 1 & 0 & 1 & 1\\
0 & 0 &-1 & 1 & 0 & 0\\
\end{array}\right| \Longrightarrow
\left| \begin{array}{cc|cc|cc}
2 & 0 & 0 & 0 & 1 & 0\\ 
0 & 0 & 0 & 0 &-1 & 1\\
\hline
1 & 0 & 2 & 0 & 0 & 0\\
-1& 1 & 0 & 0 & 0 & 0\\
\hline
0 & 0 & 1 & 0 & 2 & 0\\
0 & 0 &-1 & 1 & 0 & 0\\
\end{array}\right| %\quad \text{Finding spanning rooted spanning tree and take differences}
\end{equation*}
For each component of $g(X)$ eliminate use $\pm 1$-vectors to eliminate rows and columns using cofactor expansion. 
\begin{equation*}  \tiny
\left| \begin{array}{c|c|cccc}
2 & 0 & 0 & 0 & 1 & 0\\ 
0 & 0 & 0 & 0 &-1 & 1\\
1 & 0 & 2 & 0 & 0 & 0\\
\hline
-1& 1 & 0 & 0 & 0 & 0\\
\hline
0 & 0 & 1 & 0 & 2 & 0\\
0 & 0 &-1 & 1 & 0 & 0\\
\end{array}\right| \Longrightarrow
\left| \begin{array}{ccccc}
2 & 0 & 0 & 1 & 0\\ 
0 & 0 & 0 &-1 & 1\\
1 & 2 & 0 & 0 & 0\\
0 & 1 & 0 & 2 & 0\\
0 &-1 & 1 & 0 & 0\\
\end{array}\right| %\quad \text{For each component of } g(X) \text{ eliminate } \pm 1 \text{ rows}
\end{equation*}
\begin{equation*}  \tiny
\left| \begin{array}{cc|c|cc}
2 & 0 & 0 & 1 & 0\\ 
0 & 0 & 0 &-1 & 1\\
1 & 2 & 0 & 0 & 0\\
0 & 1 & 0 & 2 & 0\\
\hline
0 &-1 & 1 & 0 & 0\\
\hline
\end{array}\right| \Longrightarrow
\left| \begin{array}{cccc}
2 & 0 & 1 & 0\\ 
0 & 0 &-1 & 1\\
1 & 2 & 0 & 0\\
0 & 1 & 2 & 0\\
\end{array}\right| \quad 
\end{equation*}
\begin{equation*}  \tiny
\left| \begin{array}{ccc|c|}
2 & 0 & 1 & 0\\ 
\hline
0 & 0 &-1 & 1\\
\hline
1 & 2 & 0 & 0\\
0 & 1 & 2 & 0\\
\end{array}\right| \Longrightarrow
\left| \begin{array}{ccc}
2 & 0 & 1\\ 
1 & 2 & 0\\
0 & 1 & 2\\
\end{array}\right|
\end{equation*}

\end{example}

We explore more connections between the graphs $G(X)$ and $g(X)$. The following is a alternate proof of part of \autoref{thm:genconn}
\begin{lemma}
Let $\M$ be a connected matroid of rank $r$ on $n$ elements and $X \subset \Inc(\B_\M)$ be $n$ linearly independent incidence vectors of $\B_\M$. If $G(X)$ is connected then $g(X)$ is connected.
\end{lemma}
\begin{proof}
Consider the matrix $X$ and, as in the proof of \autoref{thm:genconn}, we pick a row as a root of the spanning tree of $G(X)$ and perform row operations so that the resulting matrix is our root and the rest are $\pm 1$-vectors. If $g(X)$ is not connected then let $S$ be one connected component and $T=[n] \backslash S$. Also, let $s = |S|$ and $t = |T|$. Then the following non-zero vector is orthogonal to all the elements of $X$, contradicting its full-dimensionality.
\begin{equation*}
t\sum_{i\in S} \ve e_i - s\sum_{j\in T} e_j
\end{equation*}
\end{proof}
More generally we can show the following:
\begin{lemma}
Let $\M$ be a connected matroid of rank $r$ on $n$ elements and $X \subset \Inc(\B_\M)$ be $n$ incidence vectors of $\M$. If $G(X)$ is connected, then $\rank(X) = n+1 - \#$connected components of $g(X)$.
\end{lemma}
\begin{proof}
We will show that the space of functionals constant on $X$ has dimension equal to the number $k$ of connected components of $g(X)$ minus one. Considering $X$ as a matrix with incidence row vectors we first observe that all functionals in the kernel of $X$ must give the same value to all entries ($\ve e_i$'s) in the same connected component of $g(X)$. That is, they are characterized by $k$ parameters. Thus, the dimension of the space of functionals constant on $X$ is at most $k-1$.

Consider the matrix $X$ and, as in the proof of \autoref{thm:genconn}, we pick a row as a root of the spanning tree of $G(X)$ and perform row operations so that the resulting matrix is our root root and the rest are $\pm 1$-vectors. Omitting the root row of $X$, the orthogonal complement of the remaining rows of $X$ is exactly the $k$-dimensional linear space described above. Reinserting the root row decreases the dimension of the orthogonal compliment by at most one.
\end{proof}
\begin{lemma}
Let $\M$ be a connected matroid of rank $r$ on $n$ elements and $X \subset \Inc(\B_\M)$ be $n$ incidence vectors of $\B_\M$. If $g(X)$ is connected, $\rank(X) = n$. Moreover $G(X)$ is connected and thus $|\det(X)| = r$.
\end{lemma}
\begin{proof}
Consider the matrix $X$ and, as in the proof of \autoref{thm:genconn}, we pick representative vectors $\ve v^{m_1},\ldots,\ve v^{m_p}$ for each of the components of $G(X)$. Next perform row operations so that the resulting matrix consists of the representative vectors and $\pm 1$-vectors. Moreover since $g(X)$ is connected these $\pm 1$-vectors must be a spanning tree of an orientation of $g(X)$, hence there must be $(n-1)$ $\pm 1$-vectors. This implies that we have only one representative vector and thus $G(X)$ is connected. We also see that these $\pm 1$-vectors must be linearly independent. It follows that $X$ must be linearly independent. Then $|\det(X)| =r $ follows from \autoref{cor:connuni}.
\end{proof}
In order to prove a unimodular covering exists we propose another conjecture which implies \autoref{conj:unicovering}:
\begin{conjecture} \label{conj:concovering}
Let $\M$ be a matroid of rank $r$ on $n$ elements. For every $x \in \Po_\M$ there exists $n$ linearly independent vectors $\ve v_1,\ldots,\ve v_n \in \Inc(\B_\M)$ such that $x \in \conv(\ve v_1,\ldots, \ve v_n)$ and $G(X)$ connected. Hence, $\conv(\ve v_1,\ldots, v_n)$ is unimodular.
\end{conjecture}

In support of \autoref{conj:white}, \autoref{conj:unicovering} and \autoref{conj:concovering} we now present some experimental results. We obtained bases description for all matroids of nine elements or less from \cite{Mayhew2007Matroids-with-n}. The number of matroids by rank and elements less than or equal to nine is given by the table:
\begin{equation*}
\begin{array}{l|rrrrrrrrrr}
%\begin{array}{l|r|r|r|r|r|r|r|r|r|r}
\toprule
%\text{\begin{minipage}{2cm}Elements\\ Rank\end{minipage}} & 0 & 1 & 2 & 3 & 4 & 5 & 6 & 7 & 8 & 9\\
\text{Elements} & 0 & 1 & 2 & 3 & 4 & 5 & 6 & 7 & 8 & 9\\
\text{Rank} \\
\hline
0 	& 1 	& 1 	& 1 	& 1 	& 1 	& 1 	& 1 	& 1 	& 1 	& 1        \\
1 	&   	& 1 	& 2 	& 3 	& 4 	& 5 	& 6 	& 7 	& 8 	& 9       \\
2 	&   	&   	& 1 	& 3 	& 7 	& 13 	& 23 	& 37 	& 58 	& 87       \\
3 	&   	&   	&   	& 1 	& 4 	& 13 	& 38 	& 108	& 325	& 1275       \\
4 	&   	&   	&   	&   	& 1 	& 5 	& 23 	& 108	& 940	& 190214       \\
5 	&   	&   	&   	&   	&   	& 1 	& 6 	& 37 	& 325	& 190214       \\
6 	&   	&   	&   	&   	&   	&   	& 1 	& 7 	& 58 	& 1275       \\
7 	&   	&   	&   	&   	&   	&   	&   	& 1 	& 8 	& 87       \\
8 	&   	&   	&   	&   	&   	&   	&   	&   	& 1 	& 9       \\
9 	&   	&   	&   	&   	&   	&   	&   	&   	&   	& 1       \\
\hline
\text{Total Matroids} & 1 &	2 &	4 &	8 &	17 &	38 & 	98 &	306 & 	1724 &	383172 \\
\bottomrule
\end{array}
\end{equation*}

%There are $473$ matroids of rank seven or less, $1724$ matroids of rank eight and $383172$ matroids of rank nine. 
With this data in hand we used TOPCOM \cite{Rambau:TOPCOM-ICMS:2002} to compute triangulations of all the connected matroid polytopes. From these computational experiments we can justify the following theorem:
\begin{theorem} \label{thm:uniexper}
All $1317$ matroid polytopes of connected matroids with eight or less elements have a unimodular triangulation. Moreover, the placing triangulation, determined by the order of bases as presented in \cite{Mayhew2007Matroids-with-n}, is a unimodular triangulation.
\end{theorem}
Whether the ordering is important for the placing triangulation to be unimodular for these experiments is unclear.  As a note, there certainly are non-unimodular triangulations of matroid polytopes. As a natural follow up experiment, we attempted to find as many matroids in the test set that were unimodular as well as triangulated by simplices that had connected graphs $G(X)$.
\begin{theorem} \label{thm:coverexper}
All $48$ matroid polytopes of connected matroids with six or less elements have a triangulation where for each simplex $X$, $G(X)$ is connected. Hence there exists a unimodular triangulation of each matroid polytope. 
\end{theorem}
For matroids with seven or more elements, it becomes especially time consuming to find a triangulation with $G$-connected simplices. We were able to find $36$ matroids of rank seven that had $G$-connected simplices. All matroids and triangulations for \autoref{thm:uniexper} and \autoref{thm:coverexper} can be found at \cite{HawsMatroid-Polytop}.

%\davecomment{Here I will place some computational results about unimodular triangulations discovered using TOPCOM and the Royle matroids.}

\section[Two Dimensional Faces]{Two Dimensional Faces of Matroid Polytopes}
%\davecomment{This section likely to be removed. It is an attempt to triangulate a matroid polytope using the $2$-face structure and an idea of sheering off pieces of the polytope and triangulating each piece separately}
In the pursuit of the proof of \autoref{conj:white}, \autoref{conj:unicovering}, and \autoref{conj:concovering} related to matroid polytope triangulations, coverings and decompositions, we studied the two dimensional faces of $\Po_\M$ in order to understand how such subdivisions interact with the low dimensional faces. We obtained a few results.

\begin{lemma}[See Theorem 4.8 in \cite{Alexandre-V.-Borovik2007Coxeter-matroid}] \label{lem:2faces}
    Let $\M$ be a matroid and $\Po_\M$ be its matroid polytope. The $2$-dimensional faces of $\Po_\M$ are equilateral triangles or squares.
\end{lemma}

%We show that every two dimensional square face of a matroid polytope is adjacent as in \autoref{fig:2adj}.
\begin{figure}[!htb] 
\ifpdf
    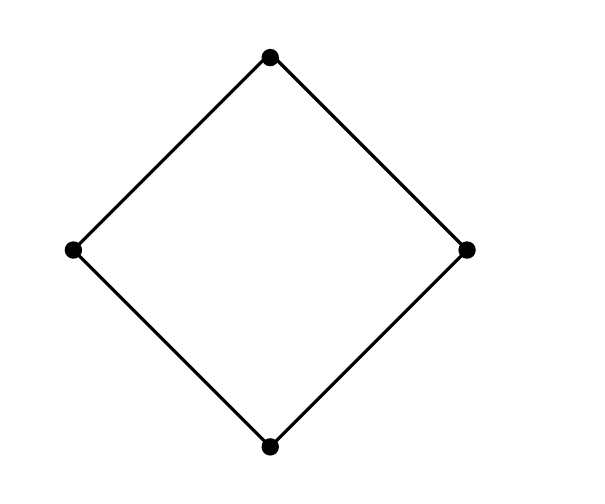
    \else
    \input{2adj.pstex_t}
    \fi
\caption{Square $2$-face adjacency structure. $\ve e_s \neq \ve e_m$, $\ve e_t \neq \ve e_l$, $\ve e_s \neq \ve e_t$, $\ve e_m \neq \ve e_l$.} \label{fig:2adj} \label{fig:2dimface}
\end{figure}

\begin{lemma} \label{lem:2face}
    Let $\M$ be a matroid of rank $r$ on $n$ elements and $F$ be a $2$-dimensional square face of $\Po_\M$ with vertices $\ve w_1, \ve w_2, \ve w_3, \ve w_4$ such that $(\ve w_1, \ve w_2)$, $(\ve w_2, \ve w_4)$, $(\ve w_1, \ve w_3)$, and $(\ve w_3, \ve w_4)$ are the edges. There exists standard unit vectors in $\ve e_s, \ve e_t, \ve e_m, \ve e_l \in \R^n$ such that $\ve w_2 = \ve w_1 + \ve e_s - \ve e_t$, $\ve w_3 = \ve w_1 + \ve e_m - \ve e_l$, and $\ve w_4 = \ve w_1 + \ve e_s - \ve e_t + \ve e_m - \ve e_l$. Moreover, and $\ve e_s \neq \ve e_m$, $\ve e_t \neq \ve e_l$, $\ve e_s \neq \ve e_t$ and $\ve e_m \neq \ve e_l$.

%    If $\ve v$ is a vertex of $F$ then $\ve v+(\ve e_s-\ve e_t), \ve v+(\ve e_m - \ve e_l), \ve v+(\ve e_s-\ve e_t) + (\ve e_m - \ve e_l)$ are the other vertices of $F$ where $\ve e_s, \ve e_t, \ve e_m, \ve e_l$ are standard unit vectors in $\R^n$ and $\ve e_s \neq \ve e_m$, $\ve e_t \neq \ve e_l$, $\ve e_s \neq \ve e_t$ and $\ve e_m \neq \ve e_l$.
\end{lemma}
\begin{proof}
%      Let $\M$ be a matroid and $F$ a $2$-dimensional face of $\Po_\M$ and $\ve w_1, \ve w_2, \ve w_3, \ve w_4$ vertices of $F$, where $\ve w_1$ and $\ve w_4$ are adjacent to $\ve w_2$ and $\ve w_3$. By \autoref{lem:adj} the vertices of $F$ differ from each other by the addition and subtraction of standard unit vectors. See \autoref{fig:2adjproof}. Since $F$ is a $2$ dimensional face, then $\ve e_s \neq \ve e_m$, otherwise $\ve w_2$ and $\ve w_3$ are adjacent by \autoref{lem:adj}, contradicting $F$ being a $2$-face. Similarly we can conclude that $\ve e_t \neq \ve e_l$, $\ve e_s \neq \ve e_q$, $\ve e_t \neq \ve e_k$, $\ve e_m \neq \ve e_o$, $\ve e_l \neq \ve e_p$, $\ve e_q \neq \ve e_o$, and $\ve e_k \neq \ve e_p$. But by vector addition, $\ve w_4 = \ve w_1 + (\ve e_s - \ve e_t) + (\ve e_q - \ve e_k) = \ve w_1 + (\ve e_m - \ve e_l) + (\ve e_o - \ve e_p)$, and therefore $\ve e_q = \ve e_m$, $\ve e_o = \ve e_s$, $\ve e_k = \ve e_l$, and $\ve e_t = \ve e_p$. See \autoref{2adj}.
%%%   
%%%       \begin{figure}[!htb] 
%%%       \inputfig{2adjproof}
%%%       \caption{$2$-face adjacency structure.} \label{fig:2adjproof}
%%%       \end{figure}
%%%   

%  {\it Alternate proof:}
Let $\M$ be a matroid and $F$ a $2$-face of $\Po_\M$ with $\ve w_1, \ve w_2, \ve w_3, \ve w_4$ vertices of $F$ such that $(\ve w_1, \ve w_2)$, $(\ve w_2, \ve w_4)$, $(\ve w_1, \ve w_3)$, and $(\ve w_3, \ve w_4)$ are the only adjacencies. Due to \autoref{lem:adj} we know if $\ve w_2 = \ve w_1+ \ve e_s - \ve e_t$ and $\ve w_4 = \ve w_1 + \ve e_s - \ve e_t + \ve e_m - \ve e_l$, then $\ve w_3 = \ve w_1 + \ve e_m - \ve e_l$. Otherwise an additional adjacency (edge) is created among $\ve w_1, \ve w_2, \ve w_3, \ve w_4$, a contradiction to $F$ a $2$-face. E.g., if $\ve w_3 = \ve w_1 + \ve e_s - \ve e_l$ then $\ve w_3$ and $\ve w_2$ are adjacent, a contradiction. Then by vector addition $\ve w_4 = \ve w_3 + \ve e_m - \ve e_l$. Additionally, $\ve e_s \neq \ve e_m$, $\ve e_t \neq \ve e_l$, $\ve e_s \neq \ve e_t$ and $\ve e_m \neq \ve e_l$. Otherwise it would create an additional adjacency, hence edge. See \autoref{fig:2dimface}.

%Then $\ve w_1$ is adjacent to $\ve w_2$ by the addition of $\ve e_s - \ve e_t$. Similarly, $\ve w_2$ is adjacent to $\ve w_4$ by the addition of $\ve e_m - \ve e_t$. Let $\ve w_3$ be the fourth vertex of $F$ adjacent to $\ve w_1$ and $\ve w_4$. Then by vector addition and \autoref{lem:adj}, $\ve w_1$ is adjacent to $\ve w_3$ by the addition of $\ve e_m - e_l$, otherwise $\ve w_2$ and $\ve w_3$ are adjacent or $\ve w_1$ and $\ve w_4$ are adjacent by \autoref{lem:adj}.

\end{proof}

\begin{figure}[!htb] 
\ifpdf
    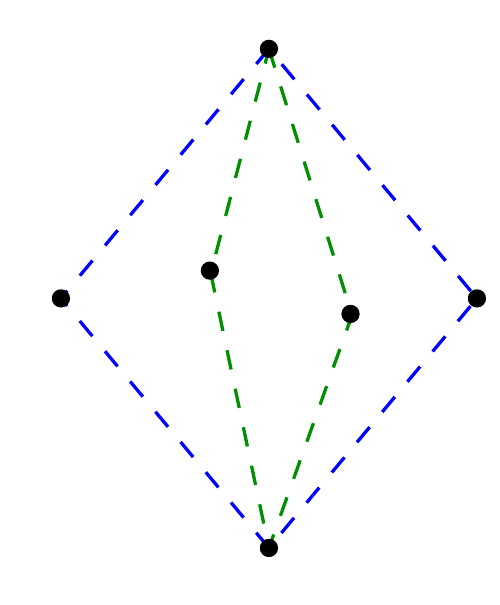
    \else
    \input{2adj6vert.pstex_t}
    \fi
\caption{Six adjacent vertices restrict any four from being a $2$-face.} \label{fig:2adj6vert}
\end{figure}

\begin{lemma} \label{lem:squareface}
    Let $\M$ be a matroid of rank $r$ on $n$ elements, $\ve w_1, \ve w_2, \ve w_3, \ve w_4 \in \Po_\M$. Let $\ve w_1$ and $\ve w_4$ be adjacent to $\ve w_2$ and $\ve w_3$ in $\Po_\M$ where $\ve w_2 = \ve w_1 + (\ve e_s - \ve e_t)$, $\ve w_3 = \ve w_1 + (\ve e_m - \ve e_l)$ and $\ve w_4 = \ve w_1 + (\ve e_s - \ve e_t) + (\ve e_m - \ve e_l)$, $\ve e_s \neq \ve e_m$, $\ve e_t \neq \ve e_l$, $\ve e_s \neq \ve e_t$ and $\ve e_m \neq \ve e_l$. Then $\ve w_1, \ve w_2, \ve w_3$ and $\ve w_4$ is a square $2$-face if and only if $\ve w_5 := \ve w_1 + (\ve e_s - \ve e_l)$ or $\ve w_6 := \ve w_1 + (\ve e_m - \ve e_t)$ are not vertices of $\Po_\M$.
\end{lemma}
\begin{proof}
    Assume that both $\ve w_5$ and $\ve w_6$ are vertices of $\Po_\M$ and $F$ is a square $2$-face with vertices $\ve w_1, \ve w_2, \ve w_3$ and $\ve w_4$.
    Let 
    \begin{equation} \label{eq:suphyp}
        f(\ve x) := a_1x_1 + a_2x_2 + \cdots + a_nx_n = c
    \end{equation}
    be a supporting hyperplane of $F$, i.e. $\ve w_1, \ve w_2, \ve w_3, \ve w_4$ satisfy \eqref{eq:suphyp}. Moreover, for all $z \in \Po_\M \setminus F$, $f(z) < c$.  Since Equation \eqref{eq:suphyp} is a supporting hyperplane and $\ve w_5, \ve w_6 \notin F$, $f(\ve w_5) < c$ and $f(\ve w_6) < c$. But $2c = f(\ve w_2 + \ve w_3) = f(\ve w_5 + \ve w_6) < 2c$, a contradiction.\\

    Assume, without loss of generality, that $\ve w_5$ is a vertex of $\Po_\M$, while $\ve w_6$ is not. We will construct a supporting hyperplane of $F$. Let
    \begin{equation} \label{eq:suphyp}
        f(\ve x) := a_1x_1 + a_2x_2 + \cdots + a_nx_n. 
    \end{equation}
    Let $a_s = a_t = 1$, $a_m = a_l = 2$, $a_i = 3$ if $i \in \supp(\ve w_1) \setminus \{s,t,m,l\}$ and $a_j = -1$ if $i \notin \supp(\ve w_1) \cup \{s,t,m,l\}$. We see that $f(\ve w_1) = f(\ve w_2) = f(\ve w_3) = f(\ve w_4) = 3(r-2) + 3$, yet $f(\ve w_5) = 3(r-1) + 2$ and for all vertices $\ve z \in \Po_\M \setminus F$, $f(\ve z) < 3(r-1) + 3$. Hence $f$ is a supporting hyperplane and $F$ is a square $2$-face. If $\ve w_5$ is also not a vertex, then $f$ is still a supporting hyperplane of $F$.
\end{proof}

\begin{corollary} \label{lem:uni2face}
    Let $\M$ be a uniform matroid of rank $r$ on $n$ elements. Then the $2$-dimensional faces of $\Po_\M$ are all equilateral triangles. 
\end{corollary}
\begin{proof}
    This follows from \autoref{lem:squareface} and that if $\ve w_2 = \ve w_1 + (\ve e_s - \ve e_t)$ and $\ve w_3 = \ve w_1 + (\ve e_m - \ve e_l)$ are vertices of $\Po_\M$ then so are $\ve w_5 := v + (\ve e_s - \ve e_l)$ and $\ve w_6 := v + (\ve e_m - \ve e_t)$ since $\M$ is a uniform matroid.
\end{proof}

\begin{proposition}
    Let $\M$ be a connected matroid of rank $r$ on $n$ elements and $v$ be a vertex of $\Po_\M$. Then all adjacent vertices of $v$ lie on a $(n-2)$-hyperplane
\end{proposition}
\begin{proof}
    Without loss of generality assume $\ve v = (\underbrace{1, \ldots, 1}_{r}, \underbrace{0, \ldots, 0}_{n-r})^\top$. Then by \autoref{lem:adj} all vertices $\ve w$ adjacent to $\ve v$ lie in the hyperplane $x_1 + \cdots + x_r = r-1$ since $\ve v - \ve w = \ve e_s - \ve e_l$ where $s \in \{ 1, \ldots, r\}$ and $ l \notin \{ 1, \ldots, r\}$.
\end{proof}

In light of \autoref{thm:genconn} we note that any triangulation of $\{\ve v, \Adj(\ve v)\}$ is unimodular where $\ve v \in \Inc(\B_\M)$. This follows from the fact that all adjacent vertices of $\ve v$ are in an $n-2$ hyperplane, and thus any full-dimensional simplex $X$ must contain $\ve v$ which implies $G(X)$ is connected.

   \chapter[% 
     Applications to Matroid Polytopes
   ]{% 
     Applications to Optimization Through the Structure of Matroid Polytopes
   }%
   \label{ch:4thChapterLabel}
   \input{Chapter4.tex}
   
   %\appendix

   %\chapter[%
   %   Short Title of Appendix A
   %]{%
   %   Long Title of Appendix A
   %}%
   %\label{ch:AppendixALabel}
   %\input{AppendixAFileName.tex}
       
   \backmatter
   
   \bibliographystyle{amsalpha-fi-arxlast}
   \bibliography{DissertationBibliography}

\newcommand{\etalchar}[1]{$^{#1}$}
\providecommand{\bysame}{\leavevmode\hbox to3em{\hrulefill}\thinspace}
\providecommand{\MR}{\relax\ifhmode\unskip\space\fi MR }
\providecommand{\MRhref}[2]{%
  \href{http://www.ams.org/mathscinet-getitem?mr=#1}{#2}
}
\providecommand{\href}[2]{#2}
\begin{thebibliography}{BLMA{\etalchar{+}}08}

\bibitem[ABB{\etalchar{+}}99]{laug}
E.~Anderson, Z.~Bai, C.~Bischof, S.~Blackford, J.~Demmel, J.~Dongarra,
  J.~Du~Croz, A.~Greenbaum, S.~Hammarling, A.~McKenney, and D.~Sorensen,
  \emph{{LAPACK} users' guide}, third ed., Society for Industrial and Applied
  Mathematics, Philadelphia, PA, 1999.

\bibitem[Bar94]{bar}
A.~I. Barvinok, \emph{Polynomial time algorithm for counting integral points in
  polyhedra when the dimension is fixed}, Mathematics of Operations Research
  \textbf{19} (1994), 769--779.

\bibitem[Bar06]{barvinok-2006-ehrhart-quasipolynomial}
\bysame, \emph{Computing the {E}hrhart quasi-polynomial of a rational simplex},
  Math. Comp. \textbf{75} (2006), no.~255, 1449--1466 (electronic).

\bibitem[BGW07]{Alexandre-V.-Borovik2007Coxeter-matroid}
A.~V. Borovik, I.~M. Gelfand, and N.~White, \emph{Coxeter matroid polytopes},
  Annals of Combinatorics \textbf{1} (2007), no.~1, 123--134.

\bibitem[BHS09]{beck-haase-sottile:theorema}
M.~Beck, C.~Haase, and F.~Sottile, \emph{{Formulas of Brion, Lawrence, and
  Varchenko on rational generating functions for cones}}, The Mathematical
  Intelligencer \textbf{31} (2009), no.~1, 9--17.

\bibitem[BJR06]{billera-2006}
L.~J. Billera, N.~Jia, and V.~Reiner, \emph{A quasisymmetric function for
  matroids}, eprint arXiv:math/0606646, 2006.

\bibitem[Bla08]{Blasiak2005The-Toric-Ideal}
J.~Blasiak, \emph{The toric ideal of a graphic matroid is generated by
  quadrics}, Combinatorica \textbf{28} (2008), no.~3, 283--297.

\bibitem[BLMA{\etalchar{+}}08]{berstein-2008-22}
Y.~Berstein, J.~Lee, H.~Maruri-Aguilar, S.~Onn, E.~Riccomagno, R.~Weismantel,
  and H.~Wynn, \emph{Nonlinear matroid optimization and experimental design},
  SIAM Journal on Discrete Mathematics \textbf{22} (2008), 901.

\bibitem[BLOW08]{berstein-2008}
Y.~Berstein, J.~Lee, S.~Onn, and R.~Weismantel, \emph{Nonlinear optimization
  for matroid intersection and extensions}, 2008.

\bibitem[BO08]{berstein-2008-5}
Y.~Berstein and S.~Onn, \emph{Nonlinear bipartite matching}, Discrete
  Optimization \textbf{5} (2008), 53.

\bibitem[BP99]{barvinok:99}
A.~I. Barvinok and J.~E. Pommersheim, \emph{An algorithmic theory of lattice
  points in polyhedra}, New Perspectives in Algebraic Combinatorics (L.~J.
  Billera, A.~Bj\"orner, C.~Greene, R.~E. Simion, and R.~P. Stanley, eds.),
  Math. Sci. Res. Inst. Publ., vol.~38, Cambridge Univ. Press, Cambridge, 1999,
  pp.~91--147.

\bibitem[BR07]{Beck2007Computing}
M.~Beck and S.~Robins, \emph{Computing the continuous discretely: Integer-point
  enumeration in polyhedra}, Springer, 2007.

\bibitem[Bri88]{Brion88}
M.~Brion, \emph{Points entiers dans les poly{\'e}dres convexes}, Ann. Sci.
  {\'E}cole Norm. Sup. \textbf{21} (1988), no.~4, 653--663.

\bibitem[BS07a]{Barvinok2007Random-weightin}
A.~Barvinok and A.~Samorodnitsky, \emph{Random weighting, asymptotic counting,
  and inverse isoperimetry}, Israel Journal of Mathematics \textbf{158} (2007),
  159--191.

\bibitem[BS07b]{beck-sottile:irrational}
M.~Beck and F.~Sottile, \emph{Irrational proofs for three theorems of
  {S}tanley}, European Journal of Combinatorics \textbf{28} (2007), no.~1,
  403--409, \mbox{math.CO/0501359}.

\bibitem[BW91]{brightwellwinkler91}
G.~Brightwell and P.~Winkler, \emph{Counting linear extensions}, Order
  \textbf{8} (1991), no.~3, 225--242.

\bibitem[BW03]{barvinok-woods-2003}
A.~I. Barvinok and K.~Woods, \emph{Short rational generating functions for
  lattice point problems}, Journal of the AMS \textbf{16} (2003), no.~4,
  957--979.

\bibitem[CS86]{Cook1986An-integer-anal}
W.~Cook and J.~F.~A. Schrijver, \emph{An integer analogue of caratheodory's
  theorem}, Journal of Combinatorial Theory \textbf{510} (1986), 179--185.

\bibitem[DF88]{dyerfrieze88}
M.~E. Dyer and A.~M. Frieze, \emph{On the complexity of computing the volume of
  a polyhedron}, SIAM J. Comput. \textbf{17} (1988), no.~5, 967--974.

\bibitem[DHLO09]{MOCHA}
J.~{De Loera}, D.~C. Haws, J.~Lee, and A.~O'Hair, \emph{{M}atroids
  {O}ptimization {C}ombinatorics {H}euristics and {A}lgorithms}, Available from
  URL {\url{http://math.ucdavis.edu/~haws/MOCHA/}}, 2009.

\bibitem[DLHH{\etalchar{+}}04]{latte2}
J.~A. De~Loera, D.~Haws, R.~Hemmecke, P.~Huggins, B.~Sturmfels, and R.~Yoshida,
  \emph{Short rational functions for toric algebra and applications}, Journal
  of Symbolic Computation \textbf{38} (2004), no.~2, 959--973.

\bibitem[DLHH{\etalchar{+}}05]{lattesoft}
J.~A. De~Loera, D.~C. Haws, R.~Hemmecke, P.~Huggins, J.~Tauzer, and R.~Yoshida,
  \emph{Software and user's guide for latte v.1.1}, 2005.

\bibitem[DLHK08]{deloera-haws-koeppe:ehrhart-matroid}
J.~A. De~Loera, D.~C. Haws, and M.~K{\"o}ppe, \emph{Ehrhart polynomials of
  matroid polytopes and polymatroids}, Discrete and Computational Geometry
  (2008).

\bibitem[DLHTY04]{latte1}
J.~A. De~Loera, R.~Hemmecke, J.~Tauzer, and R.~Yoshida, \emph{Effective lattice
  point counting in rational convex polytopes}, Journal of Symbolic Computation
  \textbf{38} (2004), no.~4, 1273--1302.

\bibitem[DLRS09]{De-Loera2006Triangulations}
J.~A. De~Loera, J.~Rambau, and F.~Santos, \emph{Triangulations: Applications,
  structures, algorithms}, Book manuscript, 2009.

\bibitem[DNH97]{Negri1997Gorenstein-Alge}
E.~De~Negri and T.~Hibi, \emph{Gorenstein algebras of {V}eronese type}, Journal
  of Algebra \textbf{193} (1997), no.~2, 629--639.

\bibitem[dPS03]{Pina2003Improved-bound-}
J.~C. de~Pina and J.~Soares, \emph{Improved bound for the carath\'eodory rank
  of the bases of a matroid}, Journal of Combinatorial Theory \textbf{B}
  (2003), no.~88, 323--327.

\bibitem[Edm03]{Edmonds2003Submodular-func}
J.~Edmonds, \emph{Submodular functions, matroids, and certain polyhedra},
  Combinatorial Optimization -- Eureka, You Shrink!: Papers Dedicated to Jack
  Edmonds. 5th International Workshop, Aussois, France, March 5--9, 2001,
  Revised Papers (M.~J\"unger, G.~Reinelt, and G.~Rinaldi, eds.), Lecture notes
  in computer science, vol. 2570, Springer-Verlag, Berlin, 2003, pp.~11--26.

\bibitem[EG00]{ehrgottsurvey}
M.~Ehrgott and X.~Gandibleux, \emph{A survey and annotated bibliography of
  multiobjective combinatorial optimization}, OR Spektrum \textbf{22} (2000),
  no.~4, 425--460.

\bibitem[Ehr96]{ehrgottmatroid}
M.~Ehrgott, \emph{On matroids with multiple objectives}, Optimization
  \textbf{38} (1996), no.~1, 73--84, Multicriteria optimization and decision
  theory (Holzhau, 1994).

\bibitem[Ele86]{elekes86}
G.~Elekes, \emph{A geometric inequality and the complexity of computing
  volume}, Discrete Comput. Geom. \textbf{1} (1986), no.~4, 289--292.

\bibitem[FH80]{Fries1980Minimum-aberrat}
A.~Fries and W.~Hunter, \emph{Minimum aberration $2^{k-p}$ designs}, Techno.
  (1980), no.~22, 601--608.

\bibitem[FS05]{Feichtner2004Matroid-polytop}
E.~M. Feichtner and B.~Sturmfels, \emph{Matroid polytopes, nested sets and
  {B}ergman fans}, Port. Math. \textbf{62} (2005), 437--468.

\bibitem[Fuk06]{cdd}
K.~Fukuda, \emph{cdd+, a {C++} implementation of the double description method
  of {Motzkin} et al.}, Available from URL
  {\url{http://www.ifor.math.ethz.ch/~fukuda/cdd_home/cdd.html}}, 2006.

\bibitem[Gea09]{GMP}
T.~Granlund and et~al., \emph{{GNU} multiple precision arithmetic library},
  2009, \url{http://gmplib.org/}.

\bibitem[GGMS87]{Gelfand1987Combinatorial-g}
I.~M. Gelfand, M.~Goresky, R.~D. MacPherson, and V.~V. Serganova,
  \emph{Combinatorial geometries, convex polyhedra, and {S}chubert cells}, Adv.
  Math. \textbf{63} (1987), 301--316.

\bibitem[GJ79]{Garey1979Computers-and-i}
M.~Garey and D.~Johnson, \emph{Computers and intractability: A guide to the
  theory of {NP}-completeness}, W.H. Freeman and Company, San Francisco, 1979.

\bibitem[GLLR79]{Graham1979Optimization-an}
R.~Graham, E.~Lawler, J.~Lenstra, and A.~{Rinnooy Kan}, \emph{Optimization and
  approximation in deterministic sequencing and scheduling: A survey}, Annals
  of Discrete Mathematics (1979), no.~5, 287--326.

\bibitem[GLM08]{Gunnels2008IBM-Research-Re}
J.~Gunnels, J.~Lee, and S.~Margulies, \emph{{IBM} research report}, Preprint
  (2008).

\bibitem[Glo86]{Glover1986Future-Paths-fo}
F.~Glover, \emph{Future paths for integer programming and links to artificial
  intelligence}, Computers and Operations Research \textbf{13} (1986),
  533--549.

\bibitem[GO97]{1997Handbook-of-dis}
J.~E. Goodman and J.~O'Rourke (eds.), \emph{Handbook of discrete and
  computational geometry}, CRC Press, Inc., Boca Raton, FL, USA, 1997.

\bibitem[Haw09]{HawsMatroid-Polytop}
D.~C. Haws, \emph{Matroid polytopes},
  \verb#http://math.ucdavis.edu/~haws/Matroids/#, 2009.

\bibitem[HH72]{Hotzmann1972On-the-tree-gra}
C.~Hotzmann and F.~Harary, \emph{On the tree graph of a matroid}, SIAM Journal
  of Applied Math \textbf{22} (1972), no.~2.

\bibitem[Hib92]{Hibi1992Algebraic-Combi}
T.~Hibi, \emph{Algebraic combinatorics on convex polytopes}, Carslaw
  Publications, Glebe, Australia, 1992.

\bibitem[Kat05]{Katzman:math0408038}
M.~Katzman, \emph{The {H}ilbert series of algebras of {V}eronese type},
  Communications in Algebra \textbf{33} (2005), 1141--1146,
  \mbox{math/0408038}.

\bibitem[KC02]{Knowles}
J.~D. Knowles and D.~W. Corne, \emph{Enumeration of pareto optimal
  multi-criteria spanning trees - a proof of the incorrectness of zhou and
  gen's proposed algorithm}, European Journal of Operational Research
  \textbf{143} (2002), no.~3, 543 -- 547.

\bibitem[Kha93]{khachiyan93}
L.~Khachiyan, \emph{Complexity of polytope volume computation}, New Trends in
  Discrete and Computational Geometry, Algorithms Combin., vol.~10, Springer,
  Berlin, 1993, pp.~91--101.

\bibitem[K{\"o}p07a]{latte-macchiato}
M.~K{\"o}ppe, \emph{{LattE macchiato}, version 1.2-mk-0.9, an improved version
  of {De Loera} et al.'s {LattE} program for counting integer points in
  polyhedra with variants of {Barvinok}'s algorithm}, Available from URL
  {\url{http://www.math.uni-magdeburg.de/~mkoeppe/latte/}}, 2007.

\bibitem[K{\"o}p07b]{koeppe:irrational-barvinok}
M.~K{\"o}ppe, \emph{A primal {B}arvinok algorithm based on irrational
  decompositions}, SIAM Journal on Discrete Mathematics \textbf{21} (2007),
  no.~1, 220--236, \mbox{math.CO/0603308}.

\bibitem[KV08]{koeppe-verdoolaege:parametric}
M.~K\"oppe and S.~Verdoolaege, \emph{Computing parametric rational generating
  functions with a primal {B}arvinok algorithm}, The Electronic Journal of
  Combinatorics \textbf{15} (2008), \#R16.

\bibitem[Law91a]{lawrence91}
J.~Lawrence, \emph{Polytope volume computation}, Math. Comp. \textbf{57}
  (1991), no.~195, 259--271.

\bibitem[Law91b]{lawrence91-2}
\bysame, \emph{Rational-function-valued valuations on polyhedra}, Discrete and
  Computational Geometry (New Brunswick, NJ, 1989/1990), DIMACS Ser. Discrete
  Math. Theoret. Comput. Sci., vol.~6, Amer. Math. Soc., Providence, RI, 1991,
  pp.~199--208.

\bibitem[Mat97]{matsui}
T.~Matsui, \emph{A flexible algorithm for generating all the spanning trees in
  undirected graphs}, Algorithmica \textbf{18} (1997), no.~4, 530--543.

\bibitem[May08]{Mayhew}
D.~Mayhew, \emph{Matroid complexity and nonsuccinct descriptions}, SIAM J.
  Discret. Math. \textbf{22} (2008), no.~2, 455--466.

\bibitem[MRar]{Mayhew2007Matroids-with-n}
D.~Mayhew and G.~Royle, \emph{Matroids with nine elements}, Journal of
  Combinatorial Theory ((To appear)).

\bibitem[Onn03]{onn-2003-17}
S.~Onn, \emph{Convex matroid optimization}, SIAM Journal on Discrete
  Mathematics \textbf{17} (2003), 249.

\bibitem[Oxl92]{Oxley1992Matroid-Theory}
J.~Oxley, \emph{Matroid theory}, Oxford University Press, New York, NY, USA,
  1992.

\bibitem[PWZ96]{Petkovsek1996AB}
M.~{Petkov\v sek}, H.~S. Wilf, and D.~Zeilberger, \emph{${A}={B}$}, AK Peters,
  Ltd., 1996.

\bibitem[Ram02]{Rambau:TOPCOM-ICMS:2002}
J.~Rambau, \emph{{TOPCOM}: Triangulations of point configurations and oriented
  matroids}, Mathematical Software---ICMS 2002 (A.~M. Cohen, X.-S. Gao, and
  N.~Takayama, eds.), World Scientific, 2002, pp.~330--340.

\bibitem[Sch86a]{schrijver}
A.~Schrijver, \emph{Theory of linear and integer programming},
  Wiley-Interscience, 1986.

\bibitem[Sch86b]{Schrijver1986Theory-of-linea}
A.~Schrijver, \emph{Theory of linear and integer programming}, John Wiley \&
  Sons, Inc., New York, NY, USA, 1986.

\bibitem[Sch03]{Schrijver2003Combinatorial-O}
A.~Schrijver, \emph{Combinatorial optimization: Polyhedra and efficiency},
  Springer, 2003.

\bibitem[Seb90]{Sebo1990Hilbert-bases-C}
A.~Seb\"o, \emph{{H}ilbert bases, {C}arath\'eodory's theorem and combinatorial
  optimization}, University of Waterlo Press, R. Kannan and W. Pullyblanks eds,
  1990.

\bibitem[Spe06]{SpeyerA-matroid-invar}
D.~E. Speyer, \emph{A matroid invariant via the {K}-theory of the
  {G}rassmannian}, eprint arXiv:math/0603551v1, 2006.

\bibitem[Sta96]{Stanley1996Combinatorics-a}
R.~P. Stanley, \emph{Combinatorics and commutative algebra: Second edition},
  2nd ed., Birkh\"auser, Boston, 1996.

\bibitem[Sta97]{Stanley1997Enumerative-Com}
R.~P. Stanley, \emph{Enumerative combinatorics}, vol.~1, Cambridge University
  Press, 1997.

\bibitem[Stu96]{Sturmfels1996Grobner-Bases-a}
B.~Sturmfels, \emph{Gr\"obner bases and convex polytopes}, University Lecture
  Series, vol.~8, American Mathematical Society, 1996.

\bibitem[Top84]{Topkis1984Adjacency-on-Po}
D.~M. Topkis, \emph{Adjacency on polymatroids}, Mathematical Programming
  \textbf{30} (1984), no.~2, 229--237.

\bibitem[VW08]{verdoolaege-woods-2005}
S.~Verdoolaege and K.~M. Woods, \emph{Counting with rational generating
  functions}, J. Symb. Comput. \textbf{43} (2008), no.~2, 75--91.

\bibitem[War85]{Warburton1985Worse-case-anal}
A.~Warburton, \emph{Worse case analysis of greedy and related heuristics for
  some min-max combinatorial optimization problems}, Mathematical Programming
  \textbf{33} (1985), no.~2, 234--241.

\bibitem[Wel76]{Welsh1976Matroid-Theory}
D.~Welsh, \emph{Matroid theory}, Academic Press, Inc., 1976.

\bibitem[Whi80]{White1980A-unique-exchan}
N.~White, \emph{A unique exchange property for bases}, Linear Algebra and its
  Applications \textbf{31} (1980), 81--91.

\bibitem[Woo04]{Woods:thesis}
K.~Woods, \emph{Rational generating functions and lattice point sets}, Ph.D.
  thesis, University of Michigan, 2004.

\end{thebibliography}
\end{document}